\documentclass{amsart}

\usepackage{amsthm}
\usepackage{amssymb}
\usepackage{mathtools}
\usepackage[shortlabels]{enumitem}
\usepackage{url}
\usepackage{bbm}
\usepackage[margin=3cm]{geometry}
\usepackage{hyperref} \hypersetup{colorlinks=true, citecolor=blue, linkcolor=black}

\usepackage{todonotes}
\DeclarePairedDelimiter\ceil{\lceil}{\rceil}
\DeclarePairedDelimiter\floor{\lfloor}{\rfloor}

\newcommand{\R}{\mathbb{R}}
\newcommand{\N}{\mathbb{N}}
\newcommand{\C}{\mathbb{C}}
\newcommand{\Z}{\mathbb{Z}}

\newcommand{\T}{\mathbb{T}}
\newcommand{\E}{\mathbb{E}}

\renewcommand{\Re}{\operatorname{Re}}
\renewcommand{\Im}{\operatorname{Im}}

\newcommand{\A}{\mathcal{A}}
\renewcommand{\d}{\mathrm{d}}
\renewcommand{\i}{\mathbf{i}}
\newcommand{\tr}{\operatorname{Tr}}

\renewcommand{\u}{\mathbf{U}}

\newcommand{\cst}{{\rm c}}
\newcommand{\Cst}{{\rm C}}

\newtheorem{thm}{Theorem}[section]
\newtheorem{prop}[thm]{Proposition}
\newtheorem{lemma}[thm]{Lemma}
\newtheorem{cor}[thm]{Corollary}
\numberwithin{equation}{section}
\newtheorem{remark}[thm]{Remark}

\title{From Berry-Esseen to super-exponential}

\date{\today}
\author{Klara Courteaut} 
\author{Kurt Johansson}
\author{Gaultier Lambert}

\begin{document}

\begin{abstract}
For any integer $m<n$, where $m$ can depend on $n$, we study the rate of convergence of $\frac{1}{\sqrt{m}}\tr \mathbf{U}^m$ to its limiting Gaussian as $n\to\infty$ for orthogonal, unitary and symplectic Haar distributed random matrices $\mathbf{U}$ of size $n$. In the unitary case, we prove that the total variation distance is less than $\Gamma(\floor{n/m}+2)^{-1}m^{-\floor{n/m}}\floor{n/m}^{1/4}\sqrt{\log n}$ times a constant. This result interpolates between the super-exponential bound obtained for fixed $m$ and the $1/n$ bound coming from the Berry-Esseen theorem applicable when $m\ge n$ by a result of Rains. We obtain analogous results for the orthogonal and symplectic groups. In these cases, our total variation upper bound takes the form $\Gamma(2\floor{n/m}+1)^{-1/2}m^{-\floor{n/m}+1}(\log n)^{1/4}$ times a constant and the result holds provided $n \geq 2m$.
For $m=1$, we obtain complementary lower bounds and precise asymptotics for the $L^2$-distances as $n\to\infty$, which show how sharp our results are. 
\end{abstract}

\maketitle
\tableofcontents

\section{Introduction and Results}

\subsection{Background and problem.} We consider random matrices from the orthogonal $\mathbb{O}(n)$, unitary $\mathbb{U}(n)$ and symplectic $\mathbb{S}\mathbbm{p}(2n)$ groups, distributed according to normalized Haar measure. This is the unique translation invariant probability measure on the group, i.e. $\mathbf{U}M \overset{d}{=} M\mathbf{U} \overset{d}{=} \mathbf{U}$ for any fixed $M$ in the same group as $\mathbf{U}$. 
In case of $\mathbb{U}(n)$, this is also known as the circular unitary ensemble or CUE. We refer to Meckes' monograph \cite{Meckes} for an introduction to random matrix theory for the classical compact groups. 



Our objects of interest are traces of powers of these matrices. 
It is well-known that these random variables are asymptotically Gaussian as a consequence of the celebrated strong Szeg\H{o} limit theorem and the Heine-Szeg\H{o} identity \eqref{HS}. 
The joint moments of these traces were first studied in \cite{DS} using representation theory, by analogy with random (uniform) permutations matrices\footnote{
Traces of powers of a permutation matrix are determined by the cycle structure of the corresponding permutation.}, and they are exactly equal to those of Gaussians, up to a surprisingly large order. 
This conveys that traces of fixed powers of random matrices from the compact classical groups converge extremely fast to their limiting Gaussians. To back this claim up, it is established that for a fixed $m\in\N$, the number of $m$-cycles of a uniform permutation converges superexponentially fast to their limiting Poisson distribution, \cite{AT}. 

The rate of convergence for traces of powers of Haar distributed random matrices was first investigated by Stein in \cite{Stein} who obtained a super-polynomial rate of convergence for a single trace of any fixed power of an orthogonal matrix.
Then, in \cite{97}, the second author proved that for any fixed (normalized) real-valued polynomial $P$, there exist $C,\delta>0$ depending on $P$ so that for a random matrix $\mathbf{U}\in \mathbb{U}(n)$, 
\begin{equation} \label{supregime}
\d_{\rm TV}\big( \tr P(\mathbf{U}) ,  \boldsymbol{\gamma}_\R \big) \le C n^{-\delta n}
\end{equation}
where $\d_{\rm TV}$ denotes the total variation distance and $ \boldsymbol{\gamma}_\R$  is a standard real-valued Gaussian. 
In the same article it was shown that for a random matrix $\mathbf{O}\in \mathbb{O}(n)$ or  $\mathbf{O}\in\mathbb{S}\mathbbm{p}(2n)$, 
\[
\d_{\rm TV}\big( \tr P(\mathbf{O}) ,  \boldsymbol{\gamma}_\R \big) \le C e^{-\delta n} .
\] 

The result for unitary matrices was then revisited and generalized to the multivariate case, i.e. considering jointly the traces of the first $m$ powers simultaneously, in \cite{JL}. 
The highest power $m$ was allowed to depend on the size $n$ of the matrix up to $m \ll \sqrt{n}$ while still having a fast rate of convergence. The multivariate super-exponential rate for orthogonal and symplectic matrices was obtained in \cite{CJ}, where $m$ was allowed to increase with $n$ up to $m \leq n^{1/4}$. 
Other results in the multivariate setting (with polynomial rate of convergence) which are based on Stein's method have been obtained in \cite{Fulman, FulmanRoss,DSt,Webb16}.
An analogous result has also been recently obtained in \cite{GR} for (uniform) random matrices from the unitary group over finite fields, in which case the rate of convergence is $O(e^{-\delta n^2})$.

Remarkably, the situation for very large powers is entirely different. Rains showed in  \cite{Rains} that for a random matrix $\mathbf{U}\in \mathbb{U}(n)$, the eigenvalues of $\mathbf{U}^n$ are statistically independent. 
Then, by the classical Berry-Esseen theorem, this implies that there exists a constant $C>0$ so that for $m\ge n$, 
\begin{equation} \label{BEregime}
\d_{\rm TV}\big( \tfrac{\tr \mathbf{U}^m}{\sqrt{n}} , \boldsymbol{\gamma}_\C \big) \le \tfrac{C}{n} 
\end{equation}
where $ \boldsymbol{\gamma}_\C$ is a standard complex-valued Gaussian. 
The rate $O(n^{-1})$ follows by symmetry and it is sharp. Analogous results hold for the other groups as well. 
The rationale behind this observation is that while the eigenvalues of $\mathbf{U}$ are rigid and evenly distributed on the unit circle \cite{Lam}, taking growing powers enforces an expansion which in effect cancels out the eigenvalue repulsion.
This phenomenon, as well as a precise description of the eigenvalues of $\mathbf{U}^m$, is described by the results of Rains, Theorem~\ref{Rains} below. 
In particular, the variance of the random variable $\tr \mathbf{U}^m$
 is exactly $\min(m,n)$ for $m\in\N$. 
Finally, let us also mention that a multivariate central limit theorem valid for linear combinations of arbitrary powers was obtained in \cite{DJ10} as a generalization of the strong Szeg\H{o} limit theorem.  

The questions that we address in this paper are twofold
\begin{itemize}
\item What is the exact dependency in $n$ of the rate of convergence for $\tr \mathbf{U}$?
\item How do the rates of convergence for $\tfrac{\tr \mathbf{U}^m}{\sqrt{m}}$ for $m$ with $m\le n$ interpolate between the super-exponential regime \eqref{supregime} and the Berry-Essen regime \eqref{BEregime}?
\end{itemize}

\subsection{Results for the unitary group.}
Let us first address theses questions precisely for a (Haar-distributed) random matrix 
$\mathbf{U}\in\mathbb{U}(n)$. We want to give bounds which hold even if $n$ is not large with explicit not too large constants. For $\tr\u$ we will also give a lower bound, which is close to the upper bound, on
the $L^1$-distance to a complex Gaussian. Furthermore in $L^2$-distance we give the precise asymptotics as a function of $n$ as $n\to\infty$.
Let $\mathbf{X}, \mathbf{Y}$ be two random variables taking values in $
\R^d$  with p.d.f.~$p_\mathbf{X}$ and $p_\mathbf{Y}$. Recall that the total variation distance between $\mathbf{X}$ and $\mathbf{Y}$ is defined as
\[
\d_{\rm TV}(\mathbf{X} ,\mathbf{Y}) = \sup_A \bigg(  \int_A p_\mathbf{X}  - \int_A p_\mathbf{Y}  \bigg) 
\]
where the sup is taken over all Borel sets $A\subset \R^d$. 
Then, we verify that 
\[
\d_{\rm TV}(\mathbf{X} ,\mathbf{Y}) 
\le 
 \| p_\mathbf{X}  - p_\mathbf{Y}  \|_{L^1} 
 \le 2 \d_{\rm TV}(\mathbf{X} ,\mathbf{Y}) .
\]
Let $\boldsymbol{\gamma}_\C$ denote a standard complex-valued Gaussian, with probability density function (p.d.f.) $\phi_\C : z \in \C \mapsto e^{-|z|^2}/{\pi}$.

\begin{thm} \label{thm:TV}
Let $p_n$ be the p.d.f.~of the complex-valued random variable $\tr\u$.
For any $n\ge 66$,
\begin{equation} \label{TV1}
\frac{1/300}{\Gamma(n+2)\sqrt{n+1}} \le \big\| p_n - \phi_\C \big\|_{L^1}  \le 19 \frac{n^{1/4}\sqrt{\log n}}{\Gamma(n+2)} . 
\end{equation}
\end{thm}

The proof of Theorem~\ref{thm:TV} is given in Section~\ref{sec:U}.
 It follows from similar estimates for $\| p_n- \boldsymbol{\gamma}_\C\|_{L^p}$ for $p=2,\infty$ which are obtained using Fourier analysis and the connection between the characteristic function of $\tr\u$ and certain Fredholm determinants.
This connection is reviewed in Section~\ref{sec:not}. 
In the course of the proof, we use \emph{Wolfram Mathematica} for numerical evaluations of several constants involved. Controlling the constants is the main reason behind the condition $n\ge 66$ and other similar conditions below.
We also comment that for $n=66$, the bound from Theorem~\ref{thm:TV} already implies that 
\[
\d_{\rm TV}(\tr\u ,\boldsymbol{\gamma}_\C) \le 4\cdot 10^{-93}
\]
which is far below Machine Epsilon (of order of $10^{-33}$ for quad(ruple) precision decimal). 

Our analysis also provides the asymptotics of the $L^2$-distance between $p_n$ and the p.d.f.~of a complex Gaussian. 

\begin{thm} \label{prop:asymp}
As $n\to\infty$,
\[
\big\| p_n - \phi_\C \big\|_{L^2}^2  \sim \frac{2 e^{4} \sqrt{\pi} }{\Gamma(n+2)^2\sqrt{n}}  . 
\]
\end{thm}

Using the result of Rains \cite{Rains}, we can also precisely estimate the total variation distance between the random variable $\frac{\tr\u^m}{\sqrt m}$ and a complex Gaussian for any power $m\in \N \cap [2,n]$.

 \begin{thm}\label{thm:inter}
 Let $n,m\in\N$ and $p_{n,m}$ be the p.d.f.~of the complex-valued random variable $\frac{\tr\u^m}{\sqrt m}$.
 Assume that $m\ge 2$ and that $n\ge \max\{700,m\}$, then with $N = \lfloor n/m \rfloor $, 
\[
 \big\| p_{n,m} - \phi_\C \big\|_{L^1} 
\le 12 \frac{(N+1)^{1/4}\sqrt{\log n}}{\Gamma(N+2) m^N} . 
 \]
\end{thm}

The proof of Theorem~\ref{thm:inter} is given in Section~\ref{sec:inter} and it relies on results used to prove Theorem~\ref{thm:TV} and Theorem~\ref{Rains}.
We emphasize that these estimates interpolate between the super-exponential rate of convergence when the degree $m$ is fixed, and the  polynomial rate of convergence \eqref{BEregime} when the ratio $N= \lfloor n/m \rfloor \in \N$ is fixed up to a $\sqrt{\log n}$ factor.

 \subsection{Results for the orthogonal and symplectic groups}
We obtain similar results for a Haar-distributed random matrix $\mathbf{O}$ from the orthogonal $\mathbb{O}(d)$ and symplectic $\mathbb{S}\mathbbm{p}(d)$ groups. Here we need to differentiate between $d$, the total number of eigenvalues and $n$, the number of non-trivial eigenvalues, i.e. the eigenvalues in the open upper half-plane. Complex eigenvalues come in conjugate pairs, and in the orthogonal case, depending on the sign of the determinant and the parity of $d$, there might be deterministic eigenvalues at $\pm 1$.

Let $\phi_\R : x \in \R \mapsto  e^{-x^2/2}/{\sqrt{2\pi}}$ be the p.d.f.~of a standard real-valued Gaussian random variable. 
Let also $q_n$ be the p.d.f.~of the real-valued random variable $\tr \mathbf{O}$, with $n$ being the number of non-trivial eigenvalues.
The upper and lower bounds on the total variation distance are given by

\begin{thm} \label{thm:TVO}
For any $n\ge 124$,
\[
\frac{1/60}{\sqrt{\Gamma(2n+1)}(2n)^{1/4}} \le \big\| q_n - \phi_\R \big\|_{L^1}  \le 5 \frac{(\log 2n)^{1/4}}{\sqrt{\Gamma(2n+1)}} . 
\]
\end{thm}

Again, we can give the precise asymptotic $L^2$-distance as $n\to\infty$.
 
\begin{thm} \label{prop:asympO}
As $n\to\infty$,
\[
\big\| q_n - \phi_\R \big\|_{L^2}^2  \sim \frac{e^{2}/\sqrt{2}}{\Gamma(d+1)\sqrt{d}}
\]
where $d$ is the total number of eigenvalues ($d=2n,\ 2n+1,\ \mathrm{or}\ 2n+2$ depending on its parity and on the sign of the determinant, see the beginning of Section \ref{sec:O}).
\end{thm} 

For traces of higher powers we have the following result which has somewhat stronger (technical) conditions compared to the unitary case. We are not able to go all the way up to $m=n$.
\begin{thm}\label{thm:interO}
Let $n,m\in\N$ and $q_{n,m}$ be the p.d.f.~of the real-valued random variable $\frac{1}{\sqrt m}\tr \mathbf{O}^m$ (with $n$ the number of non-trivial eigenvalues).
Set $N = \lfloor n/m \rfloor $. If $N$ and $m$ satisfy one of the following conditions:
\[ N \geq 5\ \mathrm{and}\ m\geq 66, \quad N \geq 4\ \mathrm{and} \ m\geq 129, \quad N \geq 2\ \mathrm{and}\ m\geq 10^4 \] 
then
\[
 \big\| q_{n,m} - \phi_\R \big\|_{L^1} 
\le 7 \frac{(\log n)^{1/4}}{\sqrt{\Gamma(2N+1)}m^{N-1}} . 
 \]
\end{thm}

The stronger requirements on $N$ and $m$ are ultimately a consequence of the fact that the Fredholm determinants arising from the characteristic function of $\tr \mathbf{O}$ are not bounded by one, unlike in the unitary case.

\subsection{Notation and ideas of the proof} \label{sec:not}
In this section, we explain the main ideas underlying the proofs of Theorem~\ref{thm:TV}. For simplicity, we focus on the case of the unitary group, the adaptation to the  orthogonal and symplectic groups are presented in Section~\ref{sec:O}. The method originates from our previous works \cite{97,JL,CJ}, but the fact that we are considering just the trace of 
a random unitary matrix means that we can do a considerably more precise asymptotic analysis and keep a very good control of the constants.

First, using Gaussian concentration bounds for $
\tr\mathbf{U}^m$ where $\mathbf{U}\in \mathbb{U}(n)$ is Haar-distributed, one can reduce the proof of Theorem~\ref{thm:TV} and~\ref{thm:inter} to controlling 
$ \big\| p_{n,m} - \phi_\C \big\|_{L^2}$; cf.~Section~\ref{sec:TV}. 

We define the characteristic function of the complex values random variable $\frac{\tr\mathbf{U}^m}{\sqrt m}$ by 
\begin{equation}
\begin{aligned}
F_{n,m}(\zeta) &= \E_n\big[e^{\i \tr \Re(\zeta\frac{\mathbf{U}^m}{\sqrt m} )} \big] \\
& = \int_\C e^{\i \Re(\zeta z)}   p_{n,m}(z) \d^2 z
\end{aligned}
\end{equation}
for $\zeta\in\C$, where $\d^2z$ denotes the Lebesgue measure on $\C$. 
Hence, $F_{n,m}(\zeta)$ is the Fourier transform of the p.d.f.~$p_{n,m}$ evaluated at $(\Re\zeta,-\Im\zeta)$ and by Plancherel's Theorem,
\begin{equation} \label{Planch}
\big\| F_{n,m} - \widehat{\phi_\C} \big\|_{L^2}^2 = 4\pi^2 \big\| p_{n,m} - \phi_\C \big\|_{L^2}^2 . 
\end{equation}
where $ \widehat{\phi_\C}$ denotes the Fourier transform of the (standard) complex  Gaussian p.d.f.,
\[
\widehat{\phi_\C}(\zeta)= \int e^{\i \Re(\zeta \overline{z})}  \phi_\C(z) \d^2z = e^{-|\zeta|^2/4} , \qquad \zeta\in \C.
\]

Consequently, the problem is to approximate the characteristic function $ F_{n,m}$. We now focus on the case $m=1$, letting $p_{n} = p_{n,1}$ and $F_{n} = F_{n,1}$. 
This is not a loss of generality since by Theorem~\ref{Rains}, one can express $ F_{n,m}(\zeta) =  \prod_{0\leq i < m}F_{N_i}\Big(\frac{\zeta}{\sqrt{m}}\Big)$ where $N_i \in \{ \lfloor n/m \rfloor, \lfloor n/m \rfloor +1\}$. 

\medskip

Let us denote by $\{e^{\i\theta_j}\}_{j=1}^n$ the eigenvalues of  $\mathbf{U}$. 
We will use the Heine-Szeg\H{o} identity: for any integrable function $\omega$ on the unit circle,
\begin{equation}\label{HS}
\E_n\bigg[\prod_{1\leq j\leq n} \omega(e^{\i\theta_j})\bigg] = \det(T_n(\omega))
\end{equation}
where $T_n(\omega)=(\hat{\omega}_{j-k})_{j,k=1}^n$ is a Toeplitz matrix\footnote{The Fourier coefficient of $\omega$ are defined by 
$\displaystyle \hat{\omega}_k = \int_\T e^{-\i k\theta} \omega(e^{\i\theta})  \frac{\d\theta}{2\pi}$ for $k\in\Z$ where $\T = \R/[2\pi]$.}.
This allows to write the Laplace transform of a general linear statistics
$\tr f(\mathbf{U}) = \sum_{j=1}^n f(e^{\i\theta_j})$ as a Toeplitz determinant;  for $\zeta\in\C$, 
\begin{equation}\label{heine-szego}
\E_n\big[e^{\i  \Re(\zeta \tr f(\mathbf{U}))}\big] = 
\det(T_n(\omega)) , \qquad \omega = e^{\i \Re(\zeta f)} . 
\end{equation}
In particular if $f(z) = z$ and $\zeta = r e^{\i \phi}$, then $\omega(e^{\i\theta}) = e^{\i r \cos(\theta+\phi)}$ and we have for $k\in\Z$, 
\[ \begin{aligned}
\widehat{\omega}_k 
&= \int_\T e^{-\i k\theta} \omega(e^{\i\theta})  \frac{\d\theta}{2\pi} = e^{\i k\phi}  \int_\T e^{-\i k\theta +\i r \cos \theta} \frac{\d\theta}{2\pi} \\ 
& = e^{\i k(\phi+\pi/2)} J_k(r)
\end{aligned}\]
where $(J_k)_{k\in\Z}$ are Bessel functions (of the first kind), cf.~DLMF formulae, \cite{DLMF}, \href{https://dlmf.nist.gov/10.9}{(10.9.2)} and \href{https://dlmf.nist.gov/10.2}{(10.2.2)}.  

Hence, by \eqref{heine-szego} and Hadamard's inequality for determinants,  we obtain an a priori bound 
\begin{equation} \label{Had}
\big|  F_n(\zeta) \big|^2  =  \big| \det(T_n(\omega)) \big|^2 \le  \prod_{j=1}^n \sum_{i=1}^n \big| \widehat{w}_{j-i}\big|^2  
= \prod_{j=1}^n \sum_{i=1}^n \big|J_{j-i}(|\zeta|)\big|^2 . 
\end{equation}
This bound will be useful to control the tail of the characteristic function $F_n$; cf.~Proposition~\ref{lem:Had1}. 

To obtain the exact asymptotics of $F_n(\zeta)$ if $\zeta$ is not too large, we use another and perhaps not as well known exact formula, known as the Borodin-Okounkov formula (sometimes also known as Geronimo-Case formula).
This formula first appeared in \cite{GC79}, and then in \cite{BO00,BW00,Bottcher} with different proofs.

\begin{thm}\label{borodin-okounkov}
Assume that $f\in L^\infty (\T)$ is complex-valued and satisfies $\sum_{k\in\Z} \lvert k\rvert \lvert \hat{f}_k \rvert^2 < \infty$ and $\hat{f}_0=0$. Let $\omega = e^f$ and 
\begin{equation}\label{mathcalA}
\mathcal{A}= \sum_{k\geq 1} k \hat{f}_{k}\hat{f}_{-k} .
\end{equation}
Then, there exists a trace class operator $K$ (depending on $\omega$) on $l^2(\N)$ such that for any $n\in\N$, 
\[
\E_n\bigg[\prod_{1\leq j\leq n} \omega(e^{\i\theta_j})\bigg] = 
 e^{\mathcal{A}} \det(\operatorname{I} - KQ_n)
\]
where the right-hand side is a Fredholm determinant and $Q_n$ denotes the orthogonal projection with kernel $\operatorname{span}(\mathrm{e}_1,\dots,\mathrm{e}_{n})$ on $l^2(\N)$. 
\end{thm}

Moreover, the operator $K$ admits an explicit representation in terms of Hankel operators \cite{BW00}.
For $\omega \in H^{1/2}(\T \rightarrow \C)$, define
\[ 
H_+(\omega) = (\widehat{\omega}_{i+j-1})_{i,j\geq 1} 
\qquad\text{and}\qquad 
H_-(\omega) = (\widehat{\omega}_{-i-j+1})_{i,j\geq 1}.
\]
The condition $\|\omega \|_{H^{1/2}} = \sqrt{ \sum_{k\in\Z} \lvert k\rvert \lvert \widehat{\omega}_k \rvert^2} < \infty$ guarantees that these (infinite) matrices define Hilbert-Schmidt operators on  $l^2(\N)$. 
In particular, if $\sum_{k\in\Z} \lvert k\rvert \lvert \hat{f}_k \rvert^2 < \infty$ as in Theorem~\ref{borodin-okounkov},  
then we can write the operator 
\begin{equation} \label{BOop} 
K = H_+(\omega_+) H_-(\omega_-) , \qquad 
\omega_{\pm}(e^{\i\theta}) = e^{ \sum_{k>0}\big( \pm \hat{f}_{-k} e^{-\i k \theta} \mp \hat{f}_k e^{\i k \theta} \big) } . 
\end{equation}
The condition $f\in H^{1/2}$ also guarantees that  $\omega_\pm \in H^{1/2}$, so that $H_{\pm}(\omega_\pm)$ are Hilbert-Schmidt and $K$ is a trace-class operator.  
That being said, we will only  apply Theorem~\ref{borodin-okounkov} in the case  $\omega(z) = e^{\i \Re(\zeta z)}$. 
Then, 
\[
\omega_{\pm}(e^{\i \theta}) = e^{ \i \big( \pm \overline{\zeta} e^{-\i \theta}  \mp  \zeta e^{\i \theta} \big)/2} = e^{\pm \Im(\zeta e^{\i\theta}) }  
\qquad\text{and}\qquad
\mathcal{A}  = - |\zeta|^2/4 . 
\]

This implies that for $\zeta = r e^{\i \phi} \in \C$ and  $k\in\N$,
\begin{equation*} 
 \begin{aligned}
(\widehat{\omega_\pm})_{\pm k} 
&=  \int_\T e^{\mp \i k\theta}  e^{\pm \Im(\zeta e^{\i\theta}) } \frac{\d\theta}{2\pi} 
 =  e^{\pm\i k \phi} \int_\T e^{\mp \i k\theta \pm r \sin \theta}  \frac{\d\theta}{2\pi}  \\
 &=e^{\pm\i k(\phi \pm \pi/2) } I_k(r),
\end{aligned}
\end{equation*}
where $(I_k)_{k\in\N}$ are modified Bessel functions (of the first kind), cf.~DLMF formulae \href{https://dlmf.nist.gov/10.32}{(10.32.3)} and \href{https://dlmf.nist.gov/10.25}{(10.25.2)}.  

\medskip

Let $K$ be the operator coming from Theorem~\ref{borodin-okounkov} associated with the symbol  $\omega(z) = e^{\i \Re(\zeta z)}$ for $\zeta = r e^{\i\phi}$.  
Let us also introduce the (infinite) matrices $\mathbb{J}(r) = \big(I_{i+j-1}(r)\big)_{i,j\geq 1}$ where $I_k(r)$ are Bessel functions evaluated at $r>0$ and   
$\Delta_{\pm}(\phi) = \operatorname{diag}(e^{\pm\i (k-1/2)(\phi \pm \pi/2)})_{k\ge 1}$ for $\phi \in [0,2\pi)$ so that we can write 
\[
H_\pm(\omega_\pm)  =  \Delta_{\pm}(\phi)   \mathbb{J}(r)  \Delta_{\pm}(\phi) .
\]
Hence, according to \eqref{BOop}, we have
\[
\Delta_+^{-1}(\phi) K\Delta_+(\phi) =  \mathbb{J}(r) \mathbb{L}   \mathbb{J}(r) \mathbb{L} ,\qquad 
\mathbb{L} =  \Delta_-(\phi)\Delta_+(\phi) = \Delta_+(\phi)\Delta_-(\phi) . 
\]
In particular, the matrix  $\mathbb{L} =  \operatorname{diag}(e^{\i k \pi})_{k\ge 1}$ is indeed independent of the argument $\phi$ of the parameter $\zeta\in\C$.
 By Theorem~\ref{borodin-okounkov}, this implies that for $n\in\N$, $\zeta\in\C$ with $|\zeta|=r$, 
\begin{equation}\label{charfunction}
 \begin{aligned}
  F_n(\zeta) & =  e^{-|\zeta|^2/4} \det(\operatorname{I} - KQ_n) \\
  & = e^{-r^2/4} \det(\operatorname{I} - Q_n K(r)Q_n) ; \qquad
  K(r) =  \mathbb{J}(r) \mathbb{L}   \mathbb{J}(r) \mathbb{L} , 
\end{aligned}
\end{equation}
where the matrices $\mathbb{J}(r) = \big(I_{i+j-1}(r)\big)_{i,j\geq 1}$. 
Note that \eqref{charfunction} follows from the cyclicity of Fredholm determinants and the fact that the projection $Q_n$ commutes with $\Delta_{\pm}(\phi)$.

\medskip

Formulae \eqref{Had} and \eqref{charfunction} are the starting point of our analysis. They  provide explicit formulae for the characteristic function of the random variable $\tr\mathbf{U}$ in terms of (modified) Bessel functions which have well-known properties and asymptotics; cf.~DLMF \href{https://dlmf.nist.gov/10}{Section 10}. 
These functions can be defined as series: for any $\nu\ge 0$ and $z\in\C$, 
\begin{equation} \label{BesselI}
I_{\nu}(z) =  \sum_{j =0}^{+\infty}  \frac{(z/2)^{2j+\nu}}{j!\Gamma(j+\nu+1)} .
\end{equation}
Then, for $k\in\Z$, $I_{-k}(z) = I_k(z) $ and  
\[
J_k(z) = e^{-\i k \pi/2} I_{k}(\i z) ;
\]
see DLMF formulae \href{https://dlmf.nist.gov/10.27}{(10.27.1)} and  \href{https://dlmf.nist.gov/10.27}{(10.27.6)}.  

Our general strategy is the following, by \eqref{charfunction} and going to polar coordinates,
\begin{equation} \label{BOF}
\big\| F_{n} - \widehat{\phi_\C} \big\|_{L^2}^2
=2 \pi \int_0^\infty \big( \det(\operatorname{I} - K(r) Q_n) -1\big)^2 \d( e^{-r^2/2}) .
\end{equation}
We show in Section~\ref{sec:trace} that if $r<2\cst n$ for a small enough $\cst>0$, the trace-norm
\[ 
\| Q_nK(r)Q_n \|_{J_1} \ll 1 ,
\] 
so the Fredholm determinant on the right-hand side of \eqref{BOF} is close to 1. 
Moreover, in this regime, we can approximate
\[
 \det(\operatorname{I} -  K(r)Q_n)
 \simeq \exp\big(- \tr(K(r)  Q_n) \big) . 
\]
By definitions, 
\begin{equation}  \label{tracef}
\tr(K(r)  Q_n) = \tr \big( Q_n \mathbb{J}(r) \mathbb{L}   \mathbb{J}(r) \mathbb{L} Q_n\big)
= \sum_{i >  n ,  j\ge 1} (-1)^{i+j} I_{i+j-1}^2(r) . 
\end{equation}
Then, we show that $\tr(K(r)  Q_n) \simeq (-1)^{n+1} I^2_n(r)$ for $r< 2\cst n$, so that 
\[
\big( \det(\operatorname{I} - K(r) Q_n) -1\big)^2 \simeq I_{n+1}^4(r)
\]
and by \eqref{BOF}
\[
\big\| F_{n} - \widehat{\phi_\C} \big\|_{L^2}^2 
\simeq 2 \pi \int_0^{2\cst n}  I_{n+1}^4(r) \d( e^{-r^2/2}) .  
\]
The asymptotics of this integral are performed in Section~\ref{sec:asymp} and this leads to the proof of Theorem~\ref{prop:asymp}. 

In the complementary regime, we show that if $n$ is large enough
\begin{equation} \label{TB}
\int_{|\zeta| \ge 2\cst n} | F_{n}(\zeta) |^2 \d^2\zeta \le \exp( - \delta n^2)
\end{equation}
for a small $\delta>0$. 
In order to obtain \eqref{TB}, we use the simple inequality \eqref{Had} in the regime $|\zeta| \gg e^{\delta n}$ which we combine with a Gaussian tail bound of the form $ | F_{n}(\zeta) |^2 \ll \exp( - \delta n^2)$ valid for all $|\zeta| \ge 2\cst n$; cf.~Proposition~\ref{prop:TB1}. 
The proof of this proposition relies on the change of variable technique, or \emph{loop equation}, originally introduced in \cite{Joh}, see also \cite{Lam}. We review the method in Section~\ref{sec:tb}.
When applied to the linear statistic $\tr\mathbf{U}$, the key new ingredient is that the error terms involve again modified Bessel functions and can be explicitly controlled.

\section{Unitary group: Proof of Theorems \ref{thm:TV} and \ref{prop:asymp}} 
\label{sec:U}

The main steps of the proof of Theorem \ref{thm:TV} are to obtain asymptotics for the characteristic function $F_n(\zeta)$ of the random variable $\tr\u$ in different regimes of $\zeta\in\C$.
By formula \eqref{charfunction}, this characteristic function is a Fredholm determinant and we argue that if the parameter $|\zeta| \le \cst n$ for a given constant $\cst>0$  where $n$ is the dimension of $\u$, this determinant almost equals $1$ and the error is controlled by  $\tr (KQ_n)$. 
Then, we obtain precise estimates for the decay of this trace, cf.~Section~\ref{sec:trace}. The regime $|\zeta| \ge \cst n$ is controlled by different methods; cf.~Section~\ref{sec:ub}. 
This leads to the following sharp bounds. 

\begin{thm} \label{thm:1}
Let $p_n$ be the density function of the random variable $\tr \u$ and  $\cst$ be the (unique) solution of the equation $\cst e^{1+\cst^2} =1$.  Then for any $n\geq 13$,
\begin{equation} \label{ub0}
\frac{1}{\Gamma(n+2)^2\sqrt{n+1}}  \le 4\pi^2 \| p_n - \phi_\C \|_{L^2}^2 \leq  
\frac{28 \sqrt\pi }{(1-\cst^4)^4}  \frac{ e^{4+ 12/n} }{\sqrt{n+1}}  \frac{1}{\Gamma(n+2)^2} +  7 n^2 e^{-\frac{n(n-2)}{9}} .
\end{equation}
The Gaussian term becomes negligible if $n\ge 63$,
\begin{equation} \label{L2sup}
\| p_n - \phi_\C \|_{L^2}^2 \leq  \frac{87}{\sqrt{n}} \frac{1}{\Gamma(n+2)^2} 
\end{equation}
and for $n\ge 66$, 
\begin{equation}\label{Linftysup}
\| p_n - \phi_\C  \|_{L^\infty}
\le  \frac{43}{\Gamma(n+2)} .  
\end{equation}
\end{thm}

We first deduce the total variation bound from Theorem \ref{thm:TV}, and then present the proof of Theorem~\ref{thm:1} in the next sections. 
A variation of the argument provides the exact asymptotics of the $L^2$-norm between the density function $p_n$ of the random variable $\tr \u$ and the Gaussian density. Namely, we prove Theorem~\ref{prop:asymp} in Section~\ref{sec:asymp}. 

\subsection{Total variation  upper--bound; Proof of Theorem~\ref{thm:TV}}
\label{sec:TV}

We need the following concentration inequality from \cite{JL}.

\begin{lemma}\label{concetrationU}
For any $L>0$ and any $m\in\N$,
\[ \mathbb{P}_n\big[ (\Re \tr \u^m, \Im \tr \u^m) \notin [-\tfrac{L}{2},\tfrac{L}{2}]^2 \big] < 4e^{-L^2/(4m)}. \]
\end{lemma}

Recall that $p_{n,m}$ denotes the p.d.f.~of the complex-valued random variable $\frac{\tr\u^m}{\sqrt m}$.
Then, by the Cauchy-Schwarz inequality and the triangle inequality, for any $m\in\N$, 
\[ \begin{aligned} 
 \big\| p_{n,m}-  \phi_\C \big\|_{L^1} 
&\le L  \big\| p_{n,m} -  \phi_\C \big\|_{L^2}  +  \mathbb{P}_n[ \tfrac{\tr \u^m}{\sqrt{m}} \notin [-\tfrac{L}{2},\tfrac{L}{2}]^2 ]  + \mathbb{P}[ \boldsymbol{\gamma}_\C \notin [-\tfrac{L}{2},\tfrac{L}{2}]^2 ].
\end{aligned}\]
A change of variables to polar coordinates gives 
\begin{equation}\label{gaussiantail}
\mathbb{P}[ \boldsymbol{\gamma}_\C \notin [-\tfrac{L}{2},\tfrac{L}{2}]^2 ] \leq 2 \int_{L/2}^\infty re^{-r^2}\d r = 2e^{-L^2/4} .
\end{equation}
This is to be compared with Lemma \ref{concetrationU}, 
\[  
\mathbb{P}_n[ \tfrac{\tr \u^m}{\sqrt{m}} \notin [-\tfrac{L}{2},\tfrac{L}{2}]^2 ]   \leq 4e^{-L^2/4}
\]
so that for any $n,m\in\N$ and $L>0$,
\begin{equation} \label{CS}
 \| p_{n,m} - \phi_\C \|_{L^1}
\le L  \big\| p_{n,m} -  \phi_\C \big\|_{L^2}  + 6e^{-L^2/4}. 
\end{equation}
Using the upper-bounds \eqref{L2sup} and \eqref{CS} with $m=1$, this implies that  for any $n\ge 63$, 
\[
 \| p_n - \phi_\C \|_{L^1} \le 
2 \bigg( \frac{9.4 L }{n^{1/4}\Gamma(n+2)} + 3 e^{-L^2} \bigg),
\]
which we minimize by choosing $L =  \sqrt{ \log \Gamma(n+2)}$.
Since $ \log \Gamma(n+2) \le n\log n$ for $n\ge 3$, 
 we  conclude that
\[
 \| p_n - \phi_\C \|_{L^1} \le  
2 \left( 9.4+\frac{3}{n^{1/4}\sqrt{\log n}} \right) \frac{n^{1/4}\sqrt{\log n}}{\Gamma(n+2)} .
\]
This completes the proof of the upper-bound \eqref{TV1}.

\medskip

For the lower-bound, observe that by H\"older's inequality and \eqref{Linftysup}, we have for $n\ge 66$,
\[
\| p_n - \phi_\C \|_{L^2}^2 \le    
\| p_n - \phi_\C \|_{L^\infty}
\| p_n - \phi_\C \|_{L^1} 
\le  \frac{300}{\Gamma(n+2)}  \| p_n - \phi_\C \|_{L^1} . 
\]
Hence, using the lower-bound \eqref{ub0} from Theorem~\ref{thm:1}, we conclude that for $n\ge 66$,
\[
 \| p_n - \phi_\C \|_{L^1}  \ge \frac{1}{300\Gamma(n+2)\sqrt{n+1}}  . 
\]


\subsection{Proof of the upper-bounds in Theorem~\ref{thm:1}}
\label{sec:ub}
By Plancherel's Theorem, it suffices to control $\| F_n - \widehat{\phi_\C} \|_{L^2}^2$ where $F_n$ is the Fourier transform of the p.d.f.~$p_n$ (i.e.~the characteristic function of the random variable $\tr\mathbf{U}$). 
Then, we divide this integral in three parts depending on the argument;
\[
|\zeta| \le 2\cst n , \qquad 2\cst n \le |\zeta| \le \Lambda , \qquad  |\zeta| \ge \Lambda , 
\]
where $\cst>0$ is the (unique) solution of the equation $\cst e^{1+\cst^2} =1$ and $\Lambda \ge 7n/4$.
In the first regime, we show that by estimating the Fredholm determinant on the right-hand side of \eqref{charfunction}, we obtain the estimate
\begin{equation} \label{ub1}
\int_{|\zeta| \le 2\cst n} \big| F_n(\zeta) - e^{-|\zeta|^2/4}   \big|^2 \d^2\zeta  
\le  \frac{28 \sqrt\pi }{(1-\cst^4)^4}  \frac{ e^{4+ 7/n} }{\sqrt{n+1}}  \frac{1}{\Gamma(n+2)^2}  ; 
\end{equation}
cf.~Sections~\ref{sec:trace} and \ref{firstregime}.
In Section~\ref{sec:tb}, we obtain a more general estimate for the tail of the characteristic function $F_n$ which is useful in the second regime.
\begin{prop} \label{prop:TB1}
For any $n\in\N$ and $\zeta \in \C$, 
\[
| F_{n}(\zeta) | \le \sqrt{2}  \exp\left(-  \frac{n^2}{4(2.05 + n^2/ |\zeta|^2)} \right) . 
\] 
\end{prop}

Observe that we immediately deduce from Proposition~\ref{prop:TB1} that for any $n\in\N$ and $|\zeta| \ge 2 \cst n$, 
\begin{equation} \label{ub2}
| F_{n}(\zeta) | \le \sqrt 2 e^{- n^2/18}  , \qquad\text{since}\quad  \frac 1{18} \le \frac{1}{4(2.05+ \cst^{-2}/4)} 
\end{equation}

For the last regime, with $\Lambda = e^{O(n)}$, we use the following complementary bound which is proved in Section~\ref{sec:Had}. 
\begin{prop} \label{lem:Had1}
For any $n\in\N$ and  $|\zeta| \ge 7n/4$, 
\begin{equation} \label{ub3}
\big|  F_n(\zeta) \big|^2   \le  \left(\frac{n}{|\zeta|}\right)^{n} . 
\end{equation}
\end{prop}

The proof of Theorem~\ref{thm:1} follows easily from combining these estimates. 
First if we integrate the estimate \eqref{ub3}, we obtain for any $n\ge 3$ and $\Lambda \ge 7n/4 $, 
\begin{equation} \label{ub5}
\int_{|\zeta| \ge \Lambda}\big|  F_n(\zeta) \big|^2 \d^2\zeta  \le 2\pi n^n  \int_{\Lambda}^{+\infty} \frac{\d r}{r^{n-1}} =  \frac{2\pi}{n-2} \frac{n^n}{\Lambda^{n-2}} .
\end{equation}

Then, we want to minimize over all $\Lambda \ge 7 n /4$, the function 
\[
\Lambda \mapsto \frac{1}{n-2} \frac{n^n}{\Lambda^{n-2}} + \Lambda^2 e^{-   n^{2}/9}  . 
\]
There is a unique critical point $\Lambda_* = 2^{-1/n}ne^{n/9}$ and it is larger than $7n/4$ for $n\geq 8$. The minimum is given by 
\[
2^{-2/n}\frac{n^3}{n-2} e^{- 2 n(n-2)/9},
\]
which by \eqref{ub2} and \eqref{ub5} implies that for $n\ge 8$, 
\[
\int_{ |\zeta| \ge 2\cst n }
\big| F_n(\zeta) \big|^2\d^2\zeta  \le  \frac{ 2\pi n^3}{n-2}  \exp\left(  -  n(n-2)/9\right) .  
\]
Also, for the Gaussian we have the estimate
\begin{equation*}
\int_{ |\zeta|\geq 2\cst n } \hspace{-.3cm}  e^{-|\zeta|^2/2}  \d^2\zeta=  \pi \int_{4 \cst^2n^2}^{+\infty} e^{-u/2} \d u
= \frac{\pi}{2} e^{- 2\cst^2n^2} .
\end{equation*}

By combining these bounds, we obtain for $n\ge 10$
\begin{equation}\label{ub4}
\int_{ |\zeta| \ge 2\cst n }
\big| F_n(\zeta) -  e^{-|\zeta|^2/4} \big|^2\d^2\zeta
\le  \frac{\pi n^2}{2}\big(\frac{1}{100}+ 5\big) \exp\left(  -  n(n-2)/9\right) 
\end{equation}

Finally, combining the estimates \eqref{ub1} and \eqref{ub4}, using \eqref{Planch} with $m=1$, this completes the proof of the upper-bound \eqref{ub0}.

The upper-bound \eqref{Linftysup} follows from an analogous argument using that  for $n\ge 2$
\begin{equation} \label{ubinfty} 
\int_{|\zeta| < 2\cst n} \big| F_n(\zeta) - e^{-|\zeta|^2/4}   \big| \d^2\zeta 
 \le \frac{5.2\, e^{2+5/n}}{(1-\cst^4)^2} \frac{2\pi}{\Gamma(n+2)}   . 
\end{equation}
If we integrate the estimates  \eqref{ub2}--\eqref{ub3}, we obtain for $n\ge 5$  and any $\Lambda \ge 7n/4$, 
\[\begin{aligned} 
\int_{|\zeta| \ge  2\cst n}
\big| F_n(\zeta) \big|\d^2\zeta   & \le 2\pi \Lambda^2 e^{- n^{2}/18} ,  \\
\int_{|\zeta| \ge \Lambda}\big|  F_n(\zeta) \big| \d^2\zeta  & \le
\frac{2 \pi }{n/2-2} \frac{n^{n/2}}{\Lambda^{n/2-2}} . 
\end{aligned}\]

In this case,  we want to minimize over all $\Lambda \ge 7n/4$, the function 
\[
\Lambda \mapsto \frac{1}{n/2-2} \frac{n^{n/2}}{\Lambda^{n/2-2}} + \Lambda^2 e^{-  c_0 n^{2}}  . 
\]
There is a unique critical point $\Lambda_* = 2^{-\frac2n}ne^{2c_{0} n}$ and it is larger than $7n/4$ for $n\geq 8$ with $c_0=1/18$. The minimum is given by 
\[
2^{-\frac4n} \frac{n^3}{n-4} e^{- c_0 n(n-4)} .
\]
This argument shows that for $n\ge 12$, 
\[
\int_{ 2\cst n \le |\zeta|}
\big| F_n(\zeta) \big|\d^2\zeta  
\le  3\pi n^2 e^{- \frac{n(n-4)}{18}} .
\]
In turn, for $n\ge 12$, 
\begin{equation*}
\int_{ |\zeta| \ge 2\cst n }
\big| F_n(\zeta) -  e^{-|\zeta|^2/4} \big|\d^2\zeta
\le 4\pi n^2 \exp\left(  - \frac{n(n-4)}{18}\right)  . 
\end{equation*}

This shows that for $n\ge 66$, 
\begin{equation*}
\big\| F_n -  \widehat{\phi_\C} \big\|_{L^1}
\le  \frac{43\cdot 2\pi}{\Gamma(n+2)}   
\end{equation*}

By Fourier's inversion formula, this implies that  for $n\ge 66$, 
\[
\| p_n - \phi_\C \|_{L^\infty} \le  \frac{43}{\Gamma(n+2)}   
\]
This proves the estimate \eqref{Linftysup}.

\subsection{Trace-norm and Fredholm determinant estimates}
\label{sec:trace}

The trace-class norm and the Hilbert-Schmidt norm of a compact operator $K$ are defined by
\begin{equation}
\| K\|_{J_1} \coloneqq \sum_{j\geq 1} s_j(K), \qquad \| K\|_{J_2} \coloneqq \Big( \sum_{j\geq 1} s_j(K)^2 \Big)^{1/2}
\end{equation}
where $\{s_j(K)\}_{j\geq1}$ are the singular values. If $A$ and $B$ are Hilbert-Schmidt (i.e. $\| A\|_{J_2}<\infty$ and $\| B\|_{J_2}<\infty$), then $AB$ is trace-class and satisfies
\begin{equation}\label{holder}
\| AB\|_{J_1} \leq \| A\|_{J_2}\| B\|_{J_2},
\end{equation}
see e.g. Lemma 2.2. in \cite{GGK}. We will also use that if $K$ is given by an infinite matrix with elements $(K_{ij})_{i,j\ge 1}$, then
\begin{equation}\label{hilbert-schmidt}
\| K\|_{J_2} = \Big( \sum_{i,j\geq 1}  \lvert K_{ij} \rvert ^2 \Big)^{1/2}.
\end{equation}
For a short proof, let $\{e_j\}_{j\geq1}$ be the standard orthonormal basis of $l_2(\N)$. Then
\[ \sum_{j\geq 1} s_j(K)^2 = \sum_{j\geq 1} \langle K^*Ke_j,e_j \rangle = \sum_{j\geq 1} \| Ke_j \|^2 = \sum_{j\geq 1}\sum_{k\geq1} \lvert\langle Ke_j,e_k\rangle\rvert^2. \]

Let us recall that according to \eqref{charfunction}: for any $n\in\N$ and $\zeta\in\C$,
\[
F_n(\zeta) =  e^{-|\zeta|^2/4} \det(\operatorname{I} - K(r)Q_n)  , \qquad
 K(\zeta) = \mathbb{J}(r) \mathbb{L}   \mathbb{J}(r) \mathbb{L}
 \]
where $\mathbb{L} =  \operatorname{diag}(e^{\i k \pi})_{k\ge 1}$, $\mathbb{J}(r) = \big(I_{i+j-1}(r)\big)_{i,j\geq 1}$  and $r=|\zeta|$. 
Hence, this Fredholm determinant depends only on the properties of modified Bessel functions  \eqref{BesselI}. 
In particular, we will make use of the following bounds in the sequel, 
\begin{equation}\label{GammaBd}
\Gamma(k+1)(k+1)^j \leq \Gamma(j+k+1) \leq \Gamma(k+1)(k+j)^j, \quad j,k\in\N
\end{equation}
Then,  it holds for all $z\in \C$ with $|z|= r$, 
\begin{equation}\label{Bessel-upperbd}
\lvert I_k(z)\rvert \leq I_k(r)
\leq \frac{(r/2)^k}{\Gamma(k+1)} \sum_{j=0}^\infty \frac{(r/2)^{2j}}{\Gamma(j+1)(k+1)^j} = \frac{|z/2|^k}{\Gamma(k+1)}e^{|z|^2/ 4(k+1)}, \qquad k\in\N. 
\end{equation}
On the other hand, by keeping only the first term of the sum \eqref{BesselI}, we have 
\begin{equation}\label{Bessel-lowerbd}
I_k(r) \geq \frac{(r/2)^k}{\Gamma(k+1)}, \qquad k\in\N.
\end{equation}
Finally, we will also use make repeated use of Stirling's approximation for the Gamma function. Equation \href{https://dlmf.nist.gov/5.6}{(5.6.1)} in DLMF gives
\begin{equation}\label{Gamma}
\sqrt{2\pi} \nu^{\nu+1/2}e^{-\nu} < \Gamma(\nu+1) < \sqrt{2\pi e} \nu^{\nu+1/2}e^{-\nu}, \quad \nu\geq 1 .
\end{equation}

\begin{lemma}\label{traceclassnorm1}
Let $\cst$ be the (unique) solution of the equation $\cst e^{1+\cst^2} =1$ and let $\Cst_0^{-1} = 2\pi(1-\cst^2)^2$. 
Numerically, we have $\Cst_0 \le .17$ and $\cst\simeq 0.33$ . 
If $r<2\cst (n+1)$, then for $n\in\N$, 
\[ \| Q_nK(r)Q_n \|_{J_1} < \frac{\Cst_0}{n+1}. \] 
\end{lemma}

\begin{proof}
By \eqref{holder}, observe that 
\[ 
\| Q_nK(r)Q_n \|_{J_1} \le 
\|Q_n \mathbb{J}(r) \|_{J_2} \| \mathbb{L} \mathbb{J}(r) \mathbb{L}Q_n\|_{J_2}  .
\] 
Then, as $\mathbb{L}  = \mathbb{L}^*$ is unitary, $\|Q_n \mathbb{J} \|_{J_2} = \| \mathbb{L} \mathbb{J} \mathbb{L}Q_n\|_{J_2}$ by \eqref{hilbert-schmidt} and 
\[
\|Q_n \mathbb{J}(r) \|_{J_2}^2 =  \sum_{j\geq 0} \sum_{k\geq 0} |I_{j+k+n+1}(r)| ^2   
\]
Hence, using the upper-bound \eqref{Bessel-upperbd}, we obtain
\[ 
\|Q_n \mathbb{J}(r) \|_{J_2}^2  \le 
\sum_{j\geq 0} \sum_{k\geq 0} e^{r^2/2(j+k+n+2)} \frac{(r/2)^{2k+2j+2n+2}}{\Gamma(j+k+n+2)^2}  
\]
Inserting the lower-bound \eqref{GammaBd} gives
\begin{equation} \label{HSnorm}
\begin{aligned}
\|Q_n \mathbb{J}(r) \|_{J_2}^2  
&\leq e^{r^2/(2n+4)} \frac{(r/2)^{2n+2}}{\Gamma(n+2)^2} \sum_{j\geq 0} \sum_{k\geq 0} \frac{(r/2)^{2k+2j}}{(n+2)^{2k+2j}}  \\ \nonumber
&\leq \frac{e^{r^2/(2n+4)}}{(1-r^2/4(n+1)^2)^2} \frac{(r/2)^{2n+2}}{\Gamma(n+2)^2} \\ 
&\leq \frac{1}{2\pi(n+1)(1-r^2/4(n+1)^2)^2}  \bigg(\frac{r e^{1+r^2/(2n+2)^2}}{2n+2}  \bigg)^{2n+2}
\end{aligned}
\end{equation}
where we used \eqref{Gamma} to obtain the last bound.
Hence, if $r \le 2 \cst (n+1)$, since the right-hand side is increasing as a function of $r$, this implies that for any $n\in\N$, 
\[
\| Q_nK(\zeta)Q_n \|_{J_1} 
\leq \frac{1}{2\pi(n+1)(1-\cst^2)^2} = \frac{\Cst_0}{n+1}  . 
\qedhere
\]
\end{proof}

\begin{lemma}\label{trace}
For any $n\geq 1$ and $ |\zeta|< \sqrt{2} (n+1)$, 
\[ \Big(1-\frac{\lvert\zeta\rvert^2}{2(n+1)^2}\Big) \frac{\lvert \zeta/2 \rvert^{2n+2}}{\Gamma(n+2)^2} \leq | \tr Q_nK(|\zeta|)Q_n | \leq 2\frac{\exp\big(\frac{\lvert\zeta\rvert^2}{2(n+1)}\big)}{\big(1-(\frac{\lvert\zeta/2\rvert}{n+1})^4\big)^2}  \frac{\lvert \zeta/2 \rvert^{2n+2}}{\Gamma(n+2)^2}. \]
\end{lemma}

\begin{proof}
Set $r=\lvert\zeta\rvert$. By formula \eqref{tracef}, 
\begin{align}
\tr(K(r)  Q_n)
&= \sum_{j,k \geq 0} (-1)^{j+k+n+1}I_{j+k+n+1}^2(r) \notag \\
&= \sum_{\ell\in\N} (-1)^{\ell+n} \ell I_{\ell+n}^2(r) \notag \\
\label{traceexact}
&= (-1)^{n+1} \sum_{\ell \ge 0} \big((2\ell+1) I_{2\ell+n+1}^2(r)-(2\ell+2)I_{2\ell+n+2}^2(r)\big)
\end{align}
The modified Bessel function has the following integral representation (see \href{https://dlmf.nist.gov/10.32}{(10.32.2)} in DLMF):
\begin{equation}
I_\nu(x) = \frac{(x/2)^\nu}{\sqrt{\pi}\Gamma(\nu+1/2)} \int_{-1}^1 (1-t^2)^{\nu-1/2} e^{xt} \mathrm{d}t, \quad \Re \nu > -\frac{1}{2}
\end{equation}
whence, for $r>0$ and $k\in\N_0$,
\begin{equation} \label{Ideac}
I_{k+2}(r) < \frac{r}{2k+3} I_{k+1}(r) .
\end{equation}
Therefore, if $r< \sqrt{2} (n+1)$, all terms of the sum on the right-hand side of \eqref{traceexact} are positive and we obtain the lower-bound, 
\begin{equation} \label{trlb}
| \tr Q_nK(r) Q_n |  > \Big(1-\frac{r^2}{2(n+1)^2}\Big) \sum_{l\ge 0} (2l+1)I_{2l+ n+1}^2(r)
\end{equation}
By keeping only the first term $(l=0)$ and using (\ref{Bessel-lowerbd}), we conclude that
\begin{align*}
| \tr Q_n K(r) Q_n | 
> \Big(1-\frac{r^2}{2(n+1)^2}\Big) \frac{(r/2)^{2n+2}}{\Gamma(n+2)^2} .
\end{align*}

 For the upper bound, we simply use that
\[ | \tr Q_n K(r) Q_n |   \leq \sum_{l\geq 0} (2l+1) I_{2l+n+1}^2(r). \] 
Next we insert the estimate (\ref{Bessel-upperbd}) followed by (\ref{GammaBd}):
\begin{align*}
|\tr Q_n K(r)Q_n| &\leq \exp\big( \frac{r^2}{2(n+1)}\big) \sum_{l\geq 0} (2l+1) \Big(\frac{(r/2)^{2l+n+1}}{\Gamma(2l+n+2)}\Big)^2 \\
&\leq \exp\big( \frac{r^2}{2(n+1)}\big) \frac{(r/2)^{2n+2}}{\Gamma(n+2)^2} \sum_{l\geq 0}  (2l+1) \Big(\frac{r}{2(n+1)}\Big)^{4l} \\
&\leq 2\frac{\exp\big( \frac{r^2}{2(n+1)}\big) }{\big(1-(\frac{r}{2(n+1)})^4\big)^2} \frac{(r /2)^{2n+2}}{\Gamma(n+2)^2}
\end{align*}
which is the desired upper bound.
\end{proof}

By adapting the previous argument, we can also obtain the asymptotics of this trace.
\begin{lemma} \label{lem:trasymp}
For any $n\in\N$, it holds uniformly for all $ r\le \epsilon n $, as $\epsilon \to 0$, 
\[
\tr \big( K(r)Q_n \big) = I_{n+1}^2(r)  \big( (-1)^{n+1} + {O(\epsilon^2)} \big) 
\]
\end{lemma}

\begin{proof}
Using the estimate \eqref{Ideac} iteratively, for any $\ell \in \N$,
\[
 I_{\ell+1+n}^2(r) \le \bigg( \frac{r}{2n+3}  \bigg)^{2\ell}  I_{n+1}^2(r)
\]
so that according to \eqref{traceexact}, summing the geometric sum, we obtain  if $r \le \epsilon n$ and $0<\epsilon<1$,
\[
\tr Q_nK(r)Q_n  =   I_{n+1}^2(r) \big( (-1)^{n+1} + O(\epsilon^2) \big) .
\qedhere
\]
\end{proof}

We now show that higher-order traces are negligible.

\begin{lemma}\label{sumoftraces}
Let $\cst$ be as in Lemma~\ref{traceclassnorm1}. 
For any $n\geq 1$, if $\lvert \zeta \rvert \leq 2\cst(n+1)$, then
\[\sum_{j\geq 2} \frac{1}{j} \lvert \tr (Q_nK(|\zeta|)Q_n)^j \rvert \leq \frac{(1+\cst^2)(n+1)}{(n+1- \cst^2)(1-2\cst^2)} | \tr Q_nK(|\zeta|)Q_n |^2 . 
\]
\end{lemma}

\begin{proof}
First, recall the following properties of modified Bessel functions;
For any $r>0$, $\nu \mapsto I_\nu(r)$ is positive and non-decreasing on $\R_+$ and 
\[
\sum_{k\ge 0} I_{\nu+k}(r) \le  \frac{I_\nu(r)}{1-r^2/4\nu}
\]
These statements follow directly from the representation \eqref{BesselI}. 
Since  $K(r)= \mathbb{J}(r) \mathbb{L}   \mathbb{J}(r) \mathbb{L}$,
\[
 \tr (Q_nKQ_n)^j = \sum_{k_1,\dots,k_{2j}\geq 0}   (\mathbb{J}\mathbb{L}) _{k_1+n , k_2+1}  
   (\mathbb{J}\mathbb{L}) _{k_2+1 , k_3+n}     \cdots 
   (\mathbb{J}\mathbb{L}) _{k_{2j-1}+n, k_{2j}+1}  
   (\mathbb{J}\mathbb{L})_{k_{2j}+1,k_1+n}
\]
 By taking absolute value, using the definition of $\mathbb{L}$ and the Hankel matrix $\mathbb{J}$, cf.~\eqref{charfunction}, we can bound
\begin{align*}
\lvert \tr (Q_nK(r)Q_n)^j \rvert &\leq \sum_{k_1,\dots,k_{2j}\geq 0}  I_{k_1+k_2+n+1}(r)I_{k_2+k_3+n+1}(r)\dots I_{k_{2j-1}+k_{2j}+n+1}(r)I_{k_{2j}+k_1+n+1}(r) \\
&\leq \sum_{k_1,\dots,k_{2j}\geq 0}  I_{k_1+k_2+n+1}(r)I_{k_2+n+1}(r)\dots I_{k_{2j-1}+k_{2j}+n+1}(r)I_{k_{2j}+n+1}(r) \\
&\leq  \frac{1}{(1-r^2/4(n+1))^j} \sum_{k_2,k_4,\dots,k_{2j}\geq 0} I_{k_2+n+1}^2(r) \cdots I_{k_{2j}+n+1}^2(r) \\
& = \bigg( \frac{1}{1-r^2/4(n+1)} \sum_{k\ge 0} I_{k+ n+1}^2(r) \bigg)^j  \\
&\le \bigg( \frac{1+ r^2/4(n+1)^2}{1-r^2/4(n+1)} \sum_{k\ge 0} I_{2k+ n+1}^2(r) \bigg)^j
\end{align*}
where we used \eqref{Ideac} to obtain the last bound. Hence, using the lower-bound \eqref{trlb} to control the sum, we conclude that for any $j\in\N$ and $r< \sqrt{2} (n+1)$, 
\[
\lvert \tr (Q_nK(r)Q_n)^j \rvert
\le \bigg( \frac{(1+ r^2/4(n+1)^2)| \tr Q_nK_{\i g} Q_n | }{(1-r^2/4(n+1))(1-r^2/2(n+1)^2)}  \bigg)^j 
\]

Observe that $\sum_{j\ge 2} \frac{x^j}{j} \le x^2$ for $x\le 1/2$, so that using the estimate from Lemma \ref{traceclassnorm1}, it holds for $r <2\cst (n+1)$
\[
\frac{(1+ r^2/4(n+1)^2)| \tr Q_nK(r)Q_n | }{(1-r^2/4(n+1))(1-r^2/2(n+1))} \le \frac{(1+\cst^2) \Cst_0}{(n+1- \cst^2)(1-2\cst^2)} \le .16
\]
by a numerical evaluation. 
Hence,
\[ \begin{aligned}
\sum_{j\geq 2} \frac{1}{j} \lvert \tr (Q_nK_{\i g}Q_n)^j \rvert 
 & \le \bigg( \frac{(1+ r^2/4(n+1)^2)| \tr Q_nK_{\i g} Q_n | }{(1-r^2/4(n+1))(1-r^2/2(n+1)^2)}  \bigg)^2 \\
 &\le \frac{(1+\cst^2)(n+1)}{(n+1- \cst^2)(1-2\cst^2)} | \tr Q_nK_{\i g} Q_n |^2
\end{aligned}\]
as claimed.
\end{proof}

\begin{cor}\label{distfrom1}
Let $\cst$ be as in Lemma~\ref{traceclassnorm1}. 
For any $n\ge 1$, if $r< 2 \cst (n+1)$, then 
\[  \Big\lvert \frac{1- \det(1-K(r) Q_n) }{\tr (Q_nK(r)Q_n)} -1\Big\rvert  \le  \frac{1}{(n+1)(1-\cst^4)^2} .  \]
\end{cor}

\begin{proof}
By Plemelj's formula (Equation (5.12) in \cite{Simon}), provided that $\| Q_nK Q_n \|_{J_1} <1$, 
\[ \det(1-KQ_n) = \exp \tr \log (1-Q_nKQ_n) = \exp \Big( -\sum_{j\geq 1} \frac{\tr (Q_nKQ_n)^j}{j}\Big)  \]
We now use the bound valid for $z\in\C$,
\[
|e^{-z}-1+z | \le \frac{|z|^2 e^{|z|}}{2}
\]
and apply it to $z=\sum_{j\geq 1} \frac{\tr (Q_nKQ_n)^j}{j} $. 
By Lemma \ref{traceclassnorm1} and  \ref{sumoftraces}, we verify numerically that for any $n\in\N$ and  $\lvert \zeta \rvert < 2\cst(n+1)$, 
\[
\big| z- \tr (Q_nKQ_n) \big| \leq  \frac{3}{2} | \tr (Q_nKQ_n) |^2
\qquad \text{and}\qquad
| \tr (Q_nKQ_n) | \le \| Q_nKQ_n \|_{J_1} < .105
\]
so that 
\[
|z| \le 1.15 | \tr (Q_nKQ_n)| . 
\]
Hence, using the previous bound, we obtain
\[ \begin{aligned}
\bigg|\frac{e^{-z}-1+  \tr (Q_nKQ_n)}{\tr (Q_nKQ_n)} \bigg| & \le   \frac{3}{2} | \tr (Q_nKQ_n) | +  \frac{|z|^2 e^{|z|}}{2| \tr (Q_nKQ_n) |}  \\ 
&\le  3 | \tr (Q_nKQ_n) | . 
\end{aligned}\]

Consequently, by Lemma \ref{trace} we obtain if $\lvert \zeta \rvert < 2\cst(n+1)$, 
\[ \begin{aligned}
 \Big\lvert \frac{1- \det(1-K Q_n)}{\tr (Q_nKQ_n)} -1\Big\rvert 
 &\le 6\frac{\exp\big(\frac{\lvert\zeta\rvert^2}{2(n+1)}\big)}{\big(1-(\frac{\lvert\zeta\rvert}{2(n+1)})^4\big)^2}  \frac{(\lvert \zeta \rvert /2)^{2n+2}}{\Gamma(n+2)^2} \\
 &\le \frac{(e^{1+\cst^2}\cst)^{2n+2}}{(n+1)(1-\cst^4)^2} . 
\end{aligned}\]
By the definition of $\cst>0$, this proves the claim. 
\end{proof}

\subsection{Exact asymptotics: Proof of Theorem~\ref{prop:asymp}} \label{sec:asymp}
As a consequence of Lemma~\ref{lem:trasymp} and Corollary~\ref{distfrom1}, it holds uniformly in the regime $r\le \epsilon n$,
\begin{equation} \label{Fredholmasymp}
1- \det(1-K(r) Q_n) = \tr\big(Q_nK(r)Q_n\big)  \big(1+\underset{n\to\infty}{o(1)}\big)
=  I_{n+1}^2(r) \big( (-1)^{n+1} + {O(\epsilon^2)} \big)   
\end{equation}
where $\epsilon\le \cst$ is small. 
This implies that 
\[
\int_{|\zeta| \le \epsilon n} \big|1- \det(1-K(|\zeta|) Q_n) \big|^2 e^{- |\zeta|^2 /2} \d^2\zeta 
=  \big(1 + {O(\epsilon^2)} \big) \int_{|\zeta| \le \epsilon n}   I_{n+1}^4(|\zeta|) e^{- |\zeta|^2 /2} \d^2\zeta  .
\]

Moreover, using the estimates from Proposition~\ref{prop:TB1} and Proposition~\ref{lem:Had1}, we immediately verify that there exists $\delta>0$ depending only on $\epsilon$ so  that as $n\to\infty$, 
\[
\int_{|\zeta| \ge \epsilon n} \big|  F_n(\zeta) \big|^2 \d^2\zeta = O(e^{-\delta n^2}) . 
\]

By \eqref{charfunction}, this implies that for any small $\epsilon>0$, there exists $\delta>0$ so that
\begin{equation} \label{asymp1}
\big\| F_n - \widehat{\phi_\C}   \big\|_{L^2}^2  =  \big(1 + {O(\epsilon^2)} \big) \int_{|\zeta| \le \epsilon n}   I_{n+1}^4(|\zeta|) e^{- |\zeta|^2 /2} \d^2\zeta
+ O(e^{-\delta n^2}) . 
\end{equation}

Now, using the bound \eqref{Bessel-upperbd} and a change of variable, for $n\ge 2$,
\[ \begin{aligned}
\int_{|\zeta| \ge \epsilon n}   I_{n+1}^4(|\zeta|) e^{- |\zeta|^2 /2} \d^2\zeta
& \le  \int_{|\zeta| \ge \epsilon n}  \frac{|\zeta/2|^{4(n+1)}}{\Gamma(n+2)^4}  e^{- |\zeta|^2(1-\frac{2}{n+1})/2}  \d^2\zeta  \\
& = \frac{2^{-2n-1}\pi}{\Gamma(n+2)^4}   \bigg(1-\frac{2}{n+1}\bigg)^{-2(n+1)}  \int_{u \ge \Lambda_n}
u^{2(n+1)} e^{- u}\d u
\end{aligned}\]
where $\Lambda_n =(\frac12-\frac{1}{n+1}) (\epsilon n)^2$. By Markov's inequality, this implies that 
\[
\int_{|\zeta| \ge \epsilon n}   I_{n+1}^4(|\zeta|) e^{- |\zeta|^2 /2} \d^2\zeta
\le \frac{2^{-2n-1}\pi e^4}{\Gamma(n+2)^4 \Lambda_n} \int_0^\infty u^{2n+3} e^{- u}\d u =  \frac{2^{-2n-1}\pi e^4 \Gamma(2n+4)}{\Gamma(n+2)^4 \Lambda_n} . 
\]
Hence, using the inequalities \eqref{Gamma}, there exists a constant $C>0$ so that 
\[
\int_{|\zeta| \ge \epsilon n}   I_{n+1}^4(|\zeta|) e^{- |\zeta|^2 /2} \d^2\zeta
\le C \frac{(2n+3)/\Lambda_n}{\sqrt{n+1} \Gamma(n+2)^2} .
\]
In particular, we have as $n\to\infty$
\[
\int_{|\zeta| \ge \epsilon n}   I_{n+1}^4(|\zeta|) e^{- |\zeta|^2 /2} \d^2\zeta
= \frac{o_\epsilon(1)}{\sqrt{n+1} \Gamma(n+2)^2} . 
\]

If we combine this estimate with \eqref{asymp1}, we have shown that for any small $\epsilon>0$, 
\begin{equation} \label{asymp2}
\big\| F_n - \widehat{\phi_\C}   \big\|_{L^2}^2  =  \big(1 + {O(\epsilon^2)} \big) \int_{\C}   I_{n+1}^4(|\zeta|) e^{- |\zeta|^2 /2} \d^2\zeta
+  \frac{o_\epsilon(1)}{\sqrt{n+1} \Gamma(n+2)^2} .
\end{equation}

To complete the proof of Proposition~\ref{prop:asymp}, it remains to obtain the asymptotics of the integral on the right-hand side of \eqref{asymp2} which are provided by the next Proposition with $\nu = n+1$. 

\begin{prop}
\[
\lim_{\nu\to\infty} \pi\sqrt{\nu}\Gamma(\nu+1)^2 \int_0^\infty  I_{\nu}^4(r) e^{- r^2 /2} \d r^2  = 2 e^{4} \sqrt{\pi} .
\]
\end{prop}

\begin{proof}
By formula \eqref{BesselI} and a change of variable,
\[ \begin{aligned}
\int_0^\infty  I_{\nu}^4(r) e^{- r^2 /2} \d r^2  
& = 4 \nu \int_0^\infty  I_{\nu}^4(2\sqrt{\nu u}) e^{-2 \nu u} \d u  \\
& = \frac{4 \nu^{2\nu+1}}{\Gamma(\nu+1)^4} \int_0^\infty \bigg( \sum_{j =0}^{\infty}  \frac{u^{j} \Gamma(\nu+1) \nu^j}{j!\Gamma(j+\nu+1)} \bigg)^4 u^{2\nu} e^{-2 \nu u} \d u 
\end{aligned}\]

Observe that $\frac{\Gamma(\nu+1) \nu^j}{\Gamma(j+\nu+1)} \le 1$ for all $j\in\N_0$ so that 
\[ \begin{aligned}
\int_0^\infty  I_{\nu}^4(r) e^{- r^2 /2} \d r^2 
\le  \frac{4 \nu^{2\nu+1}}{\Gamma(\nu+1)^4} \int_0^\infty e^{4u-2\nu \phi(u)} \d u 
\end{aligned}\]
where $\phi(u) = u - \log(u) $. This function has a unique (non-degenerate) minimum at $u=1$ so that by Laplace method, as $\nu \to\infty$, 
\[
 \int_0^\infty e^{4u-2\nu \phi(u)} \d u 
\sim e^{4-2\nu} \sqrt{\pi/\nu} .
\]

It follows that 
\begin{equation} \label{asymp3}
\limsup_{\nu\to\infty} \frac{\sqrt{\nu}\Gamma(\nu+1)^4}{\nu^{2\nu+1} e^{4-2\nu}} 
\int_0^\infty  I_{\nu}^4(r) e^{- r^2 /2} \d r^2  \le 4\sqrt{\pi}
\end{equation}
On the other-hand, for any $\kappa \in \N$
\[
\int_0^\infty  I_{\nu}^4(r) e^{- r^2 /2} \d r^2  
\ge \frac{4 \nu^{2\nu+1}}{\Gamma(\nu+1)^4} \int_0^\infty \bigg( \sum_{j \le \kappa}  \frac{u^{j} \Gamma(\nu+1) \nu^j}{j!\Gamma(j+\nu+1)} \bigg)^4 u^{2\nu} e^{-2 \nu u} \d u 
\]
Observe that by monotonicity, for $j\le \kappa$,
\[
\frac{\Gamma(\nu+1) \nu^j}{\Gamma(j+\nu+1)} \ge \frac{\Gamma(\nu+1) \nu^\kappa}{\Gamma(\kappa+\nu+1)}
\]
so that
\[
\int_0^\infty  I_{\nu}^4(r) e^{- r^2 /2} \d r^2  
\ge \frac{4 \nu^{2\nu+1}}{\Gamma(\nu+1)^4}  \bigg( \frac{\Gamma(\nu+1) \nu^\kappa}{\Gamma(\kappa+\nu+1)} \bigg)^4
 \int_0^\infty \bigg( \sum_{j \le \kappa}  \frac{u^{j}}{j!} \bigg)^4 e^{-2\nu \phi(u)} \d u  
\] 
Hence for any fixed $\kappa\in\N$, 
\begin{equation} \label{asymp4}
\liminf_{\nu\to\infty} \frac{\sqrt{\nu}\Gamma(\nu+1)^4}{\nu^{2\nu+1} e^{4-2\nu}} 
\int_0^\infty  I_{\nu}^4(r) e^{- r^2 /2} \d r^2  \ge 4\sqrt{\pi} \bigg( \sum_{j \le \kappa}  \frac{1}{j!} \bigg)^4 e^{-4} .
\end{equation}

Combining \eqref{asymp3} and  \eqref{asymp4}, letting $\kappa\to\infty$,  
this proves that 
\[
\lim_{\nu\to\infty} \frac{\sqrt{\nu}\Gamma(\nu+1)^4}{\nu^{2\nu+1} e^{4-2\nu}} 
\int_0^\infty  I_{\nu}^4(r) e^{- r^2 /2} \d r^2  = 4\sqrt{\pi}
\]
To complete, it remains to use that as $\nu\to\infty$
\[
 \frac{\Gamma(\nu+1)^2}{\nu^{2\nu+1} e^{-2\nu}} 
 \sim 2\pi \qedhere
\]
\end{proof}

\subsection{Proof of the estimates~\eqref{ub1} and \eqref{ubinfty}} 
\label{firstregime}
It follows immediately from Corollary \ref{distfrom1} that for $n\ge 1$ and $\lvert \zeta \rvert < 2 \cst (n+1)$
\[
| \det[1 -K(|\zeta|)  Q_n]-1 | \leq 1.3  | \tr (Q_n K(|\zeta|)  Q_n) |  . 
\]
Then, using formula (\ref{charfunction}), we obtain
\[\begin{aligned}
\int_{|\zeta| < 2\cst n} \big| F_n(\zeta) - e^{-|\zeta|^2/4} \big|^2 \d^2\zeta 
&=
\int_{|\zeta| <2\cst n} e^{-|\zeta|^2/2} \big( 1- \det[\operatorname{I} -  K(|\zeta]) Q_n] \big)^2 \d^2\zeta   \\
&\le 1.75
\int_{r < 2\cst n}  |\tr (Q_nK(r)Q_n)|^2e^{-r^2/2}\d^2\zeta. \end{aligned} \]
Hence, Lemma \ref{trace} gives
\[
\int_{|\zeta| < 2\cst n} \big| F_n(\zeta) - e^{-|\zeta|^2/4}   \big|^2 \d^2\zeta  \le  \frac{7}{(1-\cst^4)^4} \frac{2^{-4n-4}}{\Gamma(n+2)^4} 
\int_{\C}  e^{-|\zeta|^2(\frac{1}{2}-\frac{1}{n+1})}  |\zeta|^{4n+4}  \d^2\zeta.
\]
A change of variables to polar coordinates gives for $n\ge 2$, 
\[\begin{aligned}
\int_{|\zeta| < 2\cst n} \big| F_n(\zeta) - e^{-|\zeta|^2/4}   \big|^2 \d^2\zeta  
&\le   \frac{7\pi}{(1-\cst^4)^4} \frac{2^{-4n-3}}{\Gamma(n+2)^4}
\int_{0}^{+\infty}  e^{-u(\frac{1}{2} -\frac{1}{n+1})}  u^{2n+2}  d u \\
& =    \frac{7\pi}{(1-\cst^4)^4} \frac{2^{-2n}}{\Gamma(n+2)^4} \big(1-\tfrac{2}{n+1}\big)^{-2n-3}  \Gamma(2n+3).
\end{aligned}\]
First we have using convexity,  for $n\ge 3$, 
\[
\big(1-\tfrac{2}{n+1}\big)^{-2n-3}  = e^{-(2n+3)\log(1-\frac2{n+1})} \le   e^{4+ \frac{7}{n}} .
\]
Second, using the duplication formula for the $\Gamma$ function and Gautschi's inequality (\href{https://dlmf.nist.gov/5.5}{(5.5.5)} and \href{https://dlmf.nist.gov/5.6}{(5.6.4)} in DLMF), we have 	for $n\in\N$
\begin{align} \notag
\Gamma(2n+3)   & = \frac{2^{2n+2}}{\sqrt{\pi}}  \Gamma(n+3/2)\Gamma(n+2) \\
\label{duplication}
&\le  \frac{2^{2n+2}}{\sqrt{\pi(n+1)}}  \Gamma(n+2)^2
\end{align}
Thus, we obtain for $n\ge 3$,
\[\begin{aligned}
\int_{|\zeta| < 2\cst n} \big| F_n(\zeta) - e^{-|\zeta|^2/4}   \big|^2 \d^2\zeta  
\le \frac{28 \sqrt\pi }{(1-\cst^4)^4}  \frac{ e^{4+ 7/n} }{\sqrt{n+1}}  \frac{1}{\Gamma(n+2)^2}  .
\end{aligned}\]
This completes the proof of the estimate (\ref{ub1}).

A similar argument shows that for $n\ge 2$,
\begin{equation*} \begin{aligned}
\int_{|\zeta| < 2\cst n} \big| F_n(\zeta) - e^{-|\zeta|^2/4}   \big| \d^2\zeta 
 & \leq 1.3 \int_{|\zeta| < 2\cst n} | \tr (Q_nK(|\zeta|)Q_n) | e^{-|\zeta|^2/4} \d^2\zeta \\
 &\le   \frac{1.3 \pi \cdot 2^3}{(1-\cst^4)^2} \big(1-\tfrac{2}{n+1}\big)^{-n-2}  \frac{1}{\Gamma(n+2)^2} \int_0^\infty e^{-u} u^{n+1} d u \\
& \le \frac{5.2\, e^{2+5/n}}{(1-\cst^4)^2} \frac{2\pi}{\Gamma(n+2)}   . 
\end{aligned}
\end{equation*}
This proves \eqref{ubinfty}.

\subsection{Proof of the lower-bounds in Theorem~\ref{thm:1}}
It  follows immediately from Corollary \ref{distfrom1} that for $n\ge 11$ and $r < 2 \cst (n+1)$
\[
| \det(1+K(r) Q_n)-1 | \geq .9 | \tr (Q_nK(r)Q_n) |  . 
\]
Then, using Equation (\ref{charfunction}), we obtain
\[\begin{aligned}
\| p_n - \phi_\C \|_{L^2}^2  &\ge \int_{|\zeta| < 2\cst (n+1)} \big| F_n(\zeta) - e^{-|\zeta|^2/4} \big|^2 \d^2\zeta 
=
\int_{|\zeta| <2\cst (n+1)} e^{-|\zeta|^2/2} \big( 1- \det[\operatorname{I} -  K(|\zeta|) Q_n] \big)^2 \d^2\zeta   \\
&\ge .81
\int_{r < 2\cst (n+1)}  |\tr (Q_n(r)Q_n)|^2e^{-r^2/2}\d^2\zeta. \end{aligned} \]
Using the lower-bound from Lemma \ref{trace}, this yields
\[
\| p_n - \phi_\C \|_{L^2}^2
  \ge   \frac{.81(1-2\cst^2)^2}{\Gamma(n+2)^4} 
\int_{|\zeta| < 2\cst (n+1)} e^{-|\zeta|^2/2}  (|\zeta|/2)^{4n+4}  \d^2\zeta.
\]
A change of variables $u = |\zeta|^2/4(n+1)$ gives for $n\ge 11$, 
\[
\| p_n - \phi_\C \|_{L^2}^2
  \ge   \frac{3.7\pi (1-2\cst^2)^2}{\Gamma(n+2)^4}  (n+1)^{2n+3}
\int_0^{\cst^2(n+1)} e^{- (2n+2)(u -  \log u)} \d u.
\]
The phase $u \in \R_+ \mapsto  u - \log u$ has a unique critical point at $u=1$. For $n\ge 13$, we verify that $\cst^2(n+1) > 1.5$ and $u - \log u -1 \le (u-1)^2$  for $u \in [.5,1.5]$ by convexity so that
\[\begin{aligned}
\int_0^{\cst^2(n+1)} e^{- (2n+2)(u -  \log u)} \d u 
 & \ge e^{- (2n+2)} \int_{.5}^{1.5} e^{-(2n+2)(u -1)^2} \d u \\
 & = \frac{\sqrt{2}e^{- (2n+2)}}{\sqrt{n+1}} \int_0^{\sqrt{(n+1)/2}} e^{-u^2} \d u \\
 & \ge \frac{1.25 e^{- (2n+2)}}{\sqrt{n+1}} . 
\end{aligned}\]
This numerical estimate shows that 
\[
\| p_n - \phi_\C \|_{L^2}^2 \ge 2.83 \pi  \frac{(n+1)^{2n+3} e^{- (2n+2)}}{\Gamma(n+2)^4\sqrt{n+1}} 
\]
Using the upper-bound \eqref{Gamma}, we conclude that for $n\ge 13$, 
\[
\| p_n - \phi_\C \|_{L^2}^2 \ge \frac{2.83}{2e^{\frac1{90}}}\frac{1}{\Gamma(n+2)^2}  \ge \frac{1}{\Gamma(n+2)^2\sqrt{n+1}} . 
\]
This completes the proof of the lower-bound \eqref{ub0}.

\subsection{Proof of Proposition \ref{prop:TB1}}\label{sec:tb}
The proof is based on the following bound which first appeared in \cite{97}. The version we use is a special case of Lemma 2.9 in \cite{JL}. 
Let us denote by $\{e^{\i \theta_j}\}_{j=1}^n$ the eigenvalues of the random matrix $\mathbf{U}\in \mathbb{U}(n)$. 

\begin{lemma}\label{changeofvar1}
Let $\nu >0$ and set  $g:\T\rightarrow \R$, $g(\theta) = \Re(\zeta e^{\i\theta})$. 
Let $h:\T\rightarrow \R$ be its Hilbert transform, i.e. $h(\theta) = - \Im (\zeta e^{\i \theta})$. For any $\zeta \in \C$,  
\[
| F_n(\zeta) | \le  \E_n\bigg[ \exp\bigg( \frac{\nu^2}{2n^2} \sum_{j,k = 1}^n H(\theta_j , \theta_k) - \sum_{j=1}^n \Im g\big(\theta_j +\i \frac{\nu}{n}h(\theta_j)\big)     \bigg)  \bigg],
\]
where 
\[
H(\theta,x) = \bigg( \frac{h(\theta) - h(x)}{2\sin(\frac{\theta-x}{2})} \bigg)^2. 
\]
\end{lemma}

The proof of Lemma~\ref{changeofvar1} relies on the fact that $F_{n}(\zeta) = \E_n\big[e^{\i \sum_{j=1}^n g(\theta_j) } \big]$,  the explicit formula for the joint density of the eigenvalues of~$\mathbf{U}$ (known as Weyl's integration formula  \cite{Meckes}) and a change of variables $\theta_j \leftarrow \theta_j +\i \frac{\nu}{n}h(\theta_j)$. 
In particular, the quadratic term $H$ comes from the Jacobian of this change of variables and it turns out that the \emph{optimal choice} for $h$ is the Hilbert transform of~$g$.

\begin{lemma}\label{monotonicity}
Suppose that $f \in L^{\infty}(\T)$ is real-valued with $\|f\|_{H^{1/2}}^2  =2\mathcal{A}(f) < \infty$ where $\mathcal{A}$ is as in (\ref{mathcalA}). Then, for
any $n\in\N$,
\[ \E_n[e^{\tr f(\mathbf{U})}] \leq \exp (n\hat{f}_0+ \mathcal{A}(f)). \]
\end{lemma}
\begin{proof}
This follows directly from the Borodin-Okounkov formula (Theorem \ref{borodin-okounkov}). 
Without loss of generality, we can assume that $\hat{f}_0 = 0$. Then,
if $f$ is real-valued, we have
$\hat{f}_{-k} = \overline{\hat{f}_{k}}$ 
and, according to \eqref{BOop}, $\omega_+ = \overline{\omega_-}$ so that 
$H_+(\omega_+) = H_-(\omega_-)^*$. 
Hence, the operator $K$ is positive semi-definite and the Fredholm  determinant $\det(\operatorname{I} - KQ_n)\leq 1$.
\end{proof}

\medskip
 
Before going into the proof of Proposition~\ref{prop:TB1}, we also need the following two lemmas.

\begin{lemma} \label{lem:Psi}
Let $\Psi(\theta)= \sin(\theta) \sinh( \delta \sin \theta)$ for $\theta\in\T$ and a fixed $\delta>0$. We have 
\begin{equation} \label{psi4}
\Psi(\theta) =   I_1(\delta) + 2 \sum_{k=1}^{+\infty}  I_{2k}'(\delta) \cos( 2k\theta)
\end{equation}
where $I_\nu$ are Bessel functions -- see formula \eqref{BesselI}.  
In particular, the Fourier series \eqref{psi4} converges uniformly.  
\end{lemma}

\begin{proof}
It follows from  DLMF formula \href{https://dlmf.nist.gov/10.32}{(10.32.3)} that for $k\in\Z$
\[
\int_\T \cosh(\delta \sin \theta) \cos(k\theta)  \frac{\d\theta}{2\pi}  = \frac{1+(-1)^k}{2}I_k(\delta)  .
\]
Since for any $\delta\in\R$, $\theta \in\T \mapsto   \cosh(\delta \sin \theta)$ is even, it follows that
\[ 
 \cosh(\delta \sin \theta) =  I_0(\delta)  + 2 \sum_{k\in \N} I_{2k}(\delta) \cos( 2k\theta) .
\]
Note that for $\delta>0$, the series converges uniformly since 
 $\theta \in\T \mapsto   \cosh(\delta \sin \theta)$ is smooth, and the same holds for all its derivatives with respect to $\delta$. Hence, if we differentiate term by term, we obtain 
 \[
 \Psi(\theta) = \partial_\delta \big( \cosh(\delta \sin \theta) \big)
 = I_0'(\delta)  + 2 \sum_{k\in \N} I_{2k}'(\delta) \cos( 2k\theta) .
 \]
 The claim now follows from that $I_0' = I_1$ by 
 DLMF formula \href{https://dlmf.nist.gov/10.29}{(10.29.3)}. 
\end{proof}


\begin{lemma} \label{lem:A}
Let $\Psi$ be as in Lemma~\ref{lem:Psi} for  a fixed $\delta>0$ and $\mathcal{A}$ as in Theorem \ref{borodin-okounkov}. We have  for any $r>0$ and $x\in\R$, 
\[
\A\big( x \cos( \cdot) -  r \Psi\big) \le  \frac{x^2 + 2r^2\delta^2 \varrho(\delta)}{4} ,
\]
where $\varrho(\delta) =  e^{\delta^2/6}\frac{I_0(\delta) + J_0(\delta)}{2}$. 
\end{lemma}

\begin{proof}
First of all, according to Lemma~\ref{lem:Psi}, since $\widehat{\Psi}_1 =0$,  we check  that 
\[
\A\big( x \cos( \cdot) -  r \Psi\big) =  \frac{x^2}{4} + r^2 \A(\Psi) . 
\]
So we only need to estimate the semi--norm 
\[
\A(\Psi) = \sum_{k=1}^{+\infty}  k |\widehat{\Psi}_k|^2 = 8  \sum_{k=1}^{+\infty}  k |I_{2k}'(\delta)|^2
\]
Using the formula \eqref{BesselI} followed by \eqref{GammaBd}, we see that for any integer $k\ge 1$ and $\delta>0$, 
\begin{equation*} 
I_{2k}'(\delta) =  \frac{2}{\delta} \sum_{j =0}^{+\infty} (k+j) \frac{(\delta/2)^{2(k+j)}}{j! (j+2k)!} \le \frac{2}{\delta (2k)!}  \sum_{j =0}^{+\infty} \frac{(\delta/2)^{2(k+j)}}{j!(2k+1)^j} \le \frac{2 (\delta/2)^{2k}}{\delta (2k)!}   e^{\delta^2/12} . 
\end{equation*}
Therefore, we obtain for any $\delta>0$, 
\[
\A(\Psi) \le  \frac{32}{\delta^2}  e^{\delta^2/6} \sum_{k=1}^{+\infty}  \frac{k (\delta/2)^{4k}}{(2k)!^2} 
\le    \frac{\delta^2}{2}  e^{\delta^2/6} \sum_{k=0}^{+\infty}  \frac{(\delta/2)^{4k}}{(2k)!^2} ,
\]
where we used that $\frac{k+1}{(2k+2)^2} \le \frac 14$ for any $k\ge 0$. 
Since $\frac{I_0(\delta) + J_0(\delta)}{2} =  \sum_{k=0}^{+\infty}   \frac{(\delta/2)^{4k}}{(2k)!^2}$, by \eqref{BesselI},  this completes the proof. 
\end{proof}

We are now ready to give the proof of Proposition~\ref{prop:TB1}.  

\begin{proof}
By Lemma \ref{changeofvar1}, we have for any  $\zeta\in\C$, 
\begin{equation} \label{ub6}
| F_n(\zeta) | \le  \E_n\bigg[ \exp\bigg( \frac{\nu^2}{2n^2} \sum_{i,j = 1}^n H(\theta_i , \theta_j) - \sum_{j=1}^n \psi(\theta_j)  \bigg)  \bigg] 
\end{equation}
where 
\begin{equation} \label{psih}
\psi(\theta) = \Im g\big(\theta +\i \frac{\nu}{n}h(\theta)\big)   
\qquad\text{and}\qquad
H(\theta,x) = \bigg( \frac{h(\theta) - h(x)}{2\sin(\frac{\theta-x}{2})} \bigg)^2. 
\end{equation}
Using that $\frac{e^{\i \theta} -e^{\i x}}{2} = \i e^{\i(\theta+x)/2} \sin(\frac{\theta-x}{2}) $, we obtain 
\begin{align*}
H(\theta,x)  &= \bigg( \frac{\Im\big( \zeta ( e^{\i \theta} -e^{\i x}) \big)}{2\sin(\frac{\theta-x}{2})} \bigg)^2\\
&= \big( \Re( \zeta e^{\i(\theta+x)/2})\big)^2  \\
&=  \frac 12 \Big( \Re \big(\zeta^2 e^{\i(\theta+x)}\big) + |\zeta|^2 \Big) . 
\end{align*}
Hence, this implies that 
\begin{equation} \label{H3}
\sum_{i,j =1}^n H(\theta_i , \theta_j) = \frac{n^2|\zeta|^2}{2}  +   \frac 12 \Re \big( \zeta \tr \u\big)^2
\end{equation}
If we write $\zeta = |\zeta|e^{\i\varphi}$, let us observe that for any $\delta\in\R$,
\[
\exp\left( \delta^2 \Re( e^{\i\varphi}\tr \u )^2\right) =
 | \exp\left( \delta^2 ( e^{\i\varphi}\tr \u )^2\right) |
\le
 \int_\R e^{- x^2+  2\delta x \Re( e^{\i\varphi}\tr \u ) }  \frac{dx }{\sqrt{\pi}} .
\]
This shows that with $\delta= \frac{\nu}{n}|\zeta|$, 
\[
\exp\bigg( \frac{\nu^2}{2n^2} \sum_{i,j=1}^n H(\theta_i , \theta_j) \bigg) \le    e^{\nu^2|\zeta|^2/4} \int_\R e^{- x^2+  \delta x\Re( e^{\i\varphi}\tr \u ) }  \frac{dx }{\sqrt{\pi}} . 
\]
By \eqref{ub6}, this implies that  
\begin{equation} \label{ub7}
| F_n(\zeta) |
\le  e^{\nu^2|\zeta|^2/4} \int_{\R} 
e^{- x^2} \E_n\bigg[ \exp\Big(  {\textstyle\sum_{j=1}^n} \delta x \cos(\theta_j +\varphi)  - \psi(\theta_j) \Big) \bigg]  \frac{\d x}{\sqrt{\pi}} . 
\end{equation}

\medskip

Now, observe that since $g(z) = \frac{1}{2}(\zeta e^{\i z} +  \overline{\zeta} e^{-\i z})$ for all $z\in\C$ and $h(\theta)  =-  |\zeta| \sin(\theta+\varphi)$ for all $\theta\in\T$, we deduce from \eqref{psih} that 
\begin{equation} \label{psi3}
\begin{aligned}
\psi(\theta) & =  \frac{1}{2}  \Im \big( \zeta e^{\i\theta + \delta  \sin(\theta+\varphi) }  +  \overline{\zeta} e^{-\i\theta - \delta \sin(\theta+\varphi) }   \big)\\
& =   |\zeta|   \Psi(\theta+\varphi)  . 
\end{aligned}
\end{equation}
where $\delta= \frac{\nu}{n}|\zeta|$ and $\Psi(\theta) = \frac 12 \Im \left( e^{\i\theta + \delta\sin \theta }  + e^{-\i\theta - \delta\sin \theta  }  \right)$ as in Lemma~\ref{lem:Psi}.
Then, we deduce from \eqref{ub7}, \eqref{psi3} and the invariance by rotation of the CUE law that for any  $\delta >0$, 
\begin{equation}\label{ub8} 
| F_n(\zeta) |
\le  e^{\nu^2|\zeta|^2/4} \int_{\R} 
e^{- x^2} \E_n\bigg[ \exp\Big(  {\textstyle\sum_{j=1}^n} \delta x \cos(\theta_j)  - |\zeta| \Psi(\theta_j) \Big) \bigg]  \frac{\d x}{\sqrt{\pi}} . 
\end{equation}
Moreover, by combining Lemma~\ref{monotonicity} and  Lemma~\ref{lem:A},  we have 
\[  \begin{aligned} 
\E_n\bigg[ \exp\Big(  {\textstyle\sum_{j=1}^n} 2\delta x \cos(\theta_j)  -  |\zeta| \Psi(\theta_j) \Big) \bigg] 
&\le  \exp\Big( -  n |\zeta| \widehat{\Psi}_0  +   \A\big(\delta x \cos( \cdot) -  |\zeta| \Psi\big)\Big) \\
&\le \exp\Big( - n |\zeta| I_1(\delta)  + \delta^2 x^2/4 +    |\zeta|^2 \delta^2  \varrho(\delta)/2    \Big) . 
\end{aligned} \]
where we used that $\widehat{\Psi}_0 =   I_1(\delta)$. 
Then, we deduce from \eqref{ub8} that if $\delta<2$
\begin{align} \notag
| F_n(\zeta) |
& \le  e^{-  n |\zeta| I_1(\delta)+\nu^2|\zeta|^2/4 + |\zeta|^2 \delta^2 \varrho(\delta)/2} \int_{\R} 
e^{- x^2(1- \delta^2/4) }  \frac{\d x}{\sqrt{\pi}} \\
&\label{ub9} 
= \frac{e^{-n |\zeta| I_1(\delta)+\delta^2(n^2 +  2|\zeta|^2  \varrho(\delta))/4}}{\sqrt{1-\delta^2/4}} ,
\end{align}
where we replace the parameter $\nu$ in the last step using the condition $\nu |\zeta| =n\delta$.
We would like to minimize the right-hand side of \eqref{ub9} over all  $\delta >0 $. 
Let us observe that  from  formula \eqref{BesselI}, the function $2 I_1(\delta)/\delta \ge 1$  for all $\delta>0$. Moreover, the function  $\delta \mapsto 2\varrho(\delta)$ is smooth, increasing, and bounded from above by $2.05$ for $0<\delta \leq \tfrac{1}{2\sqrt{2}}$ (by a numerical evaluation).
Therefore,  by \eqref{ub9}, we obtain for any  $0< \delta \le \tfrac{1}{2\sqrt{2}}$, 
\begin{equation} \label{ub10}
| F_n(\zeta) | \le \frac{e^{\delta^2(n^2 +  2.05 |\zeta|^2  )/4-   n |\zeta| \delta/2}}{\sqrt{1-\delta^2/4}}  . 
\end{equation}
So it suffices to minimize over all $\delta \in [0,\tfrac{1}{2\sqrt 2}]$, the polynomial 
\[
\delta^2(n^2 + 2.05 |\zeta|^2)/2-   n |\zeta| \delta . 
\]
The minimum is attained for $\delta_* = \frac{n|\zeta|}{n^2+ 2.05 |\zeta|^2}$ (observe that $\delta_* <  \tfrac{1}{2\sqrt{2}}$ for all $\zeta\in\C$ and  $n\in\N$) and it is given by $\frac{-n^2 |\zeta|^2}{2(n^2+ 2.05 |\zeta|^2)}$. Therefore, we conclude from \eqref{ub10} that for  any  $|\zeta|  >0 $ and $n\in\N$, 
\[
| F_n(\zeta) | \le \sqrt{2}  \exp\left(-  \frac{n^2}{4(2.05 + n^2/ |\zeta|^2)} \right) . \qedhere
\] 
\end{proof}

\subsection{Proof of Proposition~\ref{lem:Had1}}
\label{sec:Had}
The proof relies on the Toeplitz determinant representation  \eqref{heine-szego} of the characteristic function $F_n$ and Hadamard's inequality.
Recall that $g(\theta)= r\sin(\theta +\varphi)$ if we write $\zeta = re^{\i(\varphi-\pi/2)}$. 
Then, for any $k\in\Z$,
\[\begin{aligned}
(\widehat{e^{\i g}})_k  & = \int_\T e^{ \i r \sin(\theta+\varphi) -\i k \theta} \frac{\d\theta}{2\pi}\\
&= e^{\i k \varphi} \int_0^\pi \cos( r\sin(\theta) -k \theta) \frac{\d\theta}{\pi}.
\end{aligned}\]
According to DLMF formula \href{https://dlmf.nist.gov/10.9}{(10.9.2)}, this implies that $(\widehat{e^{\i g}})_k =  e^{\i k \varphi }J_k(r)$ for any $k\in\Z$.
It is well--known that for any fixed $k\in\Z$,  we have the asymptotic expansion as $r\to+\infty$, 
\[
|J_k(r)|^2 =  \frac{1 + (-1)^k \sin(2r) +o(1) }{\pi r} ,
\]
see e.g. DLMF formula \href{https://dlmf.nist.gov/10.7}{(10.7.8)}. In \cite{Krasikov}, Theorem 2, Krasikov obtained the following (sharp) bound. Let $\mu = (k+1/2)(k+3/2)$. If $r>\sqrt{\mu+(\mu/2)^{2/3}}$ and $k>-1/2$, then
\[
|J_k(r)|^2 \le \frac{4(r^2-(k+1/2)(k+5/2))}{\pi (2(r^2-\mu)^{3/2} -\mu)} .
\]
Note that the function $\mu \mapsto (2(r^2-\mu)^{3/2} -\mu)$ is decreasing for $r>\sqrt{\mu+(\mu/2)^{2/3}}$ and $\mu \ge (k+1)^2$.
It follows that for all $r \ge n +(n/2)^{2/3}$  and all integer $k \in (-n,n)$, 
\[
|J_k(r)|^2 \le \frac{4 (r^2-(k+1)^2)}{\pi (2 (r^2-(k+1)^2)^{3/2} -(k+1)^2)} .
\]
This shows that 
\[
|J_k(r)|^2 \le \frac{2}{\pi r} \left( \sqrt{1-\epsilon} - \frac{\epsilon}{2r(1-\epsilon)} \right)^{-1} ; \qquad \epsilon= (\tfrac{k+1}{r})^2 <1 .
\]
The function $\epsilon \mapsto \left( \sqrt{1-\epsilon} - \frac{\epsilon}{2r(1-\epsilon)} \right)^{-1} $ is increasing on $[0,1)$ before it explodes to $+\infty$. Hence, this bound is monotone in $\epsilon>0$  and $r>0$. This implies that for any $n\ge 2$, for all $r \ge 7 n /4$  and all integer $k \in (-n,n)$, 
\begin{equation} \label{BesselJ}
|J_k(r)|^2 \le \frac{2}{\pi r} \left( \sqrt{1-\epsilon} - \frac{\epsilon}{7(1-\epsilon)} \right)^{-1} ; \qquad \epsilon= (2/3)^2 .
\end{equation}

\eqref{BesselJ} yields the numerical estimate; $|J_k(r)|^2 \le 1/r $ for all $r \ge 7 n /4$  and all integer $k \in (-n,n)$. 
 By (\ref{heine-szego}) and Hadamard's inequality,  we have
\[
\big|  F_n(\zeta) \big|^2  =  \big| \det(T_n(e^{\i g})) \big|^2 \le  \prod_{j=1}^n \sum_{i=1}^n \big| (\widehat{e^{\i g}})_{j-i}\big|^2  
= \prod_{j=1}^n \sum_{i=1}^n \big|J_{j-i}(|\zeta|)\big|^2 . 
\]
Using the uniform bound  \eqref{BesselJ}, this implies that for any $n\ge 2$ and  for all $r \ge 7 n /4$, 
\[
\big|  F_n(\zeta) \big|^2   \le  \left|\frac{n}{\zeta}\right|^{n}. 
\]
In the previous argument, we assumed that $n\ge 2$. However, 
in case $n=1$, there is an explicit formula $\widehat{p}_1(\zeta) =  J_0(|\zeta|)$ so that for all $\zeta\in\C$,
\[
|\widehat{p}_1(\zeta)|^2 \le \frac{2}{\pi |\zeta|} .
\]
This completes the proof.

\section{Orthogonal and symplectic groups: proof of Theorems \ref{thm:TVO} and \ref{prop:asympO}}
\label{sec:O}

In this section we consider Haar distributed orthogonal and symplectic matrices which we denote by $\mathbf{O}$. Observe that unlike unitary matrices, these are different from the circular ensembles. Another important difference is that in some cases there are deterministic eigenvalues at $\pm 1$ and that all random eigenvalues occur in conjugate pairs. Therefore we can make the change of variable $x_j=\cos \theta_j$ in the joint eigenvalue density, given by Weyl's integration formula (see \cite{Meckes}), to obtain 
\begin{align}\label{density}
\rho_n^{ab}(x) = \frac{1}{Z_n^{ab}}\prod_{1\leq j\leq n}(1-x_j)^a(1+x_j)^b\prod_{1\leq j<k\leq n}(x_j-x_k)^2, \quad x\in [-1,1]^n    
\end{align}
where $(a,b)=(1/2,1/2)$ for $Sp(2n)$ and $O(2n+2)^-$, $(a,b)=(-1/2,-1/2)$ for $O(2n)^+$, $(a,b)=(-1/2,1/2)$ for $O(2n+1)^-$, and $(a,b)=(1/2,-1/2)$ for $O(2n+1)^+$ (for sake of brevity we will often replace the indices $a$ and $b$ with their respective sign). Note that $n$ gives the number of non-trivial eigenvalues and $2n$ the number of random eigenvalues. We will denote by $d$ the total number of eigenvalues, i.e. the number of rows of the matrix, which is either $2n$, $2n+1$ or $2n+2$. We also write $F_n^{ab}(\xi)= \E_n^{ab}[e^{\i \xi\tr \mathbf{O}}]$ for the characteristic function of $\tr \mathbf{O}$, the sum of all random eigenvalues. 


The cosine of the random eigenangles form a determinantal point process, just as the eigenangles themselves. Therefore, by using Andr\'{e}ief's identity, it is possible to write $F_n^{ab}$ as the determinant of a matrix, more precisely a Toeplitz + Hankel matrix. We can then use an analogue of the Borodin-Okounkov formula for these types of determinants, namely the results of \cite{BE}, 
to express $F_n^{ab}$ as the characteristic function of a normal random variable multiplied by a certain Fredholm determinant, similar to (\ref{charfunction}). The details of these computations can be found in \cite{CJ}, here we merely state the results.

\begin{lemma}
\label{Kurt}
For any complex function $\psi$ on $[-1,1]$,
\begin{align*}
& \E_n^{-+}[\prod_{j=1}^n \psi(x_j)] = \det (\hat{\phi}_{j-k}+\hat{\phi}_{j+k+1})_{0\leq i,j\leq n-1} \\
& \E_n^{+-}[\prod_{j=1}^n \psi(x_j)] = \det (\hat{\phi}_{j-k}-\hat{\phi}_{j+k+1})_{0\leq i,j\leq n-1} \\
& \E_n^{++}[\prod_{j=1}^n \psi(x_j)] = \det (\hat{\phi}_{j-k}-\hat{\phi}_{j+k+2})_{0\leq i,j\leq n-1} \\
& \E_n^{--}[\prod_{j=1}^n \psi(x_j)] = \det (\hat{\phi}_{j-k}+\hat{\phi}_{j+k})_{0\leq i,j\leq n-1}
\end{align*}
where $\hat{\phi}_n$ is the $n$:th Fourier coefficient of $\psi\circ \cos$.
\end{lemma}

The next proposition, obtained in \cite{BE}, is from now on going to be referred to as the Basor-Ehrhardt formula. We consider functions in the Besov class $B_1^1$, i.e. functions $\omega$ on the unit circle which satisfy
\begin{equation}
\| \omega \|_{B_1^1} := \int_{-\pi}^\pi \frac{1}{y^2} \int_ {-\pi}^\pi \lvert \omega(e^{\i x+\i y})+\omega(e^{\i x-\i y})-2\omega(e^{\i x}) \rvert dxdy < \infty.
\end{equation}
If $\omega \in B_1^1$ we let $\omega_+$ denote its projection on $B_{1+}^1$, the subspace of $B_1^1$ for which $\omega_k=0$ for $k<0$, and we write $\tilde{\omega}(e^{\i \theta})\coloneqq \omega(e^{-\i \theta})$.

\begin{prop} \cite{BE}
\label{BE}
Denote by $Q_n$ the projection operator acting on $l_2(\mathbb{N})$ that sets the first $n$ coefficients to zero, and let $H(c)$ be the Hankel operator with symbol $c\in L^\infty(\T)$, i.e. the bounded linear operator on $l_2(\mathbb{N})$ with matrix representation $H(c) = (c_{j+k+1})_{j,k=0}^\infty$, where $c_k$ is the $k$th Fourier coefficient of $c$. Assume that $b_+\in B_{1+}^1$ and set $a_+=\exp(b_+)$, $a=a_+\tilde{a_+} = \exp (b_++\tilde{b_+})$. Then,
\begin{align*}
&\det (\hat{a}_{j-k}+\hat{a}_{j+k+1})_{0\leq i,j\leq n-1}= \\
&\qquad \exp\Big(n[\log a]_0 +{\sum_{n=0}^\infty[\log a]_{2n+1}+\frac{1}{2}\sum_{n=1}^\infty n[\log a]_{n}^2}\Big) \det(1+Q_nH(a_ +^{-1}(e^{\i \theta})\tilde{a_+}(e^{\i \theta}))Q_n) \\
&\det (\hat{a}_{j-k}-\hat{a}_{j+k+1})_{0\leq i,j\leq n-1}= \\
&\qquad \exp\Big({n[\log a]_0 -\sum_{n=0}^\infty[\log a]_{2n+1}+\frac{1}{2}\sum_{n=1}^\infty n[\log a]_{n}^2}\Big) \det(1-Q_nH(a_ +^{-1}(e^{\i \theta})\tilde{a_+}(e^{\i \theta}))Q_n)\\
&\det (\hat{a}_{j-k}-\hat{a}_{j+k+2})_{0\leq i,j\leq n-1}= \\
&\qquad \exp\Big({n[\log a]_0 -\sum_{n=1}^\infty[\log a]_{2n}+\frac{1}{2}\sum_{n=1}^\infty n[\log a]_{n}^2}\Big) \det(1-Q_nH(e^{-\i \theta}a_ +^{-1}(e^{\i \theta})\tilde{a_+}(e^{\i \theta}))Q_n)\\
&\det (\hat{a}_{j-k}+\hat{a}_{j+k})_{0\leq i,j\leq n-1}= \\
&\qquad \exp\Big({n[\log a]_0 +\sum_{n=1}^\infty[\log a]_{2n}+\frac{1}{2}\sum_{n=1}^\infty n[\log a]_{n}^2}\Big) \det(1+Q_nH(e^{\i \theta}a_ +^{-1}(e^{\i \theta})\tilde{a_+}(e^{\i \theta}))Q_n)
\end{align*}
Here $[\log a]_k$ stands for the $k$th Fourier coefficient of $\log a$. The Fredholm determinants are well-defined because each Hankel operator is trace-class.
\end{prop}

Combining the above two results gives a new expression for certain averages over the orthogonal and symplectic group, more amenable for asymptotic analysis thanks to the Fredholm determinant. In particular, we obtain the following expression for the characteristic function $F_n^{ab}$ (compare with \ref{charfunction}):

\begin{prop}\label{Fred}
Let $F_n^{ab}(\xi)= \E_n^{ab}[e^{\i \xi\tr \mathbf{O}}]$, where the expectation is with respect to (\ref{density}). Then,
\[ F_n^{ab}(\xi) = e^{\i (b-a)\xi-\xi^2/2} \det (I+Q_nK_{i\xi}^{ab}Q_n), \quad \xi\in\C\]
where, for any $z\in\C$,
\begin{align*}
&K_{z}^{-+} = (J_{j+k+1}(-2z))_{j,k\geq0}, & & K_{z}^{+-} = (-J_{j+k+1}(-2z))_{j,k\geq 0}, \\
&K_{z}^{++} = (-J_{j+k+2}(-2z))_{j,k\geq0}, & & K_{z}^{--} = (J_{j+k}(-2z))_{j,k\geq0}
\end{align*}
and $(J_k)_{k\in\N}$ are Bessel functions of the first kind.
\end{prop}

The Fredholm determinant in the above theorem converges to one as $n$ tends to infinity (see the discussion following Corollary 1.6. in \cite{CJ}), so by the continuity theorem, $\tr \mathbf{O}$ (disregarding deterministic eigenvalues) converges to a real normal random variable with mean $-1$, $1$ or $0$ and variance $1$ (adding the deterministic eigenvalues gives a mean equal to zero in all cases).

Just as for the unitary group, we obtain our bound on the $L^2$-norm from Parseval's identity and careful estimates of the characteristic function $F_n^{ab}$ of $\tr \mathbf{O}$. The bound on the total variation then follows from the Cauchy-Schwarz inequality combined with a concentration inequality for $\tr \mathbf{O}$. We still use three different techniques to study $F_n^{ab}(\xi)$, one for each regime of $\xi$. In the small regime, we use Proposition \ref{Fred} and estimate how far the Fredholm determinants are from one (the lower bound will give the lower bound on the total variation). In the intermediate regime, we make a change of variable similar to the one for the unitary case, but then the Borodin-Okounkov formula is replaced with the Basor-Ehrhardt formula. In particular, this implies that we do not have any lemma equivalent to Lemma \ref{monotonicity}, i.e. $F_n^{ab}$ is not bounded by its limit, which gives a larger bound compared to the unitary case. Finally, Hadamard's inequality gives the bound on $F_n^{ab}$ in the large regime. 

\subsection{Proof of the upper bound in Theorem \ref{thm:TVO}}

The upper bound on the total variation follows from that on the $L^2$ norm which is given in the next theorem.

\begin{thm}\label{thm:1O}
Let $p_n^{ab}$ be the probability density of $\tr \mathbf{O}- \E[\tr \mathbf{O}]$, with $a$, $b$ specifying the sign of the determinant and the parity of the size of the matrix (see \eqref{density}). Then,
\[
\|p_n^{ab}-\widehat{\phi_\R} \|_{L^2} \leq \frac{\beta_2(n)}{(2n)^{1/4}\sqrt{\Gamma(2n+1)}} + 2\frac{n}{\sqrt{n-1}}e^{-\beta_3^{ab}(n)(n^2-n)} + \frac{e^{5/8}}{2\sqrt{n}}e^{-2n^2/e^{5/2}}
\]
for all $n\geq 2$, where $\beta_2(n)$ and $\beta_3^{ab}(n)$ are bounded, and given explicitly in \eqref{beta2} and \eqref{beta3}.
Consequently, for $n\geq 122$,
\[
\|p_n^{ab}-\widehat{\phi_\R} \|_{L^2} \leq  \frac{\beta_2(n)+\frac{2}{n}}{(2n)^{1/4}\sqrt{\Gamma(2n+1)}}.
\]
\end{thm}
The next four subsections are dedicated to the proof. We also need a concentration inequality for the trace.

\begin{lemma}\label{concentration}
If $0\leq L < 4e^{-3/2}(1+n^{-2.5}e^{-n/2})n$, then
\[
\mathbb{P}\big[\lvert \tr \mathbf{O} - \E[\tr \mathbf{O}] \rvert \geq \tfrac{L}{2}\big] \leq 2\exp \Big(-\frac{L^2}{8(1+n^{-2.5}e^{-n/2})}\Big). \]
\end{lemma}

The proof is given in the next subsection. Now, let $\boldsymbol{\gamma}_\R$ denote a real standard Gaussian random variable with density given by $x \mapsto \frac{1}{\sqrt{2\pi}}e^{-x^2/2}$, and let $p_n^{ab}$ be the density of $\tr \mathbf{O}-\E[\tr \mathbf{O}]$, where $a,b$ specify the parity of the size of the matrix and the sign of the determinant. By the Cauchy-Schwarz inequality, for any $L>0$,
\[ \| p_n^{ab}- \widehat{\phi_\R} \|_{L^1(\R)} \leq \sqrt{L}\| p_n^{ab}- \widehat{\phi_\R} \|_{L^2} + \mathbb{P}\big[\lvert \tr \mathbf{O} - \E[\tr \mathbf{O}] \rvert \geq \tfrac{L}{2}\big] + \mathbb{P}[\lvert \boldsymbol{\gamma}_\R \rvert > \tfrac{L}{2}]
\]
where
\begin{equation}\label{concentrationrealgaussian}
\mathbb{P}[\lvert \boldsymbol{\gamma}_\R \rvert > \tfrac{L}{2}] = 2 \int_{L/2}^\infty e^{-\frac{x^2}{2}} \frac{\d x}{\sqrt{2\pi}}= 2 \int_{0}^\infty e^{-\frac{(x+L/2)^2}{2}} \frac{\d x}{\sqrt{2\pi}} < 2e^{-\frac{L^2}{8}} \int_0^\infty e^{-\frac{xL}{2}}\frac{\d x}{\sqrt{2\pi}} = \frac{4e^{-\frac{L^2}{8}}}{\sqrt{2\pi}L}
\end{equation}
and, if $L< 4e^{-3/2}(1+n^{-2.5}e^{-n/2})n$, 
\[ \mathbb{P}\big[\lvert \tr \mathbf{O} - \E[\tr \mathbf{O}] \rvert \geq \tfrac{L}{2}\big] \leq  2\exp \Big(-\frac{L^2}{8(1+n^{-2.5}e^{-n/2})}\Big) \]
by Lemma \ref{concentration}. Thus,
\[
\| p_n^{ab}- \widehat{\phi_\R} \|_{L^1(\R)} \leq \sqrt{L}\| p_n^{ab}- \widehat{\phi_\R} \|_{L^2} + 2\Big(1+\frac{\sqrt{2}}{\sqrt{\pi}L}\Big) \exp \Big(-\frac{L^2}{8(1+n^{-2.5}e^{-n/2})}\Big).
\]
We now insert our previous bound on the $L^2$-norm: if $n\geq 122$,
\[\| p_n^{ab}- \widehat{\phi_\R} \|_{L^1(\R)} \leq  \sqrt{L}\frac{\beta_2(n)+\frac{2}{n}}{(2n)^{1/4}\sqrt{\Gamma(2n+1)}} + 2\Big(1+\frac{\sqrt{2}}{\sqrt{\pi}L}\Big) \exp \Big(-\frac{L^2}{8(1+n^{-2.5}e^{-n/2})}\Big).\]
It remains to optimize over $L$; we choose
\[ L = 2\sqrt{2}\sqrt{(1+n^{-2.5}e^{-n/2})(\log(\Gamma(2n+1))/2+(\log n)/4))} \] 
which satisfies the assumption for large enough $n\geq 37$. Use that $\log \Gamma(2n+1) \leq 2n\log n$ for $n\geq 2$, this yields
\begin{multline*}
\| p_n^{ab}- \widehat{\phi_\R} \|_{L^1(\R)} \leq 
\Big[\sqrt{2}\Big(1+\frac{1}{n^{2.5}e^{n/2}}\Big)^{1/4}\Big(1+\frac{1}{4n}\Big)^{1/4}\Big(\beta_2(n)+\frac{2}{n}\Big) + 
\frac{3}{(n\log n)^{1/4}}\Big] \frac{(\log n)^{1/4}}{\sqrt{\Gamma(2n+1)}}.
\end{multline*} 
A numerical evaluation of the expression in brackets finishes the proof of Theorem \ref{thm:TVO}.
\subsection{Trace-norm and Fredholm determinant estimates}

Recall that the trace-class norm of a compact operator $K$ is given by
$
\| K\|_{J_1} \coloneqq \sum_{j\geq 1} s_j
$
where $\{s_j\}_{j\geq1}$ are its singular values. If $K$ is given by the infinite matrix with entries $K_{jk}$, then
\begin{equation}\label{traceclass}
\|K\|_{J_1} \leq \sum_{j\geq 1} \Big( \sum_{k\geq 1} \lvert K_{jk}\rvert ^2 \Big)^{1/2}.
\end{equation}
This can be seen from the polar decomposition of $K$: $K=U|K|$, where $U$ is a partial isometry and $|K|= (K^*K)^{1/2}$, so if $\{e_j\}_{j\geq1}$ is the standard orthonormal basis of $l_2(\N)$, then
$$ \|K\|_{J_1} = \sum_{j\geq 1} \langle \lvert K\rvert e_j,e_j \rangle = \sum_{j\geq 1} \langle Ke_j,Ue_j\rangle.  $$
The Cauchy-Schwarz inequality and the fact that $\|U\|_{op}=1$ give 
\begin{align*}
\|K\|_{J_1} \leq \sum_{j\geq 1} \|Ke_j\| \|Ue_j\| \leq \sum_{j\geq 1} \|Ke_j\| = \sum_{j\geq 1} \langle Ke_j,Ke_j \rangle^{1/2}=  \sum_{j\geq 1} \Big( \sum_{k\geq 1} \lvert\langle Ke_j,e_k\rangle\rvert ^2 \Big)^{1/2}.
\end{align*}

\begin{lemma}\label{traceclassnorm}
For any pair $(a,b)=(\pm 1/2, \pm 1/2)$, if $\lvert z \rvert < 2e^{-5/4}n$, then 
\[ \| Q_nK_{z}^{ab}Q_n \|_{J_1} < 1. \]
\end{lemma}

\begin{proof}
Consider the case $(a,b)=(-1/2,-1/2)$, i.e. $Q_nK_{z}^{--}Q_n= (J_{j+k}(-2z))_{j,k\geq n}$ according to Proposition \ref{Fred}. By (\ref{traceclass}) and (\ref{Bessel-upperbd}),
\begin{equation*}
\|Q_nK_{z}^{--}Q_n\|_{J_1} \leq \sum_{j\geq 0} \Big( \sum_{k\geq 0} |J_{j+k+2n}(-2z)| ^2 \Big)^{1/2}  \leq   \sum_{j\geq 0} \Big( \sum_{k\geq 0} e^{2|z|^2/(j+k+2n+1)} \frac{|z|^{2k+2j+4n}}{\Gamma(j+k+2n+1)^2}  \Big)^{1/2}.
\end{equation*}
Inserting the inequality $\Gamma(j+k+2n+1) \geq \Gamma(2n+1)(2n+1)^{j+k}$, $j,k\geq 0$, gives
\begin{align*}\label{J1}
\|Q_nK_{z}^{--}Q_n\|_{J_1} &\leq e^{|z|^2/(2n+1)} \frac{|z|^{2n}}{\Gamma(2n+1)} \sum_{j\geq 0} \frac{|z|^j}{(2n+1)^j} \Big( \sum_{k\geq 0} \frac{|z|^{2k}}{(2n+1)^{2k}}  \Big)^{1/2} \\ \nonumber
&< e^{|z|^2/(2n+1)} \frac{|z|^{2n}}{\Gamma(2n+1)} \Big(1-\frac{|z|}{2n+1}\Big)^{-1}\Big(1-\frac{|z|^2}{(2n+1)^2}\Big)^{-1/2}.
\end{align*}
So if $|z|<2e^{-5/4}n < 2e^{-1}n$,
\begin{equation*}
\|Q_nK_{z}^{--}Q_n\|_{J_1} \leq \frac{e^{n/2}}{\sqrt{1-e^{-2}}(1-e^{-1})} \frac{|z|^{2n}}{\Gamma(2n+1)}. 
\end{equation*}
Thus, by Stirling's approximation (\ref{Gamma}),
\begin{equation}\label{boundJ1}
\|Q_nK_{z}^{--}Q_n\|_{J_1} \leq \frac{1}{2\sqrt{1-e^{-2}}(1-e^{-1})\sqrt{\pi n}} \Big(\frac{e^{5/4}\lvert z\rvert}{2n}\Big)^{2n} <1
\end{equation}
for all $n\geq 1$ and $\lvert z\rvert<2e^{-5/4}n$. The other cases are similar.

\end{proof}

\begin{lemma}\label{sumoftraces1}
Define 
\[ \beta_1(x) = \frac{e^{\xi^2/(x+1)}}{1-\xi^2/(x+1)^2}\frac{\xi^{x}}{\Gamma(x+1)}, \quad x\in\N. \]
For any $\lvert z \rvert < 2e^{-5/4}n$ and $n\in\N$,
\[\sum_{j\geq 1} \frac{1}{j} \lvert \tr (Q_nK_{z}^{ab}Q_n)^j \rvert \leq -\log(1-\beta_1(d)) \]
where $d$ is the total number of eigenvalues.
\end{lemma}

\begin{proof}
The proof is essentially the same for all four cases so we consider only $a=b=-1/2$. We have that 
\begin{align}\label{boundontraceofpower}
\lvert \tr (Q_nK_{z}^{--}Q_n)^j \rvert &= \lvert \sum_{k_1,\dots,k_{j}\geq n} J_{k_1+k_2}(-2z)\dots J_{k_{j-1}+k_{j}}(-2z) \rvert \nonumber\\
&\leq \sum_{k_1,\dots,k_j\geq n} \lvert J_{k_1+k_2}(-2z)\dots J_{k_{j-1}+k_j}(2z) \rvert \nonumber\\
&\leq \Big(\frac{e^{|z|^2/(2n+1)}|z|^{2n}}{\Gamma(2n+1)}\Big)^j \sum_{k_1,\cdots,k_j\geq 0} \prod_{l=1}^j \big(\frac{|z|}{2n+1}\big)^{2k_l} \nonumber \\
&= \Big(\frac{e^{|z|^2/(2n+1)}}{1-|z|^2/(2n+1)^2}\frac{|z|^{2n}}{\Gamma(2n+1)}\Big)^j 
\end{align}
where we used (\ref{Bessel-upperbd}). Thus, 
\[
\sum_{j\geq 1} \frac{1}{j} \lvert \tr (Q_nK_{z}^{--}Q_n)^j \rvert \leq \sum_{j\geq 1} \frac{1}{j} \beta_1^j(2n) \\
= - \log (1-\beta_1(2n))
\]
By \eqref{density}, $a=b=-1/2$ corresponds to the case $O(2n)^+$ which has $2n$ eigenvalues in total.
\end{proof}

We obtain the following bound on the truncated $L^2$ norm of the difference between the characteristic function of $\tr \mathbf{O} - \E[\tr \mathbf{O}]$ and that of a real Gaussian.
\begin{prop}\label{regime1}
Define 
\begin{equation}\label{beta2}
\beta_2(x) = \frac{\sqrt{1-e^{-5/2}}e^{1+\frac{1}{x}}}{(1-e^{-5/2}-1/(2\sqrt{\pi x}))^{3/2}} , \quad x\in\R.
\end{equation}
Then, for any $n\in\N$,
\[
\Big( \int_{|\xi|<\frac{2n}{e^{5/4}}} |F_n^{ab}(\xi)-e^{-\xi^2/2}|^2 d\xi\Big)^{1/2} \leq \frac{\beta_2(d/2)}{d^{1/4}\sqrt{\Gamma(d+1)}}\]
with $d$ the total number of eigenvalues.
\end{prop}

\begin{proof}
Plemelj's formula, 
\[ \det (1+Q_nK_{\i \xi}^{ab}Q_n) = \exp \Big( \sum_{j\geq 1} \frac{(-1)^{j+1}}{j} \tr (Q_nK_{\i \xi}^{ab}Q_n)^j \Big) \]
holds for $\lvert \xi \rvert < 2e^{-5/4}n$ by Lemma \ref{traceclassnorm}. Set $\delta = \sum_{j\geq 1} \frac{(-1)^{j+1}}{j} \tr (Q_nK_{\i \xi}^{ab}Q_n)^j$. Then
\[ \lvert\det (1+Q_nK_{\i \xi}^{ab}Q_n)-1\rvert = \lvert e^\delta-1 \rvert < \lvert\delta\rvert e^{\lvert\delta\rvert}, \]
where
\[ \lvert \delta \rvert < \sum_{j\geq 1} \frac{1}{j} \lvert \tr (Q_nK_{\i \xi}^{ab}Q_n)^j \rvert < -\log ( 1- \beta_1(d)) \]
by Lemma \ref{sumoftraces1}.
Moreover, for $0<x<1$, $-\log(1-x)<x/\sqrt{1-x}$. Thus,
 \[ \lvert\det (1+Q_nK_{\i \xi}^{ab}Q_n)-1\rvert < ( 1- \beta_1(d))^{-3/2} \beta_1(d). \] 
Stirling's inequality (\ref{Gamma}) and the assumption that $\lvert \xi \rvert < 2e^{-5/4}n \leq e^{-5/4}d$ give
 \[ \beta_1(d) =  \frac{e^{\xi^2/(d+1)}}{1-\xi^2/(d+1)^2}\frac{\xi^{d}}{\Gamma(d+1)} < \frac{e^{d/4}}{1-e^{-5/2}} \frac{1}{\sqrt{2\pi d}}\Big(\frac{e\xi}{d}\Big)^{d} < \frac{1}{(1-e^{-5/2})\sqrt{2\pi d}}.  \]
We obtain 
\begin{align}\label{freddistance}
\lvert\det (1+Q_nK_{\i \xi}^{ab}Q_n)-1\rvert &< \frac{1}{(1-1/((1-e^{-5/2})\sqrt{2\pi d}))^{3/2}} \frac{e^{\xi^2/(d+1)}}{1-e^{-5/2}}\frac{\xi^{d}}{\Gamma(d+1)} \nonumber \\
&= \frac{\sqrt{1-e^{-5/2}}}{(1-e^{-5/2}-1/(\sqrt{2 \pi d}))^{3/2}} e^{\xi^2/(d+1)}\frac{\xi^{d}}{\Gamma(d+1)}
\end{align}
Combined with Proposition (\ref{Fred}) this gives
\begin{multline*}
 \Big( \int_{|\xi|<2n/e^{5/4}} |F_n^{ab}(\xi)-e^{-\xi^2/2}|^2 \d\xi  \Big)^{1/2} \\
\leq  \frac{\sqrt{1-e^{-5/2}}}{(1-e^{-5/2}-1/(\sqrt{2 \pi d}))^{3/2}} \frac{1}{\Gamma(d+1)} \Big(\int_{\R} e^{-\xi^2(1-\frac{2}{d+1})}|\xi|^{2d} \d\xi\Big)^{1/2}.
\end{multline*}
We make the change of variable $r=\xi^2(1-\frac{2}{d+1})$, 
\begin{align}\label{bd1}
\int_{\R} e^{-\xi^2(1-\frac{2}{d+1})}|\xi|^{2d} \d\xi = \Big(1-\frac{2}{d+1}\Big)^{-d-\frac{1}{2}}\Gamma(d+1/2)   
\end{align}
where 
$$ \Gamma(d+1/2) < \frac{\Gamma(d+1)}{\sqrt{d}} $$
by Gautschi's inequality (\href{https://dlmf.nist.gov/5.6}{(5.6.4)} in DLMF). Finally, observe that
$$ \Big(1-\frac{2}{d+1}\Big)^{-d-\frac{1}{2}} = e^{-(d+1/2)\log(1-2/(d+1))} \leq e^{(d+1/2)(2/(d+1)+4/(d+1)^2)} \leq  e^{2+4/d} $$
if $d\geq 2$.

\end{proof}

From the previous lemmas, we can also deduce the concentration result of $\tr \mathbf{O}- \E[\tr \mathbf{O}]$ that we gave in the previous subsection.
\begin{proof}[Proof of Lemma \ref{concentration}]
Consider the case $\mathbf{O} \in O(2n)^+$, i.e. $a=b=-1/2$. Then $\E[\tr \mathbf{O}] = 0$ and for any $\lambda>0$,
\[ 
\mathbb{P}\big[\lvert \tr \mathbf{O} \rvert \geq \tfrac{L}{2}\big] \leq e^{-\frac{\lambda L}{2}} (\mathbb{E}_n^{--}[e^{\lambda \tr \mathbf{O}}]+\mathbb{E}_n^{--}[e^{-\lambda \tr \mathbf{O}}])
\]
by Markov's inequality. By Proposition \ref{Fred}, for any $\lambda\in\R$,
\[\mathbb{E}_n^{--}[e^{\lambda \tr \mathbf{O}}] = F_n^{ab}(-\i\lambda) = e^{\lambda^2/2} \det (1+Q_nK_{\lambda}^{--}Q_n). \]
We bound the Fredholm determinant using Plemelj's formula: by Lemma \ref{traceclassnorm}, if $\lvert\lambda\rvert<2n/e^{5/4}$, then
\[\det (1+Q_nK_{\lambda}^{--}Q_n) = \exp\Big( \sum_{j\geq 1} \frac{(-1)^{j+1}}{j} \tr (Q_nK_{\lambda}^{--}Q_n)^j \Big) \]
so 
\begin{align*}
\lvert \det (1+Q_nK_{\lambda}^{--}Q_n) \rvert &\leq \exp\Big( \sum_{j\geq 1} \frac{1}{j} \lvert \tr (Q_nK_{\lambda}^{--}Q_n)^j \rvert \Big) \\
&\leq \Big( 1-\frac{e^{\lambda^2/(2n+1)}}{1-\lambda^2/(4n^2)}\frac{\lambda^{2n}}{\Gamma(2n+1)} \Big)^{-1}
\end{align*}
by Lemma \ref{sumoftraces1}. We use the bound $(1-x)^{-1}<\exp (x/\sqrt{1-x})$, valid for $0<x<1$.
Stirling's inequality (\ref{Gamma}) and assuming that $\lvert \lambda \rvert < 2e^{-3/2}n < n$ gives
\begin{equation*}
\frac{e^{\lambda^2/(2n+1)}}{1-\lambda^2/(4n^2)}\frac{\lambda^{2n}}{\Gamma(2n+1)} < \frac{e^{n/2}}{1-e^{-3}} \frac{1}{2\sqrt{\pi n}}\Big(\frac{e\lambda}{2n}\Big)^{2n} < \frac{1}{2(1-e^{-3})\sqrt{\pi n e^{n}}} 
\end{equation*}
but also
\begin{align*}
\frac{e^{\lambda^2/(2n+1)}}{1-\lambda^2/(4n^2)}\frac{\lambda^{2n}}{\Gamma(2n+1)} &< \frac{e^{n/2}}{1-e^{-3}} \frac{\lambda^2}{2n(2n-1)} \frac{1}{2\sqrt{\pi (n-1)}}\Big(\frac{e\lambda}{2n-2}\Big)^{2n-2}\\
& < \frac{e\lambda^2}{4(1-e^{-3})n(2n-1)\sqrt{\pi(n-1)e^n}} \\
& < \frac{e\lambda^2}{4(1-e^{-3})\sqrt{\pi}n^{2.5}e^{n/2}}
\end{align*}
for any $n\geq 2$. To simplify the constants, we have for any $n\geq 1$, 
\[ \Big(1-\frac{1}{2(1-e^{-3})\sqrt{\pi n e^{n}}} \Big)^{-1/2}\frac{e}{\sqrt{\pi}(1-e^{-3})}< 2 . \]
Thus,
\begin{equation}\label{bound_det_concentration}
\lvert \det (1+Q_nK_{\lambda}^{--}Q_n) \rvert \leq \exp \Big(\frac{\lambda^2}{2n^{2.5}e^{n/2}}\Big)
\end{equation}
if $n\geq 2$ and
\[\mathbb{P}[\lvert \tr \mathbf{O} \rvert > \frac{L}{2}] \leq 2\exp \Big(-\frac{\lambda L}{2}+\frac{\lambda^2}{2}\Big(1+\frac{1}{n^{2.5}e^{n/2}}\Big)\Big) \]
for any $\lvert\lambda\rvert < 2e^{-3/2}n$. The upper bound attains its minimum when $\lambda = L/(2+2n^{-2.5}e^{-n/2})$, and becomes the desired bound, provided $L$ satisfies the given condition. The other groups are treated similarly.
\end{proof}

\subsection{Intermediate regime}

To get an estimate of  $F_n^{ab}$ in the intermediate range of $\xi$ we will make a change of variable in the integral expression of $F_n^{ab}$, similar to the one in Proposition \ref{prop:TB1} for the unitary case. The aim is to prove

\begin{prop}\label{lambda1}
Let $\Lambda\geq 2n/e^{5/4}$. For $n\geq 2$,
\begin{align}
\int_{2n/e^{5/4}\leq|\xi|\leq\Lambda} |F_n^{ab}(\xi)|^2 \d\xi \leq 2\Lambda e^{-2\beta_3^{ab}(n)n^2}
\end{align}
where
\begin{gather}\label{beta3}
\beta_3^{--}(n)= \frac{(1-\frac{1}{2n}-\frac{8\sqrt{2}}{\pi\sqrt{3(n-1)}n})^2}{2e^{5/2}+1+\frac{8\sqrt{2}}{\sqrt{3(n-1)}}}, \ \beta_3^{++}(n)= \frac{(1+\frac{1}{2n}-\frac{8\sqrt{2}}{\pi\sqrt{3n}n})^2}{2(1+\frac{1}{n})e^{5/2}+1+\frac{8\sqrt{2}}{\sqrt{3n}}}, \\
 \beta_3^{+-}(n)=\beta_3^{-+}(n)= \frac{(1-\frac{8\sqrt{2}}{\pi\sqrt{3(n-1/2)}n})^2}{2(1+\frac{1}{2n})e^{5/2}+1+\frac{8\sqrt{2}}{\sqrt{3(n-1/2)}}}. \nonumber
\end{gather}
\end{prop}

\begin{proof}
The integral expression of $F_n^{ab}$ is given by
$$ F_n^{ab}(\xi) = \frac{1}{Z_n^{ab}}\int_{[-1,1]^n} \prod_{1\leq j\leq n} (1-t_j)^a(1+t_j)^b \prod_{1\leq j<k\leq n} (t_j-t_k)^2 \prod_{1\leq j\leq n} e^{2\i \xi t_j}dt_j. $$
The integrand has an analytic continuation in $\C \setminus (-\infty,-1]\cup[1,\infty)$, so we can deform the contour by mapping the interval $(-1,1)$ to its image under $\gamma:t\mapsto t+i\nu h(t)/n$, where $h(t)=1-t^2$, and $\nu$ is a positive parameter that we will fix later. We then make a change of variables to recover the original contour:
\begin{align*}
F_n^{ab}(\xi) &= \frac{1}{Z_n^{ab}}\int_{[-1,1]^n} \prod_{1\leq j\leq n} (1-t_j-i\frac{\nu}{n}(1-t_j^2))^a(1+t_j+i\frac{\nu}{n}(1-t_j^2))^b \\
&\prod_{1\leq j<k\leq n} (t_j-t_k+i\frac{\nu}{n}(t_j^2-t_k^2))^2 \prod_{1\leq j\leq n} e^{2\i \xi(t_j+i\nu(1-t_j^2)/n)}(1-2i\frac{\nu}{n}t_j)dt \\
&= \frac{1}{Z_n^{ab}}\int_{[-1,1]^n} \prod_{1\leq j\leq n} (1-i\frac{\nu}{n}(1+t_j))^a(1+i\frac{\nu}{n}(1-t_j))^b \\
&\prod_{1\leq j<k\leq n} (1-i\frac{\nu}{n}(t_j+t_k))^2 \prod_{1\leq j\leq n} e^{2\i \xi(t_j+i\nu(1-t_j^2)/n)}(1-2i\frac{\nu}{n}t_j)\rho_n^{ab}(t)dt
\end{align*}
The absolute value of the first product in the last integral is less than one if $a=b=-1/2$. For the three other cases we use the inequality $|1+ix|\leq e^{x^2/2}$, valid for real $x$, and the fact that $|1\pm t_j|\leq 2$. This gives
$$
|(1-i\frac{\nu}{n}(1+t_j))^a(1+i\frac{\nu}{n}(1-t_j))^b| \leq \exp\Big(\frac{\nu^2}{n^2}(\mathbbm{1}\{a>0\}+\mathbbm{1}\{b>0\})\Big).
$$
Similarly,
$$
|1-i\frac{\nu}{n}(t_j+t_k)|^2 \leq \exp\Big(4\frac{\nu^2}{n^2}\Big), \qquad
|1-2i\frac{\nu}{n}t_j| \leq \exp\Big(2\frac{\nu^2}{n^2}\Big).
$$
We now have
\begin{align}\label{final_var}
|F_n^{ab}(\xi)| \leq   \exp\Big(\nu^2(2+\frac{1}{n}(\mathbbm{1}\{a>0\}+\mathbbm{1}\{b>0\}))\Big) \E_n^{ab}[\prod_{1\leq j\leq n} e^{-2\frac{\nu\xi}{n}(1-t_j^2)}].  
\end{align}
By Lemma \ref{Kurt}, the last factor can be expressed as the determinant of a Toeplitz + Hankel matrix to which we can apply the Basor-Ehrhardt formula, Proposition \ref{BE}. The symbol in the matrix is $a(\theta)= \exp (-2\nu\xi\sin^2(\theta)/n)$ so $a_+^{-1}\tilde{a_+}(\theta)=\exp(-i\nu\xi\sin(2\theta)/n) $. For the case $a=b=1/2$ we obtain
\begin{align}\label{Fredholm}
\E_n^{++}[\prod_{1\leq j\leq n} e^{-2\frac{\nu\xi}{n}(1-t_j^2)}] &= \exp\Big(-\nu\xi-\frac{\nu\xi}{2n}+\frac{\nu^2\xi^2}{4n^2}\Big)\det(1-Q_nH(t^{-1}a_ +^{-1}\tilde{a_+})Q_n).
\end{align}
To estimate this Fredholm determinant we use the inequality
$$ |\det(1+K) \rvert \leq e^{\|K\|_{\mathcal{J}_1}} $$ 
valid for any trace-class operator $K$ (Theorem VII.3.3 in \cite{GGK}), combined with \eqref{traceclass}. We therefore need bounds on every element of the Hankel matrix in (\ref{Fredholm}), i.e on the Fourier coefficients of $a_+^{-1}\tilde{a_+}$. We use
$$ |\widehat{(e^{\i g})}_k| \leq \frac{1}{k^2}(\|g'\|_{L_2(\T)}^2+\|g''\|_{L_1(\T)}) $$
which holds for any real valued function $g\in C^2(\mathbb{T})$ and follows from two partial integrations on $\widehat{(e^{\i g})}_k$, and obtain
$$ |\widehat{(a_+^{-1}\tilde{a_+})}_k| \leq \frac{1}{k^2} \Big( \frac{8\nu\xi}{\pi n}+\frac{2\nu^2\xi^2}{n^2} \Big) $$
so
\begin{align*}
\det(1- Q_nH(t^{-1}a_+^{-1}\tilde{a_+})Q_n) &\leq \exp\Big(\sum_{j\geq n}\sqrt{\sum_{k\geq n}|\widehat{(a_+^{-1}\tilde{a_+})}_{j+k+2}|^2}\Big) \\
&\leq \exp\Big(\Big( \frac{8\nu\xi}{\pi n}+\frac{2\nu^2\xi^2}{n^2} \Big)\int_n^\infty \Big(\int_n^\infty \frac{dx}{(x+y)^4}\Big)^{1/2} dy\Big) \\
&= \exp\Big(\sqrt{\frac{2}{3n}}\Big( \frac{8\nu\xi}{\pi n}+\frac{2\nu^2\xi^2}{n^2} \Big)\Big). 
\end{align*}
Combining all estimates obtained so far gives
\begin{align*}
|F_n^{++}(\xi)| &\leq \exp\Big(-\nu\xi\Big(1+\frac{1}{2n}- \frac{8\sqrt{2}}{\pi \sqrt{3n}n}\Big)+ \nu^2\Big( 2+\frac{2}{n}+\frac{\xi^2}{4n^2}+ \frac{2\sqrt{2}\xi^2}{\sqrt{3n}n^2}\Big) \Big).
\end{align*}

It remains to choose the parameter $\nu$ to minimize the upper bound. The exponential above is of the form $\exp(-\alpha\nu+\delta\nu^2)$, which attains its minimum $\exp(-\alpha^2/4\delta)$ for $\nu=\alpha/2\delta$, which is positive for large enough $n$ (recall that $\nu\geq0$). Explicitly, 
\begin{equation}\label{lastcase}
|F_n^{++}(\xi)| \leq \exp\Big(\frac{-\frac{1}{4}(1+\frac{1}{2n}-\frac{8\sqrt{2}}{\pi\sqrt{3n}n})^2\xi^2}{2+\frac{2}{n}+\frac{\xi^2}{4n^2}+\frac{2\sqrt{2}\xi^2}{\sqrt{3n}n^2}}\Big), \qquad n\geq2    
\end{equation}
so if $|\xi|\geq 2n/e^{5/4}$ (i.e. outside of the first regime),
\begin{align}\label{intermediateO}
|F_n^{++}(\xi)| \leq \exp(-\beta_3^{++}(n)n^2 ), \qquad n\geq2.
\end{align}
where 
\begin{align}
\beta_3^{++}(n) = \frac{(1+\frac{1}{2n}-\frac{8\sqrt{2}}{\pi\sqrt{3n}n})^2}{2(1+\frac{1}{n})e^{5/2}+1+\frac{8\sqrt{2}}{\sqrt{3n}}}.
\end{align}
This proves the proposition for $(a,b)=(1/2,1/2)$. The proof for the three other cases is essentially identical.
\end{proof}

\subsection{Large regime}

We use the same method as in the unitary case.

\begin{prop}\label{lambda2}
For any pair $(a,b)=(\pm\frac{1}{2},\pm\frac{1}{2})$ we have
\begin{equation}
\int_{|\xi| \geq\Lambda} \big|F_n^{ab}(\xi)\big|^2 \d\xi \leq \frac{2}{n-1} \frac{(2n)^n}{\Lambda^{n-1}}
\end{equation}
if $n\geq 5$ and $\Lambda\ge (n+1)/\sqrt{2}$, or if $n\geq 2$ and $\Lambda\ge \frac{7}{8}n$.
\end{prop}

\begin{proof}
Consider for example the case $a=b=-1/2$. By Lemma \ref{Kurt} and Hadamard's inequality,
$$ \big|  F_n^{--}(\xi) \big|^2   \leq  \prod_{j=1}^n \sum_{i=1}^n \big| (\widehat{e^{2\i \xi \cos}})_{j-i}+(\widehat{e^{2\i \xi \cos}})_{j+i}\big|^2  
= \prod_{j=1}^n \sum_{i=1}^n |I_{j-i}(2i\xi)+I_{j+i}(2i\xi)|^2 ,$$
so by (\ref{BesselJ}), if $n\geq 5$ and $|2\xi| \geq  \sqrt{2}(n+1)$, or if $n\geq 2$ and $|2\xi|\ge \frac{7}{4}n$,
\begin{equation}\label{largeO}
\big|  F_n^{--}(\xi) \big|^2   \leq \prod_{j=1}^n \sum_{i=1}^n \frac{4}{|2i\xi|} = \Big(\frac{2n}{\lvert\xi\rvert}\Big)^n.
\end{equation}
Thus, for any $\Lambda \ge   (n+1)/\sqrt{2} $, 
\[ \int_{|\xi| \ge \Lambda}\big|  F_n^{--}(\xi) \big|^2 \d\xi  \leq 2(2n)^n  \int_{\Lambda}^{+\infty} \frac{\d r}{r^{n}}  =  \frac{2}{n-1} \frac{(2n)^n}{\Lambda^{n-1}}.
\]
\end{proof}

\subsection{Proof of Theorem \ref{thm:1O}}

We prove Theorem \ref{thm:1O} by combining Proposition \ref{regime1}, \ref{lambda1} and \ref{lambda2}. First we want to choose the parameter $\Lambda$ appearing in Proposition \ref{lambda1} and \ref{lambda2} so that it minimizes the sum of the upper bounds obtained in these propositions, i.e. we want to minimize 
$$ \phi(\Lambda) = 2\Lambda e^{-2\beta_3^{ab}(n)n^2} + \frac{2}{n-1} \frac{(2n)^n}{\Lambda^{n-1}}.$$
The unique critical point is $\Lambda_0= 2ne^{2\beta_3^{ab}(n)n}$, which satisfies the assumption of Proposition \ref{lambda2} if $n\geq2$, and gives
\begin{align}\label{Flarge1}
\phi(\Lambda) \leq 4\frac{n^2}{n-1}e^{-2\beta_3^{ab}(n)(n^2-n)}
\end{align}
for $n\geq 1$. Next,
\begin{align}\label{Glarge}
\int_{|\xi|\geq 2n/e^{5/4}} e^{-\xi^2} \d\xi = \sqrt{\pi}\mathrm{erfc}\Big(\frac{2n}{e^{5/4}}\Big)< \frac{e^{5/4}}{4n}e^{-4n^2/e^{5/2}}   
\end{align}
by \href{https://dlmf.nist.gov/7.8}{(7.8.2)} in DLMF. Equations (\ref{Flarge1}) and (\ref{Glarge}) combined with Proposition \ref{regime1} give
\[
\|F_n^{ab}-\widehat{\phi_R} \|_2\leq  \frac{\beta_2(n)}{(2n)^{1/4}\sqrt{\Gamma(2n+1)}} + 2\frac{n}{\sqrt{n-1}}e^{-\beta_3^{ab}(n)(n^2-n)} + \frac{e^{5/8}}{2\sqrt{n}}e^{-2n^2/e^{5/2}}.
\]
Parseval's identity finishes the proof of the first bound. 

To simplify the bound we use Stirling's approximation (\ref{Gamma}) which gives
\[\sqrt{\Gamma(2n+1)}e^{-\beta_3^{ab}(n)(n^2-n)} \leq (2\pi e)^{1/4}\exp(-\beta_3^{ab}(n)(n^2-n)+(n+1/4)\log(2n)-n) \]
which is less than $\sqrt{2}(2n)^{-7/4}\sqrt{1-n^{-1}}$ if $n\geq 121$, by a numerical evaluation. Since
\[ e^{-2n^2/e^{5/2}} < e^{-\beta_3^{ab}(n)(n^2-n)}, \]
we obtain
\[
\|p_n^{--}-\phi_\R \|_2\leq  \frac{\beta_2(n)+\tfrac{2}{n}}{(2n)^{1/4}\sqrt{\Gamma(2n+1)}}
\]
for all $n\geq 121$. We proceed similarly for the other cases.

\subsection{Proof of the lower bound in Theorem \ref{thm:TVO}}

\begin{lemma}\label{trace2}
Let $d$ be the total number of eigenvalues. Then
\[ \lvert \tr Q_nK_{i\xi}^{ab}Q_n \lvert \geq \Big(1-\frac{ \xi ^2}{d^2} \Big) \frac{\xi^{d}}{\Gamma(d+1)} \]
Moreover, if $a=b$ then $\tr Q_nK_{i\xi}^{ab}Q_n$ is real, otherwise it is purely imaginary.
\end{lemma}

\begin{proof}
By Proposition \ref{Fred},
\[ \tr Q_nK_{\i \xi}^{--}Q_n = \sum_{j \geq n} (-1)^{j}I_{2j}(2\xi) = (-1)^n \sum_{j \geq 0} (I_{4j+2n}(2\xi)-I_{4j+2n+2}(2\xi)) \]

The modified Bessel function has the following integral representation DLMF \href{https://dlmf.nist.gov/10.32}{(10.32.2)}
\begin{equation}\label{BesselInt}
I_\nu(x) = \frac{(x/2)^\nu}{\sqrt{\pi}\Gamma(\nu+1/2)} \int_{-1}^1 (1-t^2)^{\nu-1/2} e^{xt} \mathrm{d}t, \quad \Re \nu > -\frac{1}{2}
\end{equation}
whence, 
\[ I_{4j+2n+2}(2\xi) < \frac{\xi^2}{4n^2} I_{4j+2n}(2\xi) \]
for all $j\geq 0$ and
\[ \lvert \tr Q_nK_{\i \xi}^{--}Q_n \rvert \geq \Big(1-\frac{\xi^2}{4n^2}\Big) \sum_{j \geq 0} I_{4j+2n}(2\xi) \geq  \Big(1-\frac{\xi^2}{4n^2}\Big) I_{2n}(2\xi) \]
The bound (\ref{Bessel-lowerbd}) on Bessel functions gives
\[ \lvert \tr Q_nK_{\i \xi}^{--}Q_n \rvert \geq \Big(1-\frac{\xi^2}{4n^2}\Big) \frac{\xi^{2n}}{\Gamma(2n+1)}, \]
which finishes the proof for the first case. The other cases are treated similarly.

\end{proof}

\begin{lemma}\label{sumoftraces2}
For any $|\xi| \leq e^{-2}d$, whith $d$ the total number of eigenvalues, 
\[\sum_{j\geq 2} \frac{1}{j} \lvert \tr (Q_nK_{\i \xi}^{ab}Q_n)^j \rvert \leq \frac{1}{\sqrt{d/2}e^{d/2}} \frac{\xi^{d}}{\Gamma(d+1)}. \]
\end{lemma}

\begin{proof}
From the proof of Lemma \ref{sumoftraces1},
\[ \lvert \tr (Q_nK_{\i \xi}^{--}Q_n)^j \rvert \leq \beta_1(2n)^j \]
so
\[
\sum_{j\geq 2} \frac{1}{j} \lvert \tr (Q_nK_{\i \xi}^{--}Q_n)^j \rvert \leq \sum_{j\geq 1} \frac{1}{j}  \big(\beta_1(2n)\big)^{j+1} = - \log (1-\beta_1(2n)) \beta_1(2n). \]
We use the bound $-\log(1-x)<x/\sqrt{1-x}$.
Stirling's inequality (\ref{Gamma}) and the assumption that $\lvert \xi \rvert < 2e^{-2}n < n$ give
\begin{equation}\label{t1}
\beta_1(2n) = \frac{e^{\xi^2/(2n+1)}}{1-\xi^2/(4n^2)}\frac{\xi^{2n}}{\Gamma(2n+1)} < \frac{e^{n/2}}{1-e^{-4}} \frac{1}{2\sqrt{\pi n}}\Big(\frac{e\xi}{2n}\Big)^{2n} < \frac{1}{2(1-e^{-4})\sqrt{\pi n e^{3n}}} 
\end{equation}
but also
\begin{equation}\label{t2}
\Big( \frac{e^{\xi^2/(2n+1)}}{1-\xi^2/(4n^2)} \Big)^2 \frac{\xi^{2n}}{\Gamma(2n+1)} < \frac{e^{n}}{(1-e^{-4})^2} \frac{1}{2\sqrt{\pi n}}\Big(\frac{e\xi}{2n}\Big)^{2n} < \frac{1}{2(1-e^{-4})^2\sqrt{\pi n }e^{n}}.  
\end{equation}
Thus 
\[ \sum_{j\geq 2} \frac{1}{j} \lvert \tr (Q_nK_{\i \xi}^{--}Q_n)^j \rvert \leq \Big(1-\frac{1}{2(1-e^{-4})\sqrt{\pi n e^{3n}}}\Big)^{-1/2} \frac{1}{2(1-e^{-4})^2\sqrt{\pi n }e^{n}} \frac{\xi^{2n}}{\Gamma(2n+1)} \]
which, after some simplification, gives the inequality in the case $a=b=-1/2$.
\end{proof}

\begin{prop}\label{lowerbd}
If $|\xi| \leq e^{-2}d$,
\[
\lvert \det(1+Q_nK_{\i \xi}^{ab}Q_n)-1 \rvert \geq \gamma(d/2) \frac{\xi^{d}}{\Gamma(d+1)}
\]
where 
\[ \gamma(x) = (1-e^{-x})(1-e^{-4}-e^{-x}x^{-1/2}), \ x\in\R^+.\]

\end{prop}

\begin{proof}
We use Plemelj's formula, i.e.
\[
 \det(1+K_{i\xi}^{ab}Q_n) = \exp\Big (\sum_{j\geq 1} \frac{(-1)^{j+1}}{j} \tr(Q_nK_{\i\xi}^{ab}Q_n)^j \Big)  
\]
if $\|K_{i\xi}^{ab}Q_n\|_{J_1} < 1 $. By Lemma \ref{traceclassnorm}, this holds for $\lvert \xi \rvert < 2e^{-5/4}n$.
Also, if $|x| \leq 1$ then $|x|(1-|x|) \leq |e^x-1| $. Thus, if $\lvert \Re \delta \rvert < 1$, where $\delta=\sum_{j\geq 1} \frac{(-1)^{j+1}}{j} \tr(Q_nK_{\i\xi}^{ab}Q_n)^j$, then
\[ \lvert \det(1+K_{i\xi}^{ab}Q_n)-1 \rvert \geq \lvert \lvert\det(1+K_{i\xi}^{ab}Q_n)\rvert-1 \rvert = \lvert e^{\Re \delta} -1 \rvert \geq \lvert \Re \delta\rvert (1-\lvert \Re \delta\rvert) \geq  \lvert \Re \delta\rvert - \lvert\delta\rvert^2. \]
Now, if $a=b=-1/2$, then Lemma \ref{sumoftraces1} gives
\begin{align*}
\lvert \delta \rvert^2 &\leq \Big( \sum_{j\geq 1} \frac{1}{j} \lvert \tr(Q_nK_{\i\xi}^{ab}Q_n)^j \rvert \Big)^2 \\
&\leq  \log^2 (1-\beta_1(2n)) \\
&\leq \Big( 2\sqrt{\pi}(1-e^{-4})(1-e^{-4}-1/(2\sqrt{\pi n e^{3n}}))\Big)^{-1}\frac{1}{\sqrt{n}e^n} \frac{\xi^{2n}}{\Gamma(2n+1)} \\
&\leq \frac{1}{\sqrt{n}e^n} \frac{\xi^{2n}}{\Gamma(2n+1)}
\end{align*}
where we used (\ref{t1}), (\ref{t2}), and the fact that $-\log(1-x)<x/\sqrt{1-x}$ for $0<x<1$. Next,
\[ \lvert \Re \delta \rvert = \lvert \tr Q_nK_{\i \xi}^{--}Q_n + \Re \sum_{j\geq 2} \frac{(-1)^{j+1}}{j} \tr (Q_nK_{\i \xi}^{--}Q_n)^j \rvert \geq \lvert \tr Q_nK_{\i \xi}^{--}Q_n \rvert -\sum_{j\geq 2} \frac{1}{j} \lvert \tr (Q_nK_{\i \xi}^{--}Q_n)^j \rvert  \]
so combining Lemma \ref{trace2} with Lemma \ref{sumoftraces2} gives
\[ \lvert \Re \delta\rvert \geq \Big(1-\frac{\lvert \xi \rvert^2}{4n^2} -\frac{1}{\sqrt{n}e^n}\Big) \frac{\xi^{2n}}{\Gamma(2n+1)} \geq \Big(1-e^{-4} -\frac{1}{\sqrt{n}e^n}\Big) \frac{\xi^{2n}}{\Gamma(2n+1)} \]
if $|\xi| \leq 2e^{-2}n$. This gives the desired lower bound, since $1-e^{-4}-2e^{-n}n^{-1/2} \geq \gamma(n)$, for all $n\in\N$.

The proof for the case $a=b=1/2$ is identical, except that $n$ is replaced by $n+1$.
For the last two cases, observe that
\begin{align*}
\lvert \det(1+K_{i\xi}^{ab}Q_n)-1\rvert &= \Big( \Big(e^{\Re \delta}\cos(\Im \delta)-1\Big)^2+e^{2\Re{\delta}}\sin^2(\Im \delta)\Big)^{1/2} \\
&=\Big((e^{\Re \delta}-1)^2+2e^{\Re \delta}(1-\cos(\Im \delta))\Big)^{1/2} \\
&\geq \Big(2e^{\Re \delta}(1-\cos(\Im \delta))\Big)^{1/2} \\
&\geq e^{\Re \delta/2}  \sqrt{1-\lvert \Im \delta \rvert^2}\lvert \Im \delta \rvert
\end{align*}
To obtain a lower bound on $\Re \delta$ we simply write
\[
\Re \delta = \Re\Big( \sum_{j\geq 1} \frac{(-1)^{j+1}}{j} \tr(Q_nK_{\i\xi}^{ab}Q_n)^j \Big) \geq - \sum_{j\geq 1} \frac{1}{j} \lvert \tr(Q_nK_{\i\xi}^{ab}Q_n)^j \rvert \geq \log (1-\beta_1(2n+1)) 
\]
if $a=-b$, by Lemma \ref{sumoftraces1}. Thus,
\[ e^{\Re \delta/2} \geq \sqrt{1-\beta_1(2n+1)} \geq \Big(1-\frac{1}{2(e-e^{-3})\sqrt{\pi (n+1/2)e^{3n}}}\Big)^{1/2} \geq (1-e^{-n-1/2})^{1/2} \]
by Stirling's approximation (\ref{Gamma}) and the assumption that $\lvert\xi\rvert\leq 2e^{-2}n$. We can use the same lower bound for $\sqrt{1-\lvert \Im \delta \rvert^2}$: by Lemma \ref{sumoftraces1},
\begin{align*}
|\Im \delta| \leq \sum_{j\geq 1} \frac{1}{j} \lvert \tr(Q_nK_{\i\xi}^{ab}Q_n)^j\rvert &\leq - \log \Big( 1-\beta_1(2n+1) \Big) < \sqrt{\beta_1(2n+1)}.
\end{align*}
Finally,
\[ \lvert \Im \delta \rvert = \lvert \tr K_{i\xi}^{ab}Q_n + \Im \sum_{j\geq 2} \frac{(-1)^{j+1}}{j} \tr(Q_nK_{\i\xi}^{ab}Q_n)^j \rvert \geq \lvert \tr K_{i\xi}^{ab}Q_n \rvert -\sum_{j\geq 2} \frac{1}{j} \lvert \tr(Q_nK_{\i\xi}^{ab}Q_n)^j \rvert  \]
so by Lemma \ref{trace2} and Lemma \ref{sumoftraces2},
\[ \lvert \Im \delta \rvert \geq \Big(1-\frac{\lvert \xi \rvert^2}{4n^2} -\frac{1}{\sqrt{n+1/2}e^{n+1/2}}\Big) \frac{\xi^{2n+1}}{\Gamma(2n+2)} \geq \Big(1-e^{-4} - \frac{1}{\sqrt{n+1/2}e^{n+1/2}} \Big) \frac{\xi^{2n+1}}{\Gamma(2n+2)} \]
which proves the statement.

\end{proof}

\begin{prop}\label{prop:lowerboundO}
Let $n\geq 31$. Then
\[
\| p_n^{ab} - \phi_\R\|_{L^2}^2 \geq \frac{2\gamma(d/2)^2}{\sqrt{3\pi}e}\frac{1}{\sqrt{d}\Gamma(d+1)}
\]
where $d$ is the total number of eigenvalues. 
\end{prop}

\begin{proof}
By Proposition \ref{Fred},
\begin{align*}
 \| F_n^{--} - \widehat{\phi_\R} \|_{L^2}^2 &= \int_{\R} e^{-\xi^2} \lvert \det(1+Q_nK_{\i \xi}^{--}Q_n)-1 \rvert^2 \d \xi \\
 &\geq \int_{\lvert\xi\rvert<2n/e^2} e^{-\xi^2} \lvert \det(1+Q_nK_{\i \xi}^{--}Q_n)-1 \rvert^2 \d \xi,
\end{align*}
so Proposition \ref{lowerbd} gives
\begin{align}\label{s1}
 \| F_n^{--} - \widehat{\phi_\R} \|_{L^2}^2  \geq \frac{\gamma(n)^2}{\Gamma(2n+1)^2} \int_{\lvert\xi\rvert<2n/e^2} e^{-\xi^2} \xi^{4n}  \d \xi.
\end{align}
We estimate the integral using Laplace method. First we make the change of variable $x=\xi/\sqrt{n}$:
\begin{equation}\label{s2}
\int_{\lvert\xi\rvert<2n/e^2} e^{-\xi^2} \xi^{4n}  \d \xi = 2n^{2n+1/2} \int_0^{2\sqrt{n}/e^2} e^{nf(x)} \d x
\end{equation}
where $f(x)= -x^2+4\log x$ has a unique critical point at $x=\sqrt{2}$. Thus, for $\delta>0$
\[ \lvert f(x)-f(\sqrt{2}) \rvert \leq \sup_{x\in [\sqrt{2}-\delta,\sqrt{2}+\delta]} \frac{\lvert f''(x)\rvert\delta^2}{2} = \big(1+\frac{2}{(\sqrt{2}-\delta)^2}\big)\delta^2 \]
which is less than $3\delta^2$ if $\delta\leq \sqrt{2}-1$. We obtain, for $2e^{-2}\sqrt{n} \geq \sqrt{2}+\delta$,
\[ \int_0^{2\sqrt{n}/e^2} e^{nf(x)} \d x \geq \int_{\sqrt{2}-\delta}^{\sqrt{2}+\delta} e^{nf(x)} \d x \geq 2\delta e^{n(f(\sqrt{2})-3\delta^2)} = 2\delta \big(\frac{2}{e}\big)^{2n}e^{-3n\delta^2} \] 
which is maximized at $\delta = (6n)^{-1/2}$. We insert this value in the above expression which in turn is inserted in (\ref{s2}) and (\ref{s1}), giving
\[ \| F_n^{--} - \widehat{\phi_\R} \|_{L^2}^2  \geq \frac{2\sqrt{2}}{\sqrt{3e}}\frac{\gamma(n)^2}{\Gamma(2n+1)^2}\Big(\frac{2n}{e}\Big)^{2n}, \]
if $n\geq 31$. Stirling's approximation (\ref{Gamma}) gives
\begin{equation}\label{lower2O}
\| F_n^{--} - \widehat{\phi_\R} \|_{L^2}^2  \geq \frac{\sqrt{2}}{\sqrt{3\pi}e}\frac{\gamma(n)^2}{\sqrt{n}\Gamma(2n+1)}.
\end{equation}
The proof for the other cases is identical except that we replace $2n$ with $2n+1$ or $2n+2$.
\end{proof}

Now, observe that 
\begin{align*}
\|F_n^{--} - \widehat{\phi_\R} \|_{L^2}^2 &= \|p_n^{--} - \phi_\R\|_{L^2}^2 \\
&\leq \| p_n^{--} - \phi_\R\|_{L^\infty(\R)} \| p_n^{--} - \phi_\R\|_{L^1(\R)}.
\end{align*}
Fourier's inversion theorem gives
\[ \| p_n^{--} - \phi_\R\|_{L^\infty(\R)} \leq \frac{1}{\sqrt{2\pi}} \| F_n^{--} - \widehat{\phi_\R} \|_{L^1(\R)} \]
so to obtain the lower bound on $\| p_n^{--} - \phi_\R\|_{L^1(\R)}$ it remains to estimate the $L^1(\R$) distance between the characteristic functions. 
By \eqref{freddistance} we have, for any $n\in\N$,
\[
 \int_{|\xi|< 2n e^{-5/4}} |F_n^{--}(\xi)-e^{-\xi^2/2}| \d\xi \\
\leq  \frac{\sqrt{1-e^{-5/2}}}{(1-e^{-5/2}-1/(2\sqrt{\pi n}))^{3/2}} \frac{1}{\Gamma(2n+1)} \int_{\R} e^{-\tfrac{\xi^2}{2}(1-\frac{2}{2n+1})}|\xi|^{2n} \d\xi
\]
where
\[
\int_{\R} e^{-\tfrac{\xi^2}{2}(1-\frac{2}{2n+1})}|\xi|^{2n} \d\xi = \Big(\frac{1}{2}-\frac{1}{2n+1}\Big)^{-n-\frac{1}{2}}\Gamma(n+1/2)   
\]
and
\[ \Gamma(n+1/2)< \frac{\pi^{1/4}}{ n^{1/4}2^n} \sqrt{\Gamma(2n+1)} \]
by combining the duplication formula for the Gamma function with Gautschi's inequality DLMF (\href{https://dlmf.nist.gov/5.5}{(5.5.5)} and \href{https://dlmf.nist.gov/5.6}{(5.6.4)}). Thus, 
\begin{equation} 
\int_{|\xi|<2ne^{-5/4}} |F_n^{--}(\xi)-e^{-\xi^2/2}| \d\xi
 \leq \frac{\sqrt{2} e^{1+1/n}}{(n/\pi)^{1/4}} \frac{\beta_2(n)}{\sqrt{\Gamma(2n+1)}}   . 
\end{equation}
By integrating the estimates \eqref{intermediateO} and \eqref{largeO}, we obtain for $n\ge 5$  and any $\Lambda \ge (n+1)/\sqrt{2}$, 
\[\begin{aligned} 
\int_{2e^{-5/4} n \leq |\xi| \leq  \Lambda}
\big| F_n^{--}(\xi) \big|\d\xi   & \le 2\Lambda e^{- \beta_3^{--}(n)n^{2}} ,  \\
\int_{|\xi| \ge \Lambda}\big|  F_n^{--}(\xi) \big| \d\xi & \le
\frac{2 }{n/2-1} \frac{(2n)^{n/2}}{\Lambda^{n/2-1}} . 
\end{aligned}\]
Choosing $\Lambda = 4n^2(n-2)^{-1} \exp(-\beta_3^{--}n^2(1-2/n))$, the unique critical point, gives
\[\int_{2e^{-5/4} n \leq |\xi|}
\big| F_n^{--}(\xi) \big|\d\xi   \leq \frac{4n^2}{n-2}e^{-\beta_3^{--}(n)n^2(1-2/n)} \]
and since
\[\int_{2e^{-5/4} n \leq |\xi|}
 e^{-\xi^2/2} \d\xi   \leq \frac{e^{5/4}}{\sqrt{2}n}e^{-2e^{-5/2}n^2} < \frac{n^2}{n-2}e^{-\beta_3^{--}(n)n^2(1-2/n)} \]
we obtain
\[ \| F_n^{--} - \widehat{\phi_\R} \|_{L^1(\R)} \leq \frac{\sqrt{2} e^{1+1/n}}{(n/\pi)^{1/4}} \frac{\beta_2(n)}{\sqrt{\Gamma(2n+1)}} + \frac{5n^2}{n-2}e^{-\beta_3^{--}(n)n^2(1-2/n)} \]
for any $n\geq 5$. If $n\geq 123$,
\[ \frac{1}{\sqrt{2\pi}}\| F_n^{--} - \widehat{\phi_\R} \|_{L^1(\R)} \leq  \frac{e^{1+1/n}\beta_2(n)+\tfrac{1}{n}}{(2n)^{1/4}\sqrt{\Gamma(2n+1)}} \]
which combined with Proposition \ref{prop:lowerboundO} gives
\begin{align*}
 \| p_n^{--} - \phi_\R\|_{L^1(\R)} \geq \frac{\sqrt{2}}{\sqrt{3\pi}e}\frac{\gamma(n)^2}{e^{1+1/n}\beta_2(n)+\tfrac{1}{n}} \frac{1}{(2n)^{1/4}\sqrt{\Gamma(2n+1)}}.
\end{align*}
A numerical evaluation of the front factor gives the lower bound in Theorem \ref{thm:TVO} in the case $a=b=-1/2$, i.e. matrices with an even number of rows and with positive determinant. The other cases are treated similarly.

\subsection{Exact asymptotics: proof of Theorem \ref{prop:asympO}}
We consider the case $a=b=-1/2$.
In the proof of Proposition \ref{regime1} we showed that 
\[ \lvert\det (1+Q_nK_{\i \xi}^{--}Q_n)-1\rvert  < \lvert\delta\rvert e^{\lvert\delta\rvert} < \Big( 1- \frac{e^{\xi^2/(2n+1)}}{1-\xi^2/(4n^2)}\frac{\xi^{2n}}{\Gamma(2n+1)}\Big)^{-1} |\delta| \]
where $\delta = \sum_{j\geq 1} \frac{(-1)^{j+1}}{j} \tr (Q_nK_{\i g}^{ab}Q_n)^j$, so if $|\xi| < 2n/e^2$,
\[ \lvert\det (1+Q_nK_{\i \xi}^{--}Q_n)-1\rvert  < ( 1-e^{-n})^{-1} \Big( |\tr Q_nK_{\i \xi}^{--}Q_n| + \sum_{j\geq 2} \frac{1}{j} |\tr (Q_nK_{\i \xi}^{--}Q_n)^j| \Big). \]
On the other hand, by the proof of Proposition \ref{lowerbd},
\[ \lvert \det(1+Q_nK_{i\xi}^{--}Q_n)-1 \rvert \geq  \lvert \Re \delta\rvert - \lvert\delta\rvert^2 > |\tr Q_nK_{\i \xi}^{--}Q_n| - \sum_{j\geq 2} \frac{1}{j} |\tr (Q_nK_{\i \xi}^{--}Q_n)^j|- \lvert\delta\rvert^2. \]
By Lemma \ref{sumoftraces1} and Lemma \ref{sumoftraces2}, as $n\to\infty$,
\[
|\delta|^2 = o(1)\frac{\xi^{2n}}{\Gamma(2n+1)}, \quad
\sum_{j\geq 2} \frac{1}{j} |\tr (Q_nK_{\i \xi}^{--}Q_n)^j| = o(1)\frac{\xi^{2n}}{\Gamma(2n+1)}
\]
uniformly in $|\xi| < 2n/e^2$. But by combining \eqref{boundontraceofpower} with Lemma \ref{trace2}, 
\[ \frac{\xi^{2n}}{\Gamma(2n+1)} = O(1) \tr (Q_nK_{\i \xi}^{--}Q_n)  \]
as $n\to\infty$. Thus
\[ \lvert \det(1+Q_nK_{\i\xi}^{--}Q_n)-1 \rvert = (1+o(1))|\tr Q_nK_{\i \xi}^{--}Q_n|  \]
as $n\to\infty$, uniformly in $|\xi| < 2n/e^2$. Furthermore,
\begin{align*}
|\tr Q_nK_{\i \xi}^{--}Q_n -(-1)^nI_{2n}(2\xi)| &= |\sum_{j\geq 0} (-1)^{j+n+1} I_{2j+2n+2}(2\xi)| \\
& \leq \sum_{j\geq 0} I_{2j+2n+2}(2\xi) \\
&\leq \frac{\xi^2}{4n^2} \frac{I_{2n}(2\xi)}{1-\xi^2/(4n^2)}
\end{align*}
by using \eqref{Ideac} repeatedly. Thus $\tr Q_nK_{\i \xi}^{--}Q_n =((-1)^n+O(\epsilon^2))I_{2n}(2\xi)$ if $\xi= \epsilon n$, $\epsilon < 2e^{-2}$. Finally, by Proposition \ref{lambda1} and \ref{lambda2},
\[ \int_{|\xi|> \epsilon n} |F_n^{--}(\xi)|^2 \d\xi \leq O(e^{-\delta n^2}) \]
for some $\delta >0$, which leads us to 
\[ \|F_n^{--} - \widehat{\phi_R} \|_ 2^2 = (1+O(\epsilon^2))\int_{|\xi|< \epsilon n} I_{2n}(2\xi)^2e^{-\xi^2} \d\xi + O(e^{-\delta n^2}). \]

We can estimate the integral without having to restrict it to $|\xi| < \epsilon n$: by \eqref{Bessel-upperbd},
\[ \int_{|\xi| \geq \epsilon n } I_{2n}^2(2\xi) e^{-\xi^2} \d\xi \leq \frac{1}{\Gamma(2n+1)^2} \int_{|\xi| \geq \epsilon n} \xi^{4n} e^{-\xi^2(1-\frac{2}{2n+1})} \d\xi, \]
so by Markov's inequality, followed by the change of variables $r= \xi^2(1-\frac{2}{2n+1})$,
\begin{align*}
\int_{|\xi| \geq \epsilon n } I_{2n}^2(2\xi) e^{-\xi^2} \d\xi &\leq \frac{1}{\epsilon n\Gamma(2n+1)^2} \int_\R |\xi|^{4n+1} e^{-\xi^2(1-\frac{2}{2n+1})} \d\xi \\
&\leq \frac{1}{\epsilon n\Gamma(2n+1)}\Big(1-\frac{2}{2n+1}\Big)^{-2n-1}
\end{align*}
This implies that
\[ \|F_n^{--} - \widehat{\phi_R} \|_ 2^2 = (1+O(\epsilon^2))\int_\R I_{2n}(2\xi)^2e^{-\xi^2} \d\xi + \frac{O_\epsilon(n^{-1})}{\Gamma(2n+1)}. \]
Combining this last equality with the following lemma proves Theorem \ref{prop:asympO}.

\begin{lemma}
\[ \lim_{\nu\to\infty} \sqrt{\nu}\Gamma(\nu+1) \int_\R I_\nu^2(2\xi) e^{-\xi^2} \d\xi = \frac{e^2}{\sqrt{2}}. \]
\end{lemma}
\begin{proof}
\begin{align*}
\int_\R I_\nu^2(2\xi) e^{-\xi^2} \d\xi &= \sqrt{\nu} \int_0^\infty I_\nu^2(2\sqrt{\nu r}) \frac{e^{-\nu r}}{\sqrt{r}} \d r \\
&= \frac{\nu^{\nu+1/2}}{\Gamma(\nu+1)^2} \int_0^\infty \Big( \sum_{j\geq 0}\frac{(\nu r)^j\Gamma(v+1)}{j!\Gamma(j+\nu+1)}\Big)^2 r^{\nu-1/2} e^{-\nu r} \d r \\
&\leq \frac{\nu^{\nu+1/2}}{\Gamma(\nu+1)^2} \int_0^\infty e^{2r} r^{\nu-1/2} e^{-\nu r} \d r
\end{align*}
since $\frac{\nu ^j\Gamma(v+1)}{\Gamma(j+\nu+1)} \leq 1$ for all $j\geq 0$. A direct application of Laplace method shows that as $\nu \to \infty$,
\[ \int_0^\infty e^{2r} r^{\nu-1/2} e^{-\nu r} \d r \sim \frac{\sqrt{\pi}}{\sqrt{\nu}} e^{2-\nu} \]
so 
\[ \limsup_{\nu\to\infty} \frac{e^\nu\Gamma(\nu+1)^2}{\nu^\nu}\int_\R I_\nu^2(2\xi) e^{-\xi^2} \d\xi \leq \sqrt{\pi}e^2. \]
For the lower bound, we write
\[\int_\R I_\nu^2(2\xi) e^{-\xi^2} \d\xi 
\geq \frac{\nu^{\nu+1/2}}{\Gamma(\nu+1)^2} \int_0^\infty \Big( \sum_{j\leq\kappa} \frac{(\nu r)^j\Gamma(v+1)}{j!\Gamma(j+\nu+1)}\Big)^2 r^{\nu-1/2} e^{-\nu r} \d r  \]
where $\kappa \in\N$ and proceed as in section \ref{sec:asymp}. We obtain
\[ \lim_{\nu\to\infty} \frac{e^\nu\Gamma(\nu+1)^2}{\nu^\nu}\int_\R I_\nu^2(2\xi) e^{-\xi^2} \d\xi = \sqrt{\pi}e^2 \]
which is equivalent to the statement of the lemma.
\end{proof}

\section{$\tr \mathbf{U}^m$ and $\tr \mathbf{O}^m$: proof of Theorems \ref{thm:inter} and \ref{thm:interO}} 
\label{sec:inter}

The last section is devoted to traces of higher powers of the matrices. We use the results of Rains who showed in \cite{Rains} that the eigenvalues of $\u^m$ and $\mathbf{O}^m$ have the same distribution as those of $m$ independent matrices of smaller size. This implies that the characteristic function of $\tr \u^m$ and $\tr \mathbf{O}^m$ factorizes into a product of $m$ characteristic functions of $\tr \u$ and $\tr \mathbf{O}$, where the size of $\u$ and $\mathbf{O}$ is equal or close to $N = \floor{n/m}$. We can therefore use our estimates on $F_n$ and $F_n^{ab}$ from the first two sections to obtain estimates on the characteristic function of $\tr \u^m$ and $\tr \mathbf{O}^m$, and proceed with the proof as before.


\begin{thm}[Rains, \cite{Rains}]\label{Rains}
Fix integers $n\geq 1$, $m\geq1$ and denote by $\sim$ equality of eigenvalue distributions. We have that
\begin{align}\label{powerunitary}
U(n)^m \sim \bigoplus_{0\leq i < m} U\Big(\ceil[\Big]{\frac{n-i}{m}}\Big).
\end{align}
If $m$ is odd,
\begin{align}
O^{\pm}(2n)^m &\sim O^{\pm}\Big(2\ceil[\Big]{\frac{n}{m}}\Big) \bigoplus_{1\leq i\leq(m-1)/2} O^+\Big(\ceil[\Big]{\frac{n-i}{m}}+\ceil[\Big]{\frac{n+i}{m}}\Big) \oplus O^-\Big(\ceil[\Big]{\frac{n-i}{m}}+\ceil[\Big]{\frac{n+i}{m}}\Big)  \label{even_odd}\\
O^{\pm}(2n+1)^m &\sim \bigoplus_{0\leq i<(m-1)/2} O^+\Big(\ceil[\Big]{\frac{n-i}{m}}+\ceil[\Big]{\frac{n+1+i}{m}}\Big) \oplus O^-\Big(\ceil[\Big]{\frac{n-i}{m}}+\ceil[\Big]{\frac{n+1+i}{m}}\Big)  \\
&\oplus O^{\pm}\Big(2\ceil[\Big]{\frac{n-(m-1)/2}{m}}+1\Big) \nonumber
\end{align}
and if $m$ is even,
\begin{align}
O^{\pm}(2n)^m &\sim O^{\pm}\Big(2\ceil[\Big]{\frac{n}{m}}\Big) \bigoplus_{1\leq i\leq m/2} O^+\Big(\ceil[\Big]{\frac{n-i}{m}}+\ceil[\Big]{\frac{n+i}{m}}\Big) \oplus O^-\Big(\ceil[\Big]{\frac{n-i}{m}}+\ceil[\Big]{\frac{n+i}{m}}\Big) \\ \nonumber
&\oplus O^{\mp}\Big(2\ceil[\Big]{\frac{n-m/2}{m}}+1 \Big) \\
O^{\pm}(2n+1)^m &\sim \bigoplus_{0\leq i< m/2} O^+\Big(\ceil[\Big]{\frac{n-i}{m}}+\ceil[\Big]{\frac{n+1+i}{m}}\Big) \oplus O^-\Big(\ceil[\Big]{\frac{n-i}{m}}+\ceil[\Big]{\frac{n+1+i}{m}}\Big).
\end{align}
\end{thm}

\subsection{Unitary group: Proof of Theorem~\ref{thm:inter}}
Combining the results from Section~\ref{sec:U} with Theorem \ref{Rains} allows us to obtain a rate of convergence for $\tr \mathbf{U}^m$, for all $m\in \N \cap [2,n]$ where $\mathbf{U} \in \mathbb{U}(n)$ is Haar-distributed.
Recall that  $p_{n,m}$ denotes the p.d.f.~of the complex-valued random variable $\frac{\tr\u^m}{\sqrt m}$. 

\begin{prop}\label{prop:inter}
 Let $n,m\in\N$ and $N = \lfloor n/m \rfloor +1$. Then,
\[
\big\| p_{n,m} - \widehat{\phi_\C} \big\|_{L^2}^2 
\le \frac1{2\pi}\bigg(  \frac{0.04/\sqrt{N}}{\Gamma(N+1)^2m^{2(N-1)}} 
+ 2 \exp\big( - 0.8 N^2 m^{1-1/N} \big) 
+ \big(2+ \tfrac{1}{n}\big) m N^2  \exp(-0.02(n-2)N)  \bigg) .
\]
\end{prop}

Before going into the proof, we would like to simplify this bound by assuming that the dimension $n$ of the random matrix $\mathbf{U}$ is sufficiently large. 
First observe that by \eqref{Gamma},
\[
\frac{1/\sqrt{N}}{\Gamma(N+1)^2m^{2N}} 
\ge \frac{n^{-2N}}{2\pi e N^{3/2}} 
\]
and for all $N\ge 2$ and  $n\ge 700$, 
\begin{equation} \label{numbd1}
N^{7/2} \exp(-0.02(n-2)N)  \le 10.5 \exp(-0.019nN) \le 8.9 \, n^{-2N} .
\end{equation}
In this regime, for $m\ge 2$, 
\[
\big(2+ \tfrac{1}{n}\big) m N^2  \exp(-0.02(n-2)N)  
\le 9  \frac{m^2  n^{-2N}}{N^{3/2}} 
\le  \frac{18\pi e /\sqrt{N}}{\Gamma(N+1)^2m^{2(N-1)}}  . 
\]

Similarly for $N\ge 2$
\[
\exp\big( - 0.8 N^2 m^{1-1/N} \big) \le \exp\big( - 0.8 N^{3/2} \sqrt{n} \big)  , 
\]
and for $n\ge 700$
\[
N^{3/2} \exp\big( - 0.8 N^{3/2} \sqrt{n} \big)  \le  2 \exp\big( - 1.2 N \sqrt{n} \big) \le 2 e^{-18} \, n^{-2N} .
\]
Hence, in this regime, for $m\ge 2$, 
\begin{equation} \label{numbd2}
2 \exp\big( - 0.8 N^2 m^{1-1/N} \big) \le 2 e^{-18} \,  \frac{n^{-2N}}{N^{3/2}} 
\le  \frac{e^{-17}\pi /\sqrt{N}}{\Gamma(N+1)^2m^{2(N-1)}}  .
\end{equation}

If we combine the estimate from Proposition~\ref{prop:inter} with the numerical bounds \eqref{numbd1}--\eqref{numbd2}, we obtain for $m\ge 2$ and $n\ge 700$, 
\begin{equation} \label{trnm}
\big\| p_{n,m} - \widehat{\phi_\C} \big\|_{L^2}^2 
\le \frac{77/\pi\sqrt{ \lfloor n/m \rfloor +1}}{\Gamma( \lfloor n/m \rfloor +2)^2 m^{2\lfloor n/m \rfloor}}. 
\end{equation}

To complete the proof of Theorem~\ref{thm:inter}, note that by \eqref{CS} and using the upper-bound \eqref{trnm}, we obtain for any $n\ge 700$, 
\[
 \| p_{n,m} - \phi_\C \|_{L^1} \le 
2 \bigg( \frac{5 L }{(N+1)^{1/4}\Gamma(N+2) m^{N}} + 3 e^{-L^2} \bigg),
\]
with $N=  \lfloor n/m \rfloor$. To minimize the right-hand side, we choose $L =  \sqrt{ \log\big( \Gamma(N+2)m^{N}\big) }$.
Hence, since $L^2 \le (N+1)\log(n)$,
 we  conclude that  for $n\ge 700$, 
\[\begin{aligned}
 \| p_n - \phi_\C \|_{L^1} 
 &\le  2 \left( 5+\frac{3}{(N+1)^{1/4}\sqrt{\log n}} \right) \frac{(N+1)^{1/4}\sqrt{\log n}}{\Gamma(N+2) m^N}  \\
 & \le  12 \frac{(N+1)^{1/4}\sqrt{\log n}}{\Gamma(N+2) m^N} . 
\end{aligned}\]
This completes the proof.

\begin{proof}[Proof of Proposition~\ref{prop:inter}]
Let $N_i = \ceil{\frac{n-i}{m}}$ and observe that 
$N-1 \le  N_i  \le N$  for $i\in \{0,\dots,m-1\}$, with $N=  \lfloor n/m \rfloor+1$. 
By Theorem~\ref{Rains}, for $\zeta\in\C$,
\[
F_{n,m}(\zeta) = \prod_{0\leq i < m}F_{N_i}\Big(\frac{\zeta}{\sqrt{m}}\Big) .
\]
Then, by (\ref{charfunction}), $D_N(z) = F_N(z) e^{|z|^2/4}$ is a Fredholm determinant and we obtain for $\zeta\in\C$,
\begin{equation} \label{decomp}
F_{n,m}(\zeta) = e^{-|\zeta|^2/4} \prod_{0\leq i < m} D_{N_i}\Big(\frac{\zeta}{\sqrt{m}}\Big) .
\end{equation}

By Plemelj's formula (cf.~the proof of Corollary \ref{distfrom1}),  if $\lvert z \rvert \leq 2\cst N$, 
\[\begin{aligned}
\Big| \prod_{0\leq i < m} D_{N_i}(z) -1 \Big|
& \le  \exp \Big( \sum_{0\leq i < m}  \sum_{j\in\N} \frac{| \tr (K(z)Q_{N_i})^j|}{j}\Big) -1 \\
& \le \exp \Big( \sum_{0\leq i < m}  | \tr (K(z)Q_{N_i})| \big( 1+\tfrac32 | \tr (K(z)Q_{N_i})|  \big) \Big) -1
\end{aligned}\]
Moreover, as in the proof of Lemma~\ref{traceclassnorm1}, we can bound
\[
| \tr (K(z)Q_{N_i})| \le  \| Q_{N_i}K(z)Q_{N_i} \|_{J_1}
\le \|Q_{N_i}H_\pm(w^\pm)\|_{J_2}^2 \le \|Q_{N-1}H_\pm(w^\pm)\|_{J_2}^2 
\]
since $N_i \ge N-1$. 
Then, using the bound \eqref{HSnorm}, we obtain for $|z| \le 2\cst N$,
\[
\big| \prod_{0\leq i < m} D_{N_i}(z) -1 \Big|
\le \exp \Big( \frac{\Cst_1 m}{\Gamma(N+1)^2} (\lvert z\rvert /2)^{2N} \Big) -1 .
\]
where $\Cst_1 = \frac{e^{2\cst^2}}{(1-\cst^2)^2}(1+\frac{3\Cst_0}{4} )$.
In particular, by convexity of $\exp$, we obtain the following bounds; for $n, m \in \N$,
\begin{equation} \label{ubmult1}
\Big| \prod_{0\leq i < m} D_{N_i}\big(\tfrac{\zeta}{\sqrt m}\big) -1 \Big| 
\le\big(e^{\Cst_2}-1\big) \frac{\lvert \zeta/2 \rvert^{2N} }{\Gamma(N+1)^2m^{N-1}} 
\qquad\text{if } |\zeta| \le  2 \cst N m^{\frac{1-1/N}{2}} 
\end{equation}
(here we used that $\Gamma(N+1)^{\frac 1N} \ge e^{-1} N$ by \eqref{Gamma} for $N\ge 2$ and 
$\cst^2e^{2\cst^2+2} =1$ so that $\Cst_2 = \frac{\cst^2 e^{2}}{(1-\cst^2)^2}(1+\frac{3\Cst_0}{4}) $) and similarly, 
\begin{equation} \label{ubmult2}
\Big| \prod_{0\leq i < m} D_{N_i}\big(\tfrac{\zeta}{\sqrt m}\big) \Big| 
\le \exp \big( \tfrac{\Cst_1  (\cst e)^{2(N-1)}}{2\pi N^3} |\zeta/2|^2 \big) \qquad\text{if } |\zeta| \le  2\cst N \sqrt{m} . 
\end{equation}

Then we split
\[
\big\| F_{n,m} - \widehat{\phi_\C} \big\|_{L^2}^2
\le  I_1+ I_2 +I_3 +I_4
\]
where 
\[\begin{aligned}
I_1 &= \int_{|\zeta| \le 2  \Lambda_{N,m}} \big|F_{n,m}(\zeta)- e^{-\frac{|\zeta|^2}{4}} \big|^2 \d^2\zeta \, ;  \qquad\qquad \Lambda_{N,m} =  \cst N m^{\frac{1-1/N}{2}} , \\ 
I_2 &= \int_{\Lambda_{N,m} \le |\zeta/2| \le \cst N\sqrt m}  
\big|F_{n,m}(\zeta)|^2  \d^2\zeta + \int_{|\zeta| \ge 2\Lambda_{N,m}}e^{-\frac{|\zeta|^2}{2}} \d^2 \zeta  \, ,\\
I_3 & = \int_{\Lambda \ge |\frac{\zeta}{N\sqrt{m}}| \ge 2\cst }  
\big|F_{n,m}(\zeta)|^2  \d^2\zeta 
\qquad\text{and}\qquad
 I_4  = \int_{ |\zeta| \ge \Lambda  N \sqrt m} \big|F_{n,m}(\zeta)|^2  \d^2\zeta  \, ,
\end{aligned}\]
with $\Lambda > 2\cst$. 
\medskip

First, using the bound \eqref{ubmult1}, we have
\[\begin{aligned}
I_1 & = \int_{|\zeta| \le \Lambda_{N,m}}  \Big| \prod_{0\leq i < m} D_{N_i}\big(\tfrac{\zeta}{\sqrt m}\big) -1 \Big|^2 e^{-\frac{|\zeta|^2}{2}}\d^2\zeta  \\
& \le \frac{\big(e^{\Cst_2}-1\big)^2}{\Gamma(N+1)^4} \frac{4\pi}{m^{2(N-1)}} \int_0^\infty  u^{2N}e^{-2u} \d u \\
& = 2\pi \big(e^{\Cst_2}-1\big)^2\frac{\Gamma(2N+1)2^{-2N}}{\Gamma(N+1)^4 m^{2N}} \\
&\le \frac{2\sqrt\pi \big(e^{\Cst_2}-1\big)^2}{\sqrt{N}\Gamma(N+1)^2m^{2(N-1)}}
\end{aligned}\]
where we used the inequality \eqref{duplication} at the last step and we verify numerically that 
$\frac{(e^{\Cst_2}-1)^2}{\sqrt\pi} \le 0.04$. 
  
Second, using the bound \eqref{ubmult2}, we obtain
\[\begin{aligned}
I_2& = \int_{|\zeta| \ge 2\Lambda_{N,m}}  \bigg(
\Big| \prod_{0\leq i < m} F_{N_i}\big(\tfrac{\zeta}{\sqrt m}\big) \Big|^2 \mathbf{1}\big\{ |\zeta/2| \le \cst N \sqrt m \big\} +1 \bigg) e^{-\frac{|\zeta|^2}{2}} \d^2 \zeta \\
&\le 2  \int_{|\zeta| \ge 2\Lambda_{N,m}}\exp\Big(-|\zeta/2|^2 \big(2-  \tfrac{\Cst_1  (\cst e)^{2(N-1)}}{N^3}  \big)\Big) \d^2 \zeta
\\
&\le 8\pi \int_{u  \ge\Lambda_{N,m}^{2}} e^{-\Cst_3 u} \d u 
 = \frac{8\pi}{\Cst_3} \exp\big( -\Cst_3  N^2 m^{1-1/N} \big) \, ; \qquad \qquad
\Cst_3(N) = 2- \frac{\Cst_1 (\cst e)^{2(N-1)}}{2\pi N^3}.
\end{aligned}\]
Note that $\Cst_3(N) \le \Cst_3(2) =\frac{1+\frac{3\Cst_0}{4} }{(1-\cst^2)^2} \le 1.2 $.

Third, 
By \eqref{ub10}, the constant in the bound of Proposition \ref{prop:TB1} can directly be improved;  for $n\in\N$ with $n\ge (N-1)$ and $z\in\C$,
\begin{equation*} 
| F_n(z) | \le 1.016 
 \exp\left(-  \frac{(N-1)^2}{4(2.05 + (N-1)^2/ |z|^2)} \right) .
\end{equation*}
This implies that for any $N\ge 2$, $n\ge N-1$ and $|z| \ge  2\cst N$, 
\[
| F_n(z) | \le 1.016 e^{-0.014 N^2}
\]
Hence,
\[\begin{aligned}
 \prod_{0\leq i < m}\big| F_{N_i}(z) \big|^2
&\le    (1.016)^{2m}  e^{-0.028 m N^2} . 
\end{aligned}\]

The point is that for $n,m\in\N$  and $|\zeta| \ge  2\cst N$,
\[
\big|F_{n,m}(\zeta)|^2   \le \exp(-\Cst_4 m N^2)
\]
where $\Cst_4 = 0.028 - \frac{2\log(1.016)}{N^2} \ge .02$. This bound implies that for any $\Lambda \ge 2\cst$,
\[
I_3 = \int_{\Lambda\ge |\frac{\zeta}{N \sqrt{m}}| \ge 2\cst }  
\big|F_{n,m}(\zeta)|^2  \d^2\zeta  \le 
\pi m N^2 \Lambda^2 \exp(-\Cst_4 m N^2) . 
\]

Finally, since $N_i\in\{N-1,N\}$ and $\sum_{i<m} N_i = n$, by Proposition~\ref{lem:Had1}, it holds for any $|\zeta| \ge 7N/4$,
\[
 \prod_{0\leq i < m}\big| F_{N_i}(\zeta) \big|^2
 \le \bigg| \frac{N}{\zeta} \bigg|^n
\]
Hence, for any $\Lambda \ge 7/4$,
\[
I_4 = \int_{|\zeta| \ge  \Lambda N\sqrt m} \big|F_{n,m}(\zeta)\big|^2 \d^2\zeta \le  2\pi m N^2   \int_{u\ge \Lambda} \frac{\d u}{u^{n-1}} 
=\frac{2\pi m N^2}{n \Lambda^n}.
\]

In order to optimize the contribution from $I_3$ and $I_4$, we must minimize the function 
\[
\Lambda \mapsto \tfrac12 \Lambda^2 \exp(-\Cst_4 n N)+ \tfrac 1n \Lambda^{-n} 
\] 
over all $\Lambda \ge 2$.
The (unique) minimizer is given by $\Lambda =\max\{2,e^{\Cst_4 N} \}$ depending on $N \ge 2$.
Accordingly, we obtain
\[\begin{aligned}
I_3+I_4  & \le \pi  \big(1+ \tfrac{1}{2n}\big) m N^2  \exp(-\Cst_4 (n-2)N) \, ,  & N\ge 35 \\
I_3+I_4 & \le   4 \pi  \big(1+ \tfrac{1}{2n}\big)  m N^2 \exp(-\Cst_4 n N) \, , &N\le 34 . 
\end{aligned}\]
Putting the contributions from $I_j$ for $j\in\{1,2,3,4\}$ yields the claimed bound. 
\end{proof}

\subsection{Orthogonal and symplectic groups: Proof of Theorem \ref{thm:interO}}

We let $p_{nm}^{ab}$ denote the density function of $\frac{1}{\sqrt{m}}(\tr \mathbf{O}^m-\E[\tr \mathbf{O}^m])$.

\begin{prop}\label{prop:interO}
Set $N=\floor{n/m}$. 
If $N\geq 2$, then
\begin{multline*}
\|p_{nm}^{ab}-\|_{L^2} \leq  \frac{\delta_1(m,N)}{m^{N-1}(2N)^{1/4}\sqrt{\Gamma(2N+1)}}\\
 + 4m^{1/4}\sqrt{N+1}\Big(\frac{Nm}{Nm-1}\Big)^{1/2} \exp \Big( -\beta(m,N)m^{1-\frac{1}{N}}\Big(1-\frac{1}{Nm}\Big)N^2\Big) 
\end{multline*}
where 
\begin{equation*}
\delta_1(m,N) = \frac{e^{\frac{2}{m}(1+\frac{1}{N})+\frac{1}{\sqrt{N}}}\sqrt{1-e^{-5/2}}}{(1-e^{-5/2}-1/(2\sqrt{\pi N}))^{3/2}}, \qquad \beta(m,N) = \frac{(1-\frac{1}{2N}-\frac{8\sqrt{2}}{\pi \sqrt{3(N-1/2)}N})^2}{2e^{5/2}m^\frac{1}{N}+1+\frac{8\sqrt{2}}{\sqrt{3(N-1)}}}.
\end{equation*}

\end{prop}

As a consequence, if $N$ and $m$ satisfy one of the following conditions:
\[ N \geq 5\ \mathrm{and}\ m\geq 66, \quad N \geq 4\ \mathrm{and} \ m\geq 129, \quad N \geq 2\ \mathrm{and}\ m\geq 10^4 \] 
then
\begin{equation}\label{l2powerortho}
\|p_{nm}^{ab}-\widehat{\phi_\R}\|_{L^2} \leq  \frac{\delta_1(m,N)+\tfrac{1}{mN}}{m^{N-1}(2N)^{1/4}\sqrt{\Gamma(2N+1)}}.
\end{equation}

\begin{remark}
We need stronger conditions on $N$ and $m$ here compared to the unitary case because of the larger bound we obtained for the intermediate regime (compare \eqref{ub2} with \eqref{intermediateO}). This in turn follows from the fact that the Fredholm determinants in Section \ref{sec:O} are not bounded by one, unlike those in Section \ref{sec:U}. 
\end{remark}

To prove Theorem \ref{thm:interO} we also need a concentration inequality for $\tr \mathbf{O}^m$, but this follows easily from the proof of that of $\tr \mathbf{O}$ and Rains' theorem. Indeed, in the case of $a=b=-1/2$ and $N=n/m \in \N$, we have, by (\ref{even_odd}),
\[ E_n^{--}[e^{\frac{\lambda}{\sqrt{m}} \tr \mathbf{O}^m}]  = F_N^{++}\Big(\frac{-\i\lambda}{\sqrt{m}}\Big)\Big(F_N^{+-}\Big(\frac{-\i\lambda}{\sqrt{m}}\Big)F_N^{-+}\Big(\frac{-\i\lambda}{\sqrt{m}}\Big)\Big)^{\frac{m-1}{2}} \]
where $\lambda>0$, so Proposition \ref{Fred} gives
\[ E_n^{--}[e^{\frac{\lambda}{\sqrt{m}} \tr \mathbf{O}^m}]  = e^{\lambda^2/2}\det (I+Q_NK_{\lambda/\sqrt{m}}^{--}Q_N)(\det (I+Q_NK_{\lambda/\sqrt{m}}^{-+}Q_N)\det (I+Q_NK_{\lambda/\sqrt{m}}^{+-}Q_N))^{\frac{m-1}{2}} \]
and (\ref{bound_det_concentration}) gives
\[ E_n^{--}[e^{\frac{\lambda}{\sqrt{m}} \tr \mathbf{O}^m}]  \leq \exp\Big(\frac{\lambda^2}{2} \Big(1 + \frac{1}{N^{2.5}e^{N/2}} \Big)\Big) \]
if $N\geq 2$ and $|\lambda| \leq 2e^{-3/2} \sqrt{m}N$. Thus, by Markov's inequality,
\begin{align*}
\mathbb{P}_n^{--}\big[\lvert \tfrac{1}{\sqrt{m}}\tr \mathbf{O}^m  \rvert \geq \tfrac{L}{2}]
&\leq 2 \exp\Big(-\frac{\lambda L}{2}+\frac{\lambda^2}{2} \Big(1 + \frac{1}{N^{2.5}e^{N/2}} \Big)\Big) \\
&\leq 2 \exp\Big(-\frac{L^2}{8(1 + N^{-2.5}e^{-N/2})}\Big)
\end{align*}
by choosing $\lambda = L/(2(1 + N^{-2.5}e^{-N/2}))$. This last bound can be shown to hold for any $a$ and $b$.
We can therefore insert it in
\[ \|p_{nm}^{ab}-\widehat{\phi_\R}\|_{L^1(\R)}
\leq \sqrt{L}\| p_{nm}^{ab}- \widehat{\phi_\R} \|_{L^2} + \mathbb{P}_n^{ab}\big[\lvert \tr \mathbf{O}^m - \E[\tr \mathbf{O}^m] \rvert \geq \tfrac{L}{2}\big] + \mathbb{P}[\lvert \boldsymbol{\gamma}_\R \rvert > \tfrac{L}{2}] \]
which, as usual, follows from the Cauchy-Schwarz inequality.
We also use (\ref{concentrationrealgaussian}) and our previous bound on the $L^2$-distance (\ref{l2powerortho}) to obtain
\begin{multline*}
\|p_{nm}^{ab}-\widehat{\phi_\R}\|_{L^1(\R)}
\leq \sqrt{L}\frac{\delta_1(m,N)+1/(mN)}{m^{N-1}(2N)^{1/4}\sqrt{\Gamma(2N+1)}} + \\
\Big(2 +\frac{4}{\sqrt{2\pi} L}\Big)\exp \Big(-\frac{L^2}{8(1+N^{-2.5}e^{-N/2})}\Big).
\end{multline*}
Setting $L = 2\sqrt{2} (1+N^{-2.5}e^{-N/2})^{1/2}(\log(\Gamma(2N+1))/2+\log(N)/4+(N-1)\log m)^{1/2}$ gives
\[ \|p_{nm}^{ab}-\widehat{\phi_\R}\|_{L^1(\R)}
\leq  \delta_2(m,N)\frac{(\log n)^{1/4}}{m^{N-1}\sqrt{\Gamma(2N+1)}} \]
where
\begin{multline*}
\delta_2(m,N) = [\sqrt{2}(1+N^{-2.5}e^{-N/2})^{1/4}(\delta_2(m,N)+\tfrac{1}{mN})+2((N+\tfrac{1}{4})\log N+(N-1)\log m)^{-1/4}\\+((N+\tfrac{1}{4})\log N+(N-1)\log m)^{-3/4}](1+\tfrac{1}{4N})^{1/4}.
\end{multline*}
for integer $N$ since $ L < 4e^{-3/2}(1 + N^{-2.5}e^{-N/2})\sqrt{m}N$ for any $N\geq 2$, $m \geq 18$. Since the bound is decreasing in $N$, we can simply replace $N$ with $\floor{N}$ in the case of non-integer $N$. A numerical approximation gives the bound of Theorem \ref{thm:interO}.

\begin{proof}[Proof of Proposition \ref{prop:interO}]
Let $F_{nm}^{ab}$ be the characteristic function of $\frac{1}{\sqrt{m}}(\tr \mathbf{O}^m-\E[\tr \mathbf{O}^m])$, $F_{nm}^{ab}(\xi) =\E_n^{ab}[e^{\i \frac{\xi}{\sqrt{m}}(\tr \mathbf{O}^m+a-b)}]$ if $m$ is odd and $ F_{nm}^{ab}(\xi) = \E_n^{ab}[e^{\i \frac{\xi}{\sqrt{m}}(\tr \mathbf{O}^m+a+b)}]$ if $m$ is even.
Consider $(a,b)=(1/2,1/2)$, $m$ odd; the remaining orthogonal/symplectic cases are treated in a similar way. If $n/m$ is an integer then (\ref{even_odd}) gives
\begin{align}\label{Fnmab0}
F_{nm}^{++}(\xi) = F_N^{++}\Big(\frac{\xi}{\sqrt{m}}\Big)\Big(F_N^{+-}\Big(\frac{\xi}{\sqrt{m}}\Big)F_N^{-+}\Big(\frac{\xi}{\sqrt{m}}\Big)\Big)^{\frac{m-1}{2}}
\end{align}
so by Proposition \ref{Fred},
\begin{multline}\label{Fnmab}
F_{nm}^{++}(\xi) = e^{-\xi^2/2}\det\Big(1+(-i)^{j+k}I_{j+k+2}\Big(\frac{2\xi}{\sqrt{m}}\Big)\Big)_{N\leq j,k} \\
\times \Big( \det(1-(-i)^{j+k+1}I_{j+k+1}\Big(\frac{2\xi}{\sqrt{m}}\Big)\Big)_{N\leq j,k}\det\Big(1+(-i)^{j+k+1}I_{j+k+1}\Big(\frac{2\xi}{\sqrt{m}}\Big)\Big)_{N\leq j,k}    \Big)^{\frac{m-1}{2}}
\end{multline}
where $I_n$ is the modified Bessel function of the first kind. We estimate $\|F_{nm}^{ab}-\widehat{\phi_R}\|_{L^2}$ in the usual way by dividing the real line into three parts: $\lvert\xi\rvert< 2N e^{-5/4}m^{(1-N^{-1})/2}$, $ 2N e^{-5/4}m^{(1-N^{-1})/2} < \lvert\xi\rvert < \Lambda $ and $\Lambda < \lvert\xi\rvert$, where $\Lambda\geq 2N e^{-5/4}m^{1/2-1/(2N)}$ can be chosen arbitrarily.
First, observe that if $\{x_i\}_{i=1}^m$ is a sequence of complex numbers satisfying $\max_{1\leq i\leq m} |x_i | \leq  \theta$, $\theta \geq 1$, then
\begin{equation}\label{distto1}
|1-\prod_{i=1}^m x_i| \leq \theta^{m-1} \sum_{i=1}^m |1-x_i|.
\end{equation}
This can be proved by induction. We apply this inequality to the determinants in (\ref{Fnmab}). By \eqref{freddistance}, for any $N\in\N$,
\begin{align*}
\lvert 1- \det(1+Q_NK_{\i \xi/\sqrt{m}}^{ab}Q_N)\rvert \leq \frac{\sqrt{1-e^{-5/2}}}{(1-e^{-5/2}-1/(2\sqrt{\pi N}))^{3/2}} \frac{e^{\xi^2/(m(2N+1))}}{m^N} \frac{|\xi|^{2N}}{\Gamma(2N+1)}
\end{align*}
if $|\xi| < 2\sqrt{m}N/e^{5/4}$. Next, to bound the determinants themselves, we use that for any trace-class kernel $K$, $\det(1-K)\leq e^{\|K\|_{J_1}}$ (Theorem VII.3.3 in \cite{GGK}). By (\ref{boundJ1}),
\[
\|Q_NK_{\i \xi/\sqrt{m}}^{ab}Q_N\|_{J_1} \leq \frac{e^2}{2\sqrt{e^2-1}(e-1)\sqrt{\pi N}} \Big(\frac{e^{5/4}\lvert\xi\rvert}{2\sqrt{m}N}\Big)^{2N} <\frac{1}{\sqrt{N}m}, \qquad \lvert\xi\rvert< \frac{2m^{\frac{1}{2}(1-\frac{1}{N})}N}{e^{5/4}}.
\]
The above upper bounds are decreasing in $N$. Therefore, if $n/m$ is not an integer we can simply replace it by $N = \floor{n/m}$ in the above expressions. Thus, by (\ref{Fnmab}) and (\ref{distto1}),
\begin{align*}
 \lvert F_{nm}^{++}(\xi) - e^{-\xi^2/2} \rvert \leq e^{-\xi^2/2} me^{\frac{m-1}{\sqrt{N}m}} \frac{\sqrt{1-e^{-5/2}}}{(1-e^{-5/2}-1/(2\sqrt{\pi N}))^{3/2}} \frac{e^{\xi^2/(m(2N+1))}}{m^{N}} \frac{|\xi|^{2N}}{\Gamma(2N+1)}   
\end{align*}
provided 
$$\lvert\xi\rvert< \frac{2m^{\frac{1}{2}(1-\frac{1}{N})}N}{e^{5/4}}.$$
We obtain 
\begin{align}\label{r1}
\Big( \int_{\lvert\xi\rvert < \frac{2Nm^{1/2-1/(2N)}}{e^{5/4}}}  \lvert & F_{nm}^{++}(\xi) - e^{-\frac{\xi^2}{2}} \rvert ^2 \d\xi \Big) ^{1/2} \nonumber \\
 &\leq \frac{e^{\frac{1}{\sqrt{N}}}\sqrt{1-e^{-5/2}}}{(1-e^{-5/2}-1/(2\sqrt{\pi N}))^{3/2}} \frac{m^{1-N}}{(2N)!} \Big( \int_{-\infty}^{\infty}  e^{-(1-\frac{1}{mN})\xi^2}\lvert\xi\rvert^{4N} \d\xi \Big) ^{1/2} \nonumber \\
 &\leq \frac{\delta_1(m,N)}{m^{N-1}(2N)^{1/4}\sqrt{\Gamma(2N+1)}}.
\end{align}
by \eqref{bd1} and the following two estimates.
Next, we use our previous results for the intermediate range of $\xi$. By (\ref{lastcase}), 
\[
|F_N^{++}(\frac{\xi}{\sqrt{m}})|  \leq \exp\Big(\frac{-\frac{1}{4m}(1+\frac{1}{2N}-\frac{8\sqrt{2}}{\pi\sqrt{3N}N})^2\xi^2}{2+\frac{2}{N}+\frac{\xi^2}{4mN^2}+\frac{2\sqrt{2}\xi^2}{\sqrt{3N}mN^2}}\Big), \qquad N\geq2.  
\]
So if $|\xi| > \frac{2Nm^{\frac{1}{2}(1-N^{-1})}}{e^{5/4}}$,
\begin{align*}
|F_N^{++}\big(\frac{\xi}{\sqrt{m}}\big)|  \leq \exp\Big(\frac{-(1+\frac{1}{2N}-\frac{8\sqrt{2}}{\pi \sqrt{3N}N})^2N^2}{2e^{5/2}(1+\frac{1}{N})m^\frac{1}{N}+1+\frac{8\sqrt{2}}{\sqrt{3N}}}\Big), \qquad N\geq2.
\end{align*}
The bound is decreasing in $N$, so combined with (\ref{Fnmab0}), it gives
\begin{align}\label{l1}
\int_{ \frac{2Nm^{1/2-1/(2N)}}{e^{5/4}} <\lvert\xi\rvert<\Lambda} \lvert F_{nm}^{++}(\xi)\rvert^2 \d\xi \leq 2\Lambda  \exp \Big(-2\beta(m,N)m^{1-\frac{1}{N}}N^2\Big)
\end{align}
for $N\geq 5$, $\Lambda\geq 2N e^{-5/4}m^{1/2-1/(2N)}$.

For the last regime, we use (\ref{Fnmab0}) together with
$$ \big|  F_N^{ab}\big(\frac{\xi}{\sqrt{m}}\big) \big|^2   \leq \Big(\frac{2\sqrt{m}N}{\lvert\xi\rvert}\Big)^N $$
valid for integer $N\geq 2$ and $|\xi| \geq  \frac{7}{8}\sqrt{m}N $, which follows from Proposition (\ref{lambda2}). The right-hand side is non-increasing in $N$ for $|\xi|\geq 2\sqrt{m}N$. Thus, if $\Lambda \geq  \max\{\frac{7}{8}\sqrt{m}(N+1), 2\sqrt{m}N\}$,
\begin{align}\label{l2}
\int_{\Lambda<\lvert\xi\rvert} \lvert F_{nm}^{++}(\xi)\rvert^2 \d\xi \leq \frac{ 2 (2\sqrt{m}N)^{Nm}}{(Nm-1)\Lambda^{Nm-1}}.
\end{align}
The sum of the upper bounds in (\ref{l1}) and (\ref{l2}) is minimized at $\Lambda=2\sqrt{m}N \exp (2\beta(m,N)N/m^{\frac{1}{N}})$, but (\ref{l2}) only holds for $\Lambda \geq  \max\{\frac{7}{8}\sqrt{m}(N+1), 2\sqrt{m}N\}$, so we set 
\[ \Lambda=2\sqrt{m}(N+1) \exp(2\beta(m,N)N/m^{\frac{1}{N}}).\]
This gives
\begin{align}\label{r2}
\int_{ \Lambda <\lvert\xi\rvert} \lvert F_{nm}^{++}(\xi)\rvert^2 \d\xi \leq 4\sqrt{m}(N+1)\frac{Nm}{Nm-1} \exp \Big( -2\beta(m,N)m^{1-\frac{1}{N}}N^2\Big(1-\frac{1}{Nm}\Big)\Big).
\end{align}
Finally, note that
\begin{align}\label{r3}
\int_{ \frac{2Nm^{1/2-1/(2N)}}{e^{5/4}} \leq |\xi|} e^{-\xi^2} \d\xi = \sqrt{\pi}\mathrm{erfc}\Big(\frac{2N m^{\frac{1}{2}(1-\frac{1}{N})}}{e^{5/4}}\Big) \leq \frac{e^{5/4}}{\sqrt{2\pi}Nm^{\frac{1}{2}-\frac{1}{2N}}}\exp \Big(-\frac{2N^2m^{1-\frac{1}{N}}}{e^{5/2}}\Big).
\end{align}
Equations (\ref{r1}), (\ref{r2}) and (\ref{r3}) give the final bound.
\end{proof}

\end{document}